\documentclass[english,11pt,leqno]{article}
\usepackage[latin1]{inputenc}

\usepackage{fullpage}

\usepackage{tipa}
\usepackage{dsfont}
\usepackage{tikz-cd}
\usepackage{float}
\usepackage{graphicx}

\usepackage{amsxtra}
\usepackage{amsmath}
\usepackage{amssymb}
\usepackage{amsfonts}
\usepackage[all,2cell]{xy}
\usepackage{mathrsfs}
\usepackage{amsthm}
\usepackage{enumitem}
  \usepackage{epigraph}
\usepackage{hyperref}
\usepackage{subfigure}
\usepackage{mathabx}
\usepackage{euscript}
\usepackage[bbgreekl]{mathbbol}
\usepackage{pgfplots}
\usepackage{mathtools}
 
\usetikzlibrary{patterns}
\usetikzlibrary{calc}

\usepackage[OT2,T1]{fontenc}
\DeclareSymbolFont{cyrletters}{OT2}{wncyr}{m}{n}
\DeclareMathSymbol{\Sha}{\mathalpha}{cyrletters}{"58}
 
\makeatletter
\hypersetup{
    unicode=true,           
    pdftoolbar=true,        
    pdfmenubar=true,        
    pdffitwindow=false,     
    pdfstartview={FitH},    
    pdftitle={\@title},     
    pdfauthor={\@author},   
    pdfsubject={},          
    pdfcreator={},          
    pdfproducer={},         
    pdfkeywords={},         
    pdfnewwindow=true,      
    colorlinks,             
    linkcolor=black,        
    citecolor=black,        
    filecolor=black,        
    urlcolor=black          
}
\makeatother

\newtheorem{thm}[equation]{Theorem}
\newtheorem{cor}[equation]{Corollary}
\newtheorem{lem}[equation]{Lemma}
\newtheorem{prop}[equation]{Proposition}

\newtheoremstyle{example}{\topsep}{\topsep}%
     {}
     {}
     {\bfseries}
     {.}
     {2pt}
     {\thmname{#1}\thmnumber{ #2}\thmnote{ #3}}

   \theoremstyle{example}
   
   \newtheorem{Defi}[equation]{Definition}
   \newtheorem{rem}[equation]{Remark}
   \newtheorem{rems}[equation]{Remarks}
   
   \newtheorem{exas}[equation]{Examples}
   \newtheorem{ex}[equation]{Example}

\numberwithin{equation}{section}

\renewcommand{\theparagraph}{\Alph{paragraph}.} 
\def\AAA{\mathbb{A}}

\def\CC{\mathbb{C}}
\def\DD{\mathbb{D}}
\def\EE{\mathbb{E}}
\def\FF{\mathbb{F}}
\def\GG{\mathbb{G}}
 
 \def\KK{{\mathbb{K}}}
 
\def\PP{\mathbb{P}}
\def\RR{\mathbb{R}}
\def\ZZ{\mathbb{Z}}

\def\ben{\mathfrak{b}}
\def\den{\mathfrak{d}}

\def\fen{\mathfrak{f}}
\def\gen{\mathfrak{g}}
\def\hen{\mathfrak{h}}

\def\len{\mathfrak{l}}

\def\pen{\mathfrak{p}}
\def\sen{\mathfrak{s}}

\def\uen{\mathfrak{u}}
\def\ven{\mathfrak{v}}

\def\Cen{\mathfrak{C}}
\def\Den{\mathfrak{D}}
\def\Een{\mathfrak{E}}

\def\Ien{\mathfrak{I}}

\def\Qen{\mathfrak{Q}}
\def\Ren{\mathfrak{R}}
\def\Sen{\mathfrak{S}}

\def\Ac{\mathcal{A}}
\def\Bc{\mathcal{B}}
\def\Cc{\mathcal{C}}

\def\Dc{\mathcal{D}}
\def\Ec{\mathcal{E}}
\def\Fc{\mathcal{F}}
\def\Gc{\mathcal{G}}

\def\Lc{\mathcal{L}}
\def\Mc{\mathcal{M}}
\def\Nc{\mathcal{N}}
\def\Hc{\mathcal{H}}
\def\Oc{\mathcal{O}}

\def\Qc{\mathcal{Q}}

\def\Sc{\mathcal{S}}
\def\Tc{\mathcal{T}}
\def\Uc{{\mathcal {U}}}
\def\Vc{\mathcal{V}}

\def\ba{{\mathbf{a}}}
\def\bb{{\mathbf{b}}}
\def\bc{{\mathbf{c}}}
\def\bd{{\mathbf{d}}}

\def\bm{{\mathbf{m}}}
\def\bn{{\mathbf{n}}}
\def\bp{{\mathbf{p}}}
 
 \def\bx{{\mathbf{x}}}

 \def\be{\begin{equation}}
 \def\bem{\begin{matrix}}
 \def\Bic{{\on{Bic}}}
 \def\bpm{\begin{pmatrix}}
 \def\Br{{\on{Br}}}
 \def\Bun{{\on{Bun}}}

 \def\codim{{\on{codim}}}
 \def\Coh{{\on{Coh}}}
  
 \def\Cone{{\on{Cone}}}
 \def\constr{{\on{constr}}}

 \def\coor{{\on{coor}}}

 \def\del{{\partial}}
  
  \def\Dsim{{\Delta_\simp}}

  \def\ee{\end{equation}}
 
 \def\enm{\end{matrix}}
 \def\epm{\end{pmatrix}}
 \def\eps{{\varepsilon}}

 \def\Fun{{\on{Fun}}}

 \def\gr{{\on{gr}}}

 \def\hCW{{W\backslash \hen_\CC}}

   \def\hRW{{W\backslash\hen_\RR}}
   \def\hra{\hookrightarrow}
  \def\Hor{{\on{Hor}}}
  \def\hW{{W\backslash \hen}}

 \def\Id{\operatorname{Id}\nolimits}
 \def\Ind{{\on{Ind}}}

 \newcommand{\isom}{\overset{\sim}{=}}
  \def\Im{{\on{Im}}}

 \def\k{\mathbf k}
  \def\Ker{\operatorname{Ker}\nolimits}

  \def\Lie{{\on{Lie}}}
  \def\Lin{{\on{Lin}}}
 \def\lla{\longleftarrow}
   \def\lra{\longrightarrow}

 \def\MBS{{\on{MBS}}}

\def\on{\operatorname}
\def\ol{\overline}
\def\op{{\on{op}}}

\def\orr{{\on{or}}}

 \def\Perv{{\on{Perv}}}
 \def\phi{{\varphi}}
 \def\pt{{\on{pt}}}

\def\Re{{\on{Re}}}
\def\reg{{\on{reg}}}
 \def\Rep{\on{{Rep}}}
  \def\Res{\on{{Res}}}

  \def\sgn{{\on{sgn}}}
  \def\sh{{\on{sh}}}
 \def\Sh{{\on{Sh}}}

 \def\simp{{\on{sim}}}

 \def\sup{{\on{Sup}}}
 \def\Supp{{\on{Supp}}}
   \def\Sym{{\on{Sym}}}
   
\def\teq{{\, \trianglelefteq\, }}

\def\ul{\underline}

\def\Vect{\on{Vect}}
\def\Ver{{\on{Ver}}}

\def\wt{\widetilde}

\def\1{{\mathbf{1}}}
\def\2{{\mathbf{2}}}
\def\(({(\hskip -1mm (}
\def\)){)\hskip -1mm )}
\def\-{{\setminus}}
 \def\be{\begin{equation}}
\def\ee{\end{equation}}
\def\ed{\end{document}}

\title{ Parabolic induction and   perverse sheaves on $ \hW$  }
\author{ Mikhail Kapranov, Vadim Schechtman
 }
 
\begin{document}

\maketitle

\thanks{ }

 \setlength{\epigraphwidth}{0.6\textwidth}
 \epigraph{ Alle Gestalten sind \"ahnlich und keine gleichet den andern,  \\
 und so deutet das Chor auf ein geheimes Gesetz\ldots}
        { Goethe, {\em Die Metamorphose der Pflanzen }}

\begin{abstract}
 For a complex reductive Lie group $G$ with Lie algebra $\gen$, Cartan subalgebra $\hen\subset \gen$
 and Weyl group $W$, we describe the category of perverse sheaves on $\hW$ smooth
 w.r.t the natural stratification. The answer is given in terms of mixed Bruhat sheaves, which
 are  certain mixed sheaf-cosheaf data on cells of a natural cell decomposition of $\hW$.
 Using the parabolic Bruhat decomposition, we relate mixed Bruhat sheaves with   properties
 of various procedures of parabolic induction and restriction that
  connect different Levi subgroups in $G$.
 \end{abstract}
 
\tableofcontents

\pagebreak

\addtocounter{section}{-1}

\section{Introduction}

\paragraph {}For a complex reductive Lie algebra $\gen$,  the quotient $\hW$ of the Cartan subalgebra
by the Weyl group is isomorphic to an affine space and carries a natural stratification $\Sc^{(0)}$.
For example, for $\gen=\gen\len_n$ we get the space of monic polynomials  
\[
f(x) = x^n+a_1x^{n-1}+\cdots + a_n
\]
 and $\Sc^{(0)}$ is given by the singularities of the
discriminantal hypersurface $\Delta(f)=0$. 

\vskip .2cm

The first goal of this paper is to give an elementary (quiver-type)  description
of $\Perv(\hW)$, the category of perverse sheaves on $\hW$ smooth with respect to $\Sc^{(0)}$. 
Our main result, Theorem \ref{thm:main},  identifies  $\Perv(\hW)$, with the category of
objects of mixed functoriality with respect  to a natural cell decomposition $\{U_\bm\}_{\bm\in\Xi}$
of $\hW$ refining $\Sc^{(0)}$, objects which we call {\em mixed Bruhat sheaves}. 
A mixed Bruhat sheaf $E$ consists of vector spaces $E(\bm)$, one for each cell $U_\bm$
 and behaves like a cellular sheaf   with respect to one part of cell inclusions
 $\ol U_\bm \supset U_\bn$.  It behaves like a cellular cosheaf  
 with respect to  another part, see Definition \ref{def:MBS}. 
 
 \vskip .2cm 
 
 This mixed nature harmonizes well with the  self-dual, intermediate position of perverse sheaves
 themselves, half-way between the abelian categories of 
 sheaves and of cosheaves (i.e.,  objects of the derived category
 Verdier dual to sheaves).  At the more geometric level, this corresponds to
 the position of intersection homology  as half-way between cohomology and homology. 

\vskip .2cm

The indexing set $\Xi$ for our cell decomposition (and for mixed Bruhat sheaves) is the
 2-sided Coxeter complex of Petersen \cite{petersen}. For $\gen=\gen\len_n$
this is the set of contingency matrices of content $n$ studied in statistics \cite{diaconis, {KS-contingency}}.

\paragraph{} Our second and broader goal is to relate $\Perv(\hW)$ to a classical subject
of representation theory which well predates perverse sheaves: the ``algebra of parabolic
induction''.  By this we mean the entire package of results related to principal series 
(parabolically-induced)
representations of reductive groups (in all contexts: finite field, real, p-adic, adelic, automorphic),
their intertwiners,  Eisenstein series, constant terms of automorphic forms and other
procedures that pass from one Levi subgroup to another. 
 
 \vskip .2cm

The rules of this ``algebra'' 
(which underlies  the philosophy of cusp forms of Gelfand and Harish-Chandra)
are familiar to all practitioners of representation theory, but  it 
was somehow considered {\em sui generis}, their interpretation in terms of something else
being unclear or unlooked for.  For groups $\on{GL}_n$, one  can interpret some of the rules
 in terms of braided Hopf algebras
(via the concept of the Hall algebra) \cite{k-eis-series, {k-schiffman-vass}} and braided monoidal categories \cite{joyal-street, soibelman}. For a more general reductive group $G$ with Lie algebra $\gen$,
this is not possible.

\vskip .2cm

Our observation 
is that perverse sheaves  (and their categorical analogs, perverse schobers
\cite{KS-schobers})  on $\hW$  provide a conceptual encoding of this peculiar algebra,
giving it a name, so to speak. 
We illustrate this on two simplest examples in \S \ref{sec:Fq} and sketch some other examples
in \S \ref{sec:further}. A germ of this connection can be seen in  the fact that the  principal
series intertwiners form a representation of the braid group $\Br_\gen=\pi_1(\hW^\reg)$
and so give a local system on the generic stratum $\hW^\reg\subset\hW$. 
The examples we consider indicate  that the correspondence between
perverse sheaves on $\hW$ and the algebra of parabolic induction
 should hold in many different contexts. 
 
\vskip .2cm

Using $\on{GL}_n$ as a  point of departure, our approach can be seen as importing,
into the general  theory of representations and automorphic forms, 
the new 2-dimensional point of view on Hopf algebras coming from their relation
to $E_2$-algebras in the work of Lurie. 

\paragraph{} The reason for the relation between $\Perv(\hW)$ and parabolic induction
comes from the elementary but remarkable  matching between elements  $\bm\in \Xi$,
i.e., cells of $U_\bm \subset \hW$ and {\em Bruhat orbits}, i.e., $G$-orbits $O_\bm\subset F_I\times F_J$ in 
the pairwise products of all possible flag varieties for $G$. This matching is   just
the  parabolic Bruhat decomposition; the remarkable fact is that some topological  relations among the cells  $U_\bm$
have, as  counterparts,  algebro-geometric relations among the orbits $O_\bm$. 
For example, the property for a cell inclusion $\ol U_\bm\supset U_\bn$ is
{\em anodyne}  (i.e.,  that both cells lie in the same stratum of $\Sc^{(0)}$) corresponds to
the property that  the projection of orbits $p_{\bm,\bn}: O_\bm\to O_\bn$ 
is a fiber bundle of 
 affine spaces and so for example, gives an isomorphism in the
category of Voevodsky motives \cite{beil-volog}. 

\vskip .2cm

The data appearing in the theory of parabolic induction are usually labelled by
the standard Levis, i.e., by subsets $I\subset\Dsim$ of simple roots. This gives a 
{\em bicube},  i.e., a diagram of $2^{|\Dsim|}$ vector spaces or categories related
by maps or functors back and forth along the inclusions $I\subset J$ (e.g., induction
versus restriction,
Eisenstein series versus constant term), see \S \ref{sec:perv=mbs}\ref{par:bicube}
In  \S \ref{sec:Fq}-\ref{sec:further} we extend some  of these bicube
diagrams to   arbitrary $\bm\in\Xi$, i.e., to  an arbitrary Bruhat orbit $O_\bm$.
Informally,  such a larger diagram has the parabolic intertwiners already ``pre-installed'',
since among the orbits we find the correspondences used to define the intertwiners. 

\paragraph{} The organization of the paper is as follows. In Section \ref{sec:2-cox} we recall
the $2$-sided Coxeter complex $\Xi$
and  introduce the cell decomposition  of $\hW$ into cells $U_\bm$
labelled by $\Xi$. The definition (and thus the whole approach of the paper)
 involves separating the real and imaginary parts of a point of $\hen= \hen_\RR\oplus i\hen_\RR$. 
Thus, for $\gen=\gen\len_n$, a   cell in $\hW=\Sym^n(\CC)$ consists of polynomials
whose zeroes follow a given pattern of coincidences 
among their real and imaginary parts given by a contingency matrix. 

\vskip .2cm

In Section \ref{sec:perv=mbs} we define mixed Bruhat sheaves and formulate the main result,
Theorem   \ref{thm:main}. We also explain the relation of mixed Bruhat sheaves with bicubes
 and work out  the examples of $\sen\len_2$ and $\sen\len_3$. 
 
 \vskip .2cm
 
 Sections \ref{sec:cous-bruhat} and \ref{sec:cous-perv} are devoted to the proof of Theorem
 \ref{thm:main}. The proof is based on the  technique of  {\em Cousin complexes}, used
 in different forms in \cite{KS-hyp-arr, KS-shuffle,  {KS-hyp-arr-II}}. They are certain explicit
 complexes of sheaves whose terms are constructible with respect to an intermediate
 real stratification $\Sc^{(1)}$ of $\hW$ (in fact, also a cell decomposition),  
 whose strata we call {\em Fox-Neuwirth-Fuchs cells}. 
 As in the classical cell decompositions of configuration spaces \cite{fox-neuwirth, fuchs},
 the definition of these cells involves choosing
  a preference of the real parts over the imaginary parts. 
 At the same time, these complexes
 represent (i.e., are isomorphic to) perverse sheaves from $\Perv(\hW)$, so their
 cohomology sheaves are $\Sc^{(0)}$-constructible. Our proof of this cohomological
 $\Sc^{(0)}$-constructibility is based on the remarkable property
 \[
 \Sc^{(1)} \vee \tau\Sc^{(1)} = \Sc^{(0)}. 
 \]
Here $\tau\Sc^{(1)}$ is a stratification similar to $\Sc^{(1)}$ but with the roles of the real and 
imaginary parts interchanged, and the statement means that $\Sc^{(0)}$ is the smallest
stratification of which  $\Sc^{(1)}$ and $\tau\Sc^{(1)}$ are both refinements. 

\vskip .2cm

In fact in Sections \ref{sec:2-cox} to \ref{sec:cous-perv} we are working in a more general situation when $W$ is an arbitrary 
finite reflection group, non necessarily the Weyl group of a root system.

\vskip .2cm

In Sections \ref{sec:GBO} and \ref{sec:BO=MBC},
we discuss the geometry of the Bruhat orbits $O_\bm$
and its relation to the properties of the corresponding labels $\bm$, which are themselves
  $W$-orbits in coset spaces of the form $W/W_I\times W/W_J$. 
   In particular, we establish
the ``$\AAA^1$-equivalence''   property (Proposition \ref{prop:anodyne-affine}) of the orbit projection corresponding to an anodyne
inclusion $\bm\geq\bn$.

\vskip .2cm

In Section \ref{sec:Fq} we consider the simplest example of a motivic Bruhat sheaf
coming from Bruhat orbits: the sheaf of appropriate spaces of functions on $\FF_q$-points. 
We also work out an even easier ``$\FF_1$-version'', when we consider functions
on the $W$-orbits  $\bm$ themselves. The resulting perverse sheaves
are then identified in terms of  representations of the Hecke algebras and Coxeter groups. 

\vskip .2cm 

The concluding Section \ref{sec:further} sketches some further constructions in the same
spirit which we plan to develop in subsequent papers. In particular, we discuss
a natural categorical generalization of mixed Bruhat sheaves. 

\vskip .2cm

Finally, the Appendix collects  notations and conventions related to constructible sheaves and 
stratifications that are used in the main body of the paper.

\paragraph{} It would be interesting to understand the relation between $\Perv(\hW)$ and
parabolic induction in a more direct, intrinsic way. 

\vskip .2cm

Geometrically, parabolic subgroups in $G$
correspond to  various ways of approaching the infinity either in $G$ itself, or in its
arithmetic quotients. 
So one can think of    realizing the cell decomposition  $ \{U_\bm\}$  of $\hW$
as some combinatorial
complex describing regions at infinity in $G$ or in a related space. The closest picture of this kind known to us is the wonderful compatification
 $\ol G\supset G$, see \cite{deconcini, bernstein-bezr}. 
This is a smooth projective $(G\times G)$-variety with   $(G\times G)$-orbits $X_I$ labelled by  $I\subset\Dsim$.
For $G$ of adjoint type, $X_I$ fibers over $F_I\times F_I$ with the fiber being the adjoint quotient
of the corresponding Levi.  So considering the action of the diagonal $G\subset G\times G$ on $\ol G$
(which extends the action of $G$ on itself by conjugation) does  lead to the appearance of Bruhat
orbits but only those in  $F_I\times F_I$, not those in arbitrary $F_I\times F_J$. 

\vskip .2cm

In a  somewhat different direction,  
it seems interesting to understand the relation of  the cell
decomposition $\{U_\bm\}$ 
with the characteristic map
\[
\chi: \gen\lra \hW. 
\]
In particular, the topology of the regions $\chi^{-1}(U_\bm)\subset\gen$
seems worth studying: 
for example,    the way
they approach the nilpotent cone $\Nc=\chi^{-1}(0)$ in the case
where $\gen$ is semisimple.

\paragraph{} 
Our latest interest in  these questions was triggered by discussions with R. Bezrukavnikov.
We would also like to thank A. Beilinson and M. Finkelberg for useful correspondence and J. Tao for sending us his very interesting preprints
which develop
some topics of the present article. 
We are grateful to the anonymous referees for their comments that have led to several improvements
in  the paper. 
The research of M.K. was supported by World Premier International Research Center Initiative (WPI Initiative), 
 MEXT, Japan.


\section{The 2-sided Coxeter complex as a cell decomposition of $\hW$. }\label{sec:2-cox}

\paragraph{Notation.} We consider the classical context of finite reflection groups,
cf. \cite{bourbaki}, Ch. IV, \S 1. More precisely, we introduce the following notation;

\vskip .2cm

$\hen_\RR$ is a finite dimensional real  vector space equipped with a Euclidean metric. 

\vskip .2cm

$W\subset O( \hen_\RR)$ is a finite reflection group. In particular,  $W$ is a finite Coxeter group. 

\vskip .2cm

$\Hc$ is  the  collection of  reflection hyperplanes for $W$. Thus $\Hc$ is an arrangement of
linear hyperplanes in $\hen_\RR$ invariant under the action of $W$. For a reflection $s\in W$ we denote
by $H_s\in\Hc$ the corresponding reflection hyperplane.

\vskip .2cm

$\Cc$ is the Coxeter complex, i.e., the decomposition of $\hen_\RR$ into
faces of $\Hc$, see \cite {KS-hyp-arr}  \S 2A. We think of $\Cc$ as both
a geometric decomposition of $\hen_\RR$ and as a poset $(\Cc,\leq)$ of
faces ordered by inclusion of closures.  As usual, the notation $a<b$  means that
$a\leq b$ and
$a\neq b$. 

\vskip .2cm  

We choose a {\em dominant chamber} $C^+\in \Cc$, a connected component in the complement
$$
\hen_\RR \setminus\bigcup_{H\in\Hc} H.
$$
As is well known, the reflections with respect to walls of $C^+$ form a system of generators $S\subset W$, 
and $(W, S)$ is a Coxeter system, cf. \cite{bourbaki}, Ch. IV, \S 1.

\vskip .2cm

$C^+$ is an open simplicial cone with a face $C^+_I$ for each $I\subset S$.  
The notation is normalized so that $C^+_I$ is open in the subspace $H_I=\bigcap_{s\in I} H_s$. 
 By $\ol C^+$ we denote the closure of $C^+$, i.e., the union of all the $C^+_I$. 
 The $W$-action on $\Cc$ induces an identification
\[
\bigsqcup_{I\subset S} W/W_I  \,\buildrel \simeq \over\lra  \, \Cc,
\quad wW_I \mapsto w(C^+_I). 
\]
 We will write $\bc, \bd$, etc. for cosets $wW_I$ and $A_\bc, A_\bd $ etc. 
 for the corresponding faces, i.e., elements of $\Cc$. Thus we have two ways  of
 viewing elements of $\Cc$: one as cosets, the other as
  faces of the cell decomposition
 of $\hen_\RR$.  Both ways will be useful depending on the situation. 
 
 \vskip .2cm
 
  $\hen := \hen_\RR \otimes_\RR\CC $ is  the complexification of $\hen_\RR$. We denote by
$\Hc_\CC = \{H_\CC := H\otimes_\RR\CC,\ H\in \Hc\}$:   the complexified arrangement
in  $\hen$.

\vskip .2cm

$p: \hen\lra \hW$ is  the canonical projection of $\hen$ to its quotient  by $W$.   
By Chevalley's theorem, $\hW$ is
an algebraic  variety isomorphic to an affine space.  As is well known,  $\hW$ is in fact defined over $\ZZ$. 

\vskip .2cm

$K=\hRW$ is a closed curvilinear cone in the real affine space $(\hW)(\RR)$,
for which $p: \ol C^+\to K$ is a homeomorphism. We denote the faces of $K$ by
 $K_I=p(C^+_I)\subset K$,
$I\subset S$, 
. 

\vskip .2cm

$\Sc^{(0)}_\Hc$ is  the stratification of $\hen$ into generic parts of the flats of $\Hc$,
see  \cite[\S 2D] {KS-hyp-arr} . 
This stratification is $W$-invariant and so induces a complex stratification of $\hW$
which we denote $\Sc^{(0)}$. Thus  each stratum of $\Sc^{(0)}$ is the image
of the generic part of a flat of $\Hc$.

\paragraph{The two-sided Coxeter complex.} 
We further write:

\vskip .2cm

$\Xi = W\backslash (\Cc\times\Cc)$ for the two-sided Coxeter complex of Petersen
  \cite{petersen}. We view it as a poset of $W$-orbits, 
    with the order on $\Xi$ induced
  by the product order on $\Cc\times\Cc$. Thus 
   \[
 \Xi \,=\, \bigsqcup_{I,J\subset S} \Xi(I,J), \quad \text{where} \quad
 \Xi(I,J) \,:=\, W\backslash ((W/W_I)\times W/W_J)). 
 \]
 The sets $\Xi(I,J)$ are
connected by the {\em horizontal} and {\em vertical contraction maps}
\be\label{eq:contr-maps}
\begin{gathered}
\phi'_{(I_1, I_2|J)}: \Xi(I_1, J) \lra \Xi(I_2, J), \quad I_1\subset I_2,
\\
\phi''_{(I|J_1, J_2)}:  \Xi(I, J_1) \lra\Xi(I, J_2), \quad J_1\subset J_2.
\end{gathered} 
\ee
They are induced by the natural projections $W/W_{I_1}\to W/W_{I_2}$
and $W/W_{J_1}\to W/W_{J_2}$ respectively. The two types of maps 
commute with each other: For $I_1\subset I_2$ and $J_1\subset J_2$
the diagram below commutes: 
\be\label{eq:leq'-leq''-proj}
\xymatrix{
\Xi(I_1, J_1) \ar[rr]^{\phi'_{(I_1, I_2|J_1)}} 
\ar[d]_{\phi''_{(I_1|J_1, J_2)}}
&& \Xi(I_2, J_1)
\ar[d]^{\phi''_{(I_2|J_1, J_2)}}
\\
\Xi(I_1, J_2) \ar[rr]_{\phi'_{(I_1, I_2|J_2)}} && \Xi(I_2, J_2). 
}
\ee
 Indeed, this diagram   is obtained from the $W$-equivariant commutative square
 \be\label{eq:cart-diagr-W}
 \xymatrix{
 (W/W_{I_1}) \times (W/W_{J_1})\ar[d]_{\pi''}  \ar[r]^{\pi'}& (W/W_{I_2})\times (W/W_{J_1})
 \ar[d]^{\pi''} 
 \\
 (W/W_{I_1})\times (W/W_{J_2}) \ar[r]^{\pi'}& (W/W_{I_2})\times (W/W_{J_2}) 
 }
 \ee
  by taking
 $W$-quotients.  

 A typical element of  $\Xi(I,J)$ will be denoted $\bm = W(\bc, \bd)$,
where $\bc\in W/W_I$ and $\bd\in W/W_J$.
 We write $\bm \geq'\bn$, if $\bn$ is obtained from $\bm$ by a horizontal contraction,   and $\bm\geq'' \bn$, if  $\bn$ is obtained from $\bm$  by a vertical contraction.
 The following is obvious from the definition of the order on $\Xi$:
 
 \begin{prop}
 For $\bm, \bn\in \Xi$ the following are equivalent:
 \begin{itemize}
 \item[(i)] $\bm\geq \bn$.
 
 \item[(ii)] There exists a (unique) $\bm'$ such that $\bm\geq' \bm'\geq''\bn$.
 
 \item[(iii)] There exists a (unique) $\bn'$ such that $\bm\geq'' \bn'\geq'\bn$.   \qed
 \end{itemize}
 \end{prop}
 
 In particular, for $\bm\in\Xi(I_1, J_1)$ and $ \bn\in\Xi(I_2, J_2)$ the inequality
 $\bm\geq\bn$ implies that $I_1\subset I_2$ and $J_1\subset J_2$, and the arrows in
 \eqref{eq:cart-diagr-W} define a (surjective) $W$-invariant map
 \be\label{eq:pi-m-n}
 \pi_{\bm,\bn}: \bm\lra\bn, 
 \ee
 where we regard $\bm$ and $\bn$ as $W$-orbits in the corresponding objects
 of the diagram  \eqref{eq:cart-diagr-W}. These maps are transitive, i.e.,
 define a covariant functor from the category formed by the poset $(\Xi, \geq)$ 
 to the category of sets with $W$-action. 
 
 \paragraph{Mixed supremum in $\Xi$.} 
 Let $\bm', \bn\in\Xi$. Their {\em mixed supremum} is the subset
 \be\label{eq:sup}
 \sup(\bm', \bn) \,=\, \bigl\{ \bm\in\Xi \mid \, \bm' \leq'' \bm \geq' \bn\bigr\}. 
 \ee
 For this set to be nonempty, it is necessary that there be $I_1\subset I_2$ and
 $J_1\subset J_2$ such that 
  $ \bm'\in \Xi(I_1, J_2)$ and  $\bn\in \Xi(I_2, J_1)$, 
 in which case
 \[
  \sup(\bm', \bn)\,=\, \bigl\{ \bm\in\Xi(I_1, J_1) \bigl| \,\, 
  \phi''_{(I_1|J_1, J_2)}(\bm) = \bm', \,\, \phi'_{(I_1, I_2|J_1)} = \bn\bigr\}. 
 \]
 and furthermore, that 
 $
 \phi'_{(I_1, I_2|J_2)}(\bm') \,=\, \phi''_{(I_2|J_1, J_2)}(\bn).
 $
 Denoting this common value by $\bn'$, we have $\bm'\geq' \bn' \leq''\bn$.

 The following is then straightforward. 
 \begin{prop}\label{prop:sup=FP}
 $\sup(\bm',\bn)$ is the set of $W$-orbits in the fiber product 
  $
 \bm'\times_{\bn'}\bn  
 $,
 where we consider $\bm', \bn',\bn$ as  $W$-orbits
  in the corresponding objects of the diagram
  \eqref{eq:cart-diagr-W}. \qed
 \end{prop}

 \paragraph{The cell decomposition of $\hW$ and anodyne inequalities.}
 We further write:
 
 \vskip .2cm
 
 $\Sc^{(2)}_\Hc$ for the cell decomposition of $\hen = i\hen_\RR\oplus \hen_\RR$
 into the   {\em product cells} $iC+D$ where $C,D\in\Cc$ are faces of the arrangement $\Hc$. 
 Recall that each face of $\Hc$ has the form $C=A_\bc$ for some  $\bc\in W/W_I$ and 
  $I\subset S$. 

 \vskip .2cm
 
 $\Sc^{(2)} = p(\Sc^{(2)}_\Hc)$ for the decomposition of $\hW$ into the images
 \[
 U_\bm = p(iA_\bc+A_\bd)
 \]
 for $\bm = W(\bc, \bd)\in\Xi$. 
 
 \begin{prop}\label{prop:p-iA+B-homeo}
 The restriction of $p$ to the closure of each $iA_\bc+A_\bd$ is a homeomorphism
 onto its image. Therefore $\Sc^{(2)}$ is a quasi-regular cell decomposition of
 $\hW$ into the cells $U_\bm$, $\bm\in\Xi$, refining the complex stratification
 $\Sc^{(0)}$. 
 \end{prop}
 
 Thus we have three ways of viewing elements of $\Xi$: first, as left $W$-cosets in 
 various
 $(W/W_I)\times
 (W/W_j)$, second, as left $W$-cosets in the product cell decomposition $\Sc^{(2)}_\Hc$
 of $\hen_\CC$ and third, as cells of the cell decomposition $\Sc^{(2)}$ of $\hW$. 
 We will use all three of these ways, depending on the situation.
 
 \vskip .2cm

 \noindent{\sl Proof of Proposition \ref{prop:p-iA+B-homeo}:}  
 We need to show the following: if $z,z'\in \ol U_\bm =  i\ol A_\bc + \ol A_\bd$
 and $p(z)=p(z')$, meaning  $z'=wz$ for some $w\in W$, then $z=z'$. 
 A
  similar statement for  the $W$-action on $\hen_\RR$
 is standard: if $C\in\Cc$ and $x,x'\in \ol C$ are such that $x'=wx$ for some
 $w\in W$, then $x=x'$. This is because we can reduce to $C=C^+$ and 
 $p$ maps the closure
 of any face bijectively onto the closure of a face of $K$. 
 
 Writing $z=ix+y$, $z'=ix'+y'$ with $x,y,x',y'\in\hen_\RR$,
 the condition $z'=wz$ means
 $x'=wx$ and $y'=wy$.  So  by the above,  $x'=x$ and  $y'=y$, meaning $z'=z$. \qed
 
 \vskip .2cm 
 
 We note that $\bm\geq\bn$ if and only if $U_\bm\supset\ol U_\bn$. 
 
 \begin{Defi}
  An inequality $\bm\geq\bn$ (or $\bm\geq' \bn$ or $\bm\geq'' \bn$)
  is called {\em anodyne} if $U_\bm$ and $U_\bn$ lie in the same stratum of
  $\Sc^{(0)}$. 
 \end{Defi}
 
 \begin{prop}\label{prop:ano=bij}
 Let $\bm\geq\bn$. The following are equivalent:
 \begin{itemize}
 \item[(i)] The  inequality $\bm\geq\bn$ is anodyne.
 
 \item[(ii)]   $|\bm| = |\bn|$ when we regard
  $\bm, \bn$
 as subsets in the poset $\Cc\times\Cc$.

 \item [(iii)] The map $\pi_{\bm,\bn}: \bm\to\bn$ is a bijection. 
 \end{itemize}

 \end{prop}
 
 \noindent{\sl Proof:}  Let us prove (i)$\iff $(ii). 
 
 For $\bc\in \Cc$ let $W^\bc\subset W$ be the stabilizer of
 $\bc$ in $W$.  Let $\bm= W(\bc, \bd)$.  Then $|\bm| = |W|/|W^\bc\cap W^\bd|$. 
 Now, $W^\bc\cap W^\bd$ is the stabilizer of the cell $iA_\bc + iA_\bd\subset \hen$
 or, what is the same by Proposition \ref{prop:p-iA+B-homeo}, the stabilizer
 of any point in this cell. 
 
 Since $\bm\geq\bn$, we can represent $\bn=W(\ba,\bb)$
 with $\bc\geq\ba$ and $\bd\geq\bb$, i.e., $A_\bc\supset \ol A_\ba$
 and $A_\bd\supset \ol A_\bb$. Let $S_\bm$ and $S_\bn$
 be the strata of $\Sc^{(0)}$ containing $U_\bm$ and $U_\bn$, so that
 $S_\bm\supset \ol S_\bn$.  That is, $S_\bm= p(L_{\bc, \bd})$
 and $S_\bn = p(L_{\ba, \bb})$, where $L_{\bc, \bd}$ is the stratum
 of $\Sc^{(0)}_\Hc$  containing $iA_\bc+A_\bd$, and similarly with $L_{\ba, \bb}$. 
 
 Note that $L_{\bc, \bd}\supset \ol L_{\ba, \bb}$ and, moreover, 
 $L_{\bc, \bd} = L_{\ba, \bb}$ if and only if $S_\bm=S_\bn$. 
  Note further that the points in $ L_{\bc, \bd}$ all have the stabilizer 
 $W^\bc\cap W^\bd$ while points in $L_{\ba, \bb}$ all  have  same
 stabilizer $W^\ba\cap W^\bb$. 
 
 Now,   generic points of a (strictly) smaller complex
 flat of the root arrangement, have (strictly) larger stabilizer in $W$. 
 Therefore
 we have $|W^\ba\cap W^\bb|\geq |W^\bc\cap W^\bd|$ with equality if and only if
 that $ L_{\bc, \bd}=  L_{\ba, \bb}$, i.e., $S_\bm=S_\bn$, i.e., that
 the inequality $\bm\geq\bn$ is anodyne. This proves (i)$\iff $(ii). 
 The equivalence (ii)$\iff $(iii) is clear since $p_{\bm,\bn}$ is a surjective map. 
  \qed

  \begin{prop}\label{prop:anodyne-unique}
 Suppose that  $\bm'\geq' \bn' \leq''\bn$  and
 at least one of the  inequalities
  is anodyne. Then $\sup(\bm',\bn)$ consists of exactly one element $\bm$. 
  If $\bm'\geq'\bn'$ is anodyne, then $\bm\geq'\bn$ is anodyne. If $\bn\geq''\bn'$ is
 anodyne,  then $\bm\geq''\bm'$ is anodyne. 
  \end{prop}
 
 \noindent{\sl Proof:} 
 We apply Proposition \ref{prop:sup=FP}. Note that the diagram 
 \eqref{eq:cart-diagr-W}  is Cartesian, being the external Cartesian product
 of two arrows
 \[
 \{ W/W_{I_1} \lra W/W_{J_1}\} \times  \{ W/W_{I_2} \lra W/W_{J_2}\}. 
 \]
 Suppose that $\bm'\geq' \bn'$ is anodyne. Then  the map $\pi_{\bm,\bn}: \bm'\to\bn'$
 (induced by $\pi'$) is a bijection. So  the arrow $\rho'$
   in the Cartesian product diagram
  \[
  \xymatrix{
  \bm'\times_{\bn'} \bn\ar[d]  \ar[r]^{\hskip .5cm \rho'}&\bn\ar[d]
  \\
  \bm' \ar[r]_{\pi'}&\bn'
  }
  \]
  is a bijection, i.e., $\bm'\times_{\bn'} \bn$ consists of one
  $W$-orbit. So  by Proposition \ref{prop:sup=FP},  
  there is only one possible $\bm$ and 
  by Proposition \ref{prop:ano=bij}, 
    $\bm\geq'\bn$
  is anodyne.
  The case when $\bn\geq''\bn'$ is anodyne, is treated similarly. \qed


\section{Main result: $\Perv(\hW)$ and mixed Bruhat sheaves}\label{sec:perv=mbs}
 
 \paragraph{Mixed Bruhat sheaves.}\label{par:MBS}

Let $\k$ be a field and  $\Vect_\k$   the category of finite-dimensional $\k$-vector spaces.

\begin{Defi}\label{def:MBS}
A {\em mixed Bruhat sheaf} of type $(W,S)$ is  the data of a  finite-dimensional $\k$-vector
spaces $E(\bm)$ for each  $\bm\in\Xi$ and linear maps
\[
 \del'_{\bm,\bn} = \del'_{\bm,\bn,E}: E(\bm)\lra E(\bn)
 \]
 for   $\bm\geq' \bn$ 
 and
\[
\del''_{\bm,\bn} = \del''_{\bm,\bn, E} : E(\bn)\lra E(\bm) 
 \]
 for  $\bm\geq'' \bn$,
satisfying the conditions:
\begin{itemize}
\item[(MBS1)] The $\del'_{\bm, \bn}$ are transitive, i.e., form a covariant functor from
the category formed by the poset $(\Xi, \geq')$  to $\Vect_\k$. Similarly, the
$\del''_{\bm,\bn}$ are transitive, i.e., form a contravariant functor $(\Xi, \geq'')\to\Vect_\k$. 

\item[(MBS2)] The $\del'$-  and $\del''$-maps commute with each other. That is, suppose
$\bm'\geq' \bn' \leq''\bn$. Then
\[
\del''_{\bn, \bn'} \del'_{\bm', \bn'} \,=
\sum_{\bm\in\sup(\bm',\bn)}  \del'_{\bm,\bn,} \del''_{\bm, \bm'}. 
\]

\item[(MBS3)] If $\bm\geq' \bn$ is an anodyne inequality, then $\del'_{\bm,\bn}$ is an isomorphism.
If  $\bm\geq'' \bn$ is an anodyne inequality, then $\del''_{\bm,\bn}$ is an isomorphism.

\end{itemize}
We denote by $\MBS=\MBS_{(W,S)}$ the category of mixed Bruhat sheaves of type $(W,S)$. 
\end{Defi}

\begin{prop}\label{prop:E-LS}
Let $X$ be a stratum of $\Sc^{(0)}$.
 The part of 
  a mixed Bruhat sheaf $E$ consisting of $E(\bm)$ with $U_\bm\subset X$ and
  the maps $(\del')^{-1}, \del''$ between such $E(\bm)$  gives a local system $\Lc_{E,X}$
  on $X$. 
   \end{prop}
   
     In particular, taking $X= W\backslash \hen^\reg$ to be the open stratum,
  we get, very directly, 
   a representation of the braid group $\Br_\gen = \pi_1(W\backslash \hen^\reg)$. 
   
   \vskip .2cm

   \noindent{\sl Proof:} 
   To define a local system $\Gc$ on $X$, we need to give vector spaces $\Gc_\bm$
   for each $U_\bm\subset X$ and  generalization maps
   $\gamma_{\bn, \bm}: \Gc_\bn\to \Gc_\bm$ for any inclusion $U_\bn\subset \ol U_\bm$
   with $U_\bm\subset X$, which are isomorphisms and satisfy the transitivity conditions,
   see Proposition \ref {prop:cell-sheaves}. 
   
   We put $\Gc_\bm=E(\bm)$. When 
   $U_\bn\subset \ol U_\bm$
   with $U_\bm\subset X$, we have an anodyne inequality $\bn\leq \bm$.
   Now, if  $\bn\leq'\bm$, then  we set $\gamma_{\bn,\bm} = (\del'_{\bm,\bn})^{-1}$,
   the inverse of the isomorphism $\del'_{\bm,\bn}$, see (MBS3)
   If $\bn\leq'' \bm$, then we set $\gamma_{\bn,\bm}=\del''_{\bm, \bn}$, also an isomorphism by
   (MBS3). 
   Proposition 
    \ref {prop:anodyne-unique} together with (MBS2)
    implies that the two types of isomorphisms $\gamma_{\bn,\bm}$
    thus defined  commute with each other. Together with the transitivity of the
    $\del'$ and o $\del''$, this implies that these two types of  $\gamma_{\bn,\bm}$ extend
    uniquely to  any $\bn\leq \bm$ such that $U_\bm\subset X$ and 
    that they are transitive. \qed 
   
    \vskip .2cm 
    
    Let $\tau:\Cc\times\Cc\to\Cc\times\Cc$ be the permutation of the factors. It induces
an involution $\tau:\Xi\to\Xi$. For $\bm\in\Xi$ we write $\bm^\tau=\tau(\bm)$. Thus, if
$\bm=W(\bc,\bd)$, the $\bm^\tau=W(\bd,\bc)$. 
The category $\MBS$ carries a perfect duality $E\mapsto E^\tau$, where the
``dual'' mixed Bruhat sheaf $E^\tau$ has 
\be\label{eq:MBS-dual}
 E^\tau(\bm) \,=\, E(\bm^\tau)^*,\quad 
\del'_{\bm,\bn, E^\tau} \,=\, (\del''_{\bn^\tau, \bm^\tau, E})^*, \quad 
\del''_{\bm,\bn, E^\tau} \,=\, (\del'_{\bn^\tau, \bm^\tau, E})^*. 
 \ee

\paragraph{Perverse sheaves and the main result.}
Let 
\[
\Perv(\hW)\, = \, \Perv(\hW, \Sc^{(0)})\, \subset \, D^b_{\Sc^{(0)}}\Sh(\hW)
\]
 be  the category of   (middle-perversity) perverse sheaves of $\k$-vector
spaces on $\hW$ which are (cohomologically)  constructible with respect to $\Sc^{(0)}$.
In this paper we use the standard normalization of perversity conditions from
\cite{BBD} so that a  local system on a smooth $d$-dimensional subvariety $Z$
is perverse, if put in degree $-{\dim_\CC Z}$.   
Thus, explicitly, $\Fc\in  D^b_{\Sc^{(0)}}\Sh(\hW)$ is perverse iff:
\begin{itemize}
\item[($\Perv^-$)] For each $q$, the sheaf $\ul H^q(\Fc)$ is supported on a complex analytic
subvariety of complex dimension $\leq -q$ (so that $\ul H^q(\Fc)=0$ for $q>0$).

  \item[($\Perv^+$)] Condition ($\Perv^-$) holds also for the Verdier dual complex $\DD(\Fc)$. 
\end{itemize}

\noindent 
By definition, the Verdier duality $\DD$ preserves $\Perv(\hW)$. 

Consider the involution
\be\label{eq:tau-h}
\tau: \hen\lra\hen, \quad x+iy \mapsto y+ix, \,\,\, x,y\in\hen_\RR. 
\ee
This involution commutes with the $W$-action and preserves the strata of $\Sc^{(0)}_\Hc$.
Therefore it descends to an involution  of $\hW$ which we denote by the
same letter $\tau: \hW\to\hW$ and which preserves the strata of $\Sc^{(0)}$. 
This means that the pullback $\tau^*$ preserves $\Perv(\hW)$. We define
the {\em twisted dual} of $\Fc\in\Perv(\hW)$ to be
\be\label{eq:tw-verdier}
\Fc^\tau = \tau^*(\DD(\Fc)) = \DD(\tau^*\Fc). 
\ee
Let $\sgn: W\to \k^*$ be the sign character of $W$, defined by $\sgn(s)=-1, \ s\in S$.
As we have the surjection
\[
\Br_{(W,S)} := \pi_1(W\backslash \hen^\reg) \lra W, 
\]
$\sgn$ is also a character of $\Br_{(W,S)}$. In particular, we have the $1$-dimensional local system
$\Lc_\sgn$ of $\k$-vector spaces on $W\backslash \hen^\reg$.

What follows is the main result of this paper. 
 
 \begin{thm}\label{thm:main}
 \begin{enumerate}
 \item[(a)] We have an equivalence of categories $\GG:  \MBS_{(W,S)} \to  \Perv(\hW, \Sc^{(0)})  $
 taking the duality \eqref{eq:MBS-dual} on $\MBS_{(W,S)}$ to
  the twisted Verdier duality \eqref{eq:tw-verdier}.

 \item[(b)] For  mixed Bruhat sheaf $E$, the  restriction of the perverse sheaf $\GG(E)$  to  the open stratum
  $W\backslash \hen^\reg$, is isomorphic to the shifted  local system $\Lc_E\otimes \Lc_\sgn[r]$. 
  Here $r=\dim_\CC\hen$ and $\Lc_E$ is the local system associated to $E$
  by Proposition \ref{prop:E-LS}. 
  \end{enumerate}
 \end{thm}
 
 The proof will be given in Sections \ref{sec:cous-bruhat}-\ref{sec:cous-perv}.

 \paragraph{Examples. The  bicube point of view.} \label{par:bicube}
 Many examples of mixed Bruhat sheaves appear, most immediately,  in the form of simpler diagrams
 which we call  { bicubes}. More precisely, 
 let $T$ be a finite set. By a {\em $T$-bicube} we mean a diagram $Q=(Q_I, u_{IJ}, v_{IJ})$,
 where:
 \begin{itemize}
 \item  $Q_I\in\Vect_\k$ is a vector space,  given for any subset $I\subset T$.
 
 \item $v_{IJ}: Q_I\to Q_J$ and $u_{IJ}: Q_J\to Q_I$ are linear maps given for any $I\subset J\subset T$
 and satisfying the transitivity properties:
 \[
 v_{II}=\Id, u_{II}=\Id, \quad v_{IK} = v_{JK} v_{IJ}, \,\, u_{IK} = u_{IJ} u_{JK}
  \]
  for any   $I\subset J \subset K$. 
 \end{itemize}
 In other words, a bicube consists of two commutative cubes superimposed on the same set of vertices
 such that the arrows in the two cubes going in the opposite directions. We denote by $\Bic_T$ the
 category of $T$-bicubes. 
 
 \vskip .2cm
 
 To a mixed Bruhat sheaf $E\in\MBS_{(W,S)}$, we attach an $S$-bicube $Q=\Qc(E)$ as follows.
 Recall that we identify $\Xi$ with $W\backslash (\Cc\times\Cc)$.
  For $I\subset S$ put
 \[
 \bm''_I = W(K_I,0), \quad \bm_I = W(K_I,K_I), \quad  \bm'_I = W(0,K_I)\quad \in\,\, \Xi. 
 \]
 Then $\bm'_I\leq ' \bm_I \geq'' \bm''_I$, both inequalities being anodyne. 
 We
 put $Q_I = E(\bm'_I)$
 and write $\phi_I: Q_I\to E(\bm''_I)$ for the composition
 \[
 \xymatrix{
 Q_I = E(\bm'_I)  \ar[rr]^{\hskip .5cm (\del'_{\bm_I, \bm'_I})^{-1}}  && E(\bm_I)  
 \ar[rr]^{\del''_{\bm_I, \bm''_I}}&& E(\bm''_I). 
 }
 \]
 If  $I\subset J$, then   $\bm'_I \geq' \bm'_J$
 and $\bm''_I \geq'' \bm''_J$. We define 
 \begin{align*}
 v_{IJ} &= \del'_{\bm'_I, \bm'_J}: Q_I \lra Q_J, 
 \\
  u_{IJ} &= \phi_I^{-1}\circ  \del''_{\bm''_I, \bm''_J}
 \circ \phi_J: Q_J\lra E_Q. 
 \end{align*}
  Transitivity of the $\del'$  and $\del''$ implies that $\Qc(E)$ is indeed a bicube. This gives a functor
  \[
 \Qc:  \MBS_{(W,S)} \lra \Bic_S, \quad E\mapsto \Qc(E). 
  \]
  
  \
  \begin{rem}\label{rem:faithful}
  Recently J. Tao, in his preprint
  ``An alternative definition of mixed Bruhat sheaves''
   has proved this
 this functor is fully faithful, i.e. that  any mixed Bruhat sheaf
  can be recovered from the corresponding bicube, see also \cite{KS-prob} where the essential image of 
  $\Qc$ is described for the $A_n$ series.
\end{rem}  
  
  \ 
  
  The following examples
  illustrate  fully faithfulness of $\Qc$  for $(W,S)$ be of types $A_r$, $r = 1, 2$.

 \begin{ex}\label{ex:sl2-1}
 Let $(W,S)$ be of type $A_1$. Then $\hen=\CC$, $\hen_\RR=\RR$, 
 and the arrangement of hyperplanes $\Hc$ in $\hen_\RR=\RR$ consists of one hyperplane 
 $\{0\}$.
 The Coxeter complex $\Cc$ consists of $\RR_{<0}, \{0\}$ and $\RR_{>0}$. 
 The Weyl group $W = \{1,s\}$, where $s: \hen\to\hen$ takes $z\mapsto -z$.
 The quotient $\hW$ is identified with $\CC$ by the function $z^2$. 
  Thus $\Perv(\hW,\Sc^{(0)})=\Perv(\CC,0)$ is the classical category of perverse
 sheaves on $\CC$ whose only possible singularity is $0$. 
The cell decompositions $\Sc_\Hc^{(2)}$ of $\hen=\CC$ and $\Sc^{(2)}$ of $\hW=\CC$
 are depicted in Fig.\ \ref {fig.sl2-1}.

  \begin{figure}[H]
 \centering
 
 \begin{tikzpicture}[scale=.4, baseline=(current bounding box.center)]
 
 \node at (0,0){$\bullet$};
  \draw (0,-4) -- (0,4); 
  \draw (4,-0) -- (-4,0); 
  \draw (4,4) --(4, -4) -- (-4,-4) -- (-4,4) -- (4,4); 
  \node at (-.5, -0.5) {$0$}; 
 \node at (3.5, 3.5){$\hen$};  
  \end{tikzpicture}
  \quad\quad 
 { \huge  $\buildrel z^2\over  \lra$}
   \quad\quad 
   \begin{tikzpicture}[scale=.4, baseline=(current bounding box.center)]
 
 \node at (0,0){$\bullet$};
  \draw (4,-0) -- (-4,0); 
  \draw (4,4) --(4, -4) -- (-4,-4) -- (-4,4) -- (4,4); 
  \node at (-.5, -0.5) {$0$}; 
 \node at (2.8, 3.3){$\hW$};  
 
  \end{tikzpicture}
 \caption{ Cell decompositions $\Sc^{(2)}_\Hc$ and $\Sc^{(2)}$ for $W$ of type $A_1$.}\label{fig.sl2-1}
 \end{figure}
 
  For simplicity we label the five cells of  $\Sc^{(2)}$
 by their representative points $0, \pm 1, \pm i$. These representative points provide therefore
 a labelling of the poset $\Xi$. We recall that an element
 $C_1\times C_2$ in $\Cc\times\Cc$ is interpreted as the product cell
 $iC_1+C_2\subset \hen_\CC$. With this understanding, 
  $\Xi$ consists of the elements
 \[
 \begin{gathered}
 \bm_0 = W(0,0), 
 \\
 \bm_{-1} = W(\RR_{>0}\times\{0\}) = W(\RR_{<0}\times\{0\}), \quad \bm_{1} = W(\{0\}\times
 \RR_{>0}) = W(\{0\}\times \RR_{<0}), 
 \\
 \bm_i = W(\RR_{>0}\times\RR_{>0}) = W(\RR_{<0}\times \RR_{<0}) ,\quad
 \bm_{-i} = W(\RR_{>0}\times\RR_{<0}) = W(\RR_{<0}\times \RR_{>0}).
 \end{gathered} 
 \]
 The partial order on $\Xi$ in this labelling is given by the arrows in the  diagram
 
 \def\ano{{\on{ano}}}
 
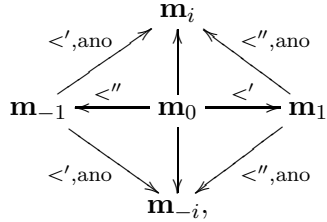
 \[
 \xymatrix{
 &\bm_i&
 \\
 \bm_{-1} \ar[ru]^{<', \ano}   \ar[rd]_{<', \ano}   & \bm_0  \ar[r]^{<'}  \ar[l]_{<''}  \ar[u] \ar[d]& \bm_1 \ar[lu]_{<'', \ano}  \ar[ld]^{<'', \ano} 
 \\
 & \bm_{-i}, & 
 }
 \]
in which the arrows corresponding to inequalities  $<'$ or $<''$ are marked
with these signs, and anodyne inequalities are marked ``$\ano$''.  
 A mixed Bruhat sheaf is
 a diagram
 \be\label{eq:MBS-sl2}
 \xymatrix{
 & E_i \ar[rd]^a& 
 \\
 E_{-1}\ar[rd]_d
  \ar[ru]^b \ar[r]^u& E_0 \ar[r]^v& E_1
 \\
 & E_{-i} \ar[ru]_c& 
 }
 \ee
 where $a,b,c,d$ are  isomorphisms and $uv=ab+cd$.
 A description of $\Perv(\CC,0)$ in
 terms  of  such 
  diagrams  is equivalent to the classical description in terms of diagrams 
      \be\label{eq:phi-psi}
      \xymatrix{
  \Phi
 \ar@<.4ex>[r]^{v}&\Psi
\ar@<.4ex>[l]^{u}   
} 
\ee
where $T_\Psi:=\Id_\Psi-vu$   is an isomorphism. 
That is, $E_0$ is identified with $\Phi$ and the other $4$ spaces are identified with $\Psi$.
Note that the diagram \eqref{eq:phi-psi} is an $S$-bicube, and passing from 
\eqref{eq:MBS-sl2} to  \eqref{eq:phi-psi} is a particular case of the functor $\Qc$.
So in this case $\Qc$ is fully faithful and its essential image is characterized by the  
invertibility
of $T_\Psi$. 
 \end{ex}

  \begin{ex}
  Let $(W,S)$ be of type $A_2$.
   In this case $W=S_3$, and the stratification $\Sc^{(0)}$
  of $\hW\simeq \CC^2$
  is given by the  cuspidal cubic
  \[
Y =  \bigl\{ (a,b)\in\CC^2 | \,\, 4a^3 + 27 b^2 = 0\bigr\}.
\]
 That is, the strata are $\CC^2\- Y$, $Y\-\{0\}$ and $\{0\}$. 
   Further, in this  case 
   $S$ consists of two elements which we denote
  $s_1$ and $s_2$. Accordingly, the bicube  associated to $E\in\MBS_{(W,S)}$
  is labelled by subsets of $\{1,2\}$ and has the (bisquare)  form
   
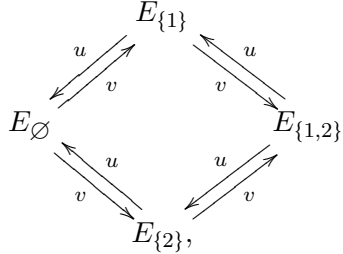
\[
 \xymatrix{
 & E_{\{1\}}
 \ar@<-.5ex>[dr]_{v}
 \ar@<-.5ex>[dl]_{u} &
 \\
 E_{\emptyset} 
\ar@<-.5ex>[ur]_{v}
 \ar@<-.5ex>[dr]_{v}
 && E_{\{1, 2\}}
 \ar@<-.5ex>[dl]_{u}
 \ar@<-.5ex>[ul]_{u}
 \\&E_{\{2\}}, 
 \ar@<-.5ex>[ur]_{v}
 \ar@<-.5ex>[ul]_{u}&
 }
 \]
  where we have omitted the indexing of the $v$- and $u$-maps. 
   A description of $\Perv(\hW)$ in this case was given in 
   \cite {GM-cusp} (see also  \cite  {KS-shuffle} \S 5.3), again in terms of bisquares 
      satisfying certain conditions.  Using this description, one checks directly
      that the functor $\Qc$ is 
      fully faithful in this case as well. 
      \end{ex}


\section{The Cousin complex of a mixed Bruhat sheaf}\label{sec:cous-bruhat}

Our proof of Theorem \ref{thm:main} is, similarly to \cite {  {KS-hyp-arr}, KS-hyp-arr-II}
 based on  associating to a mixed Bruhat sheaf $E$ a certain complex of sheaves 
$\Ec^\bullet = \Ec^\bullet(E)$ that we call its {\em Cousin complex}. A priori, $\Ec^\bullet$
is only $\RR$-constructible but it turns to be (cohomologically) constructible with respect to $\Sc^{(0)}$
and, moreover, a perverse sheaf. Here we describe this construction. 

\paragraph{Imaginary strata in $\hW$.} Since the action of $W$ on $\hen$ is induced
by its action on $\hen_\RR$, the 
 ``imaginary part'' map $\Im: \CC\to\RR$
 induces the map 
 \[
 \Ien: \hW\lra W\backslash \hen_\RR. 
 \]
 Recall that $W\backslash \hen_\RR= \bigsqcup_{I\subset S} K_I$. 
The sets
 \[
 X_I^\Im = \Ien^{-1}(K_I) \xhookrightarrow{ k''_I}  \hW.
 \]
   form a locally closed decomposition (see
   Definition \ref{def:lcd}) 
   of $\hW$ which we call the 
 {\em imaginary decomposition} and denote by $\Sc^\Im$. 
 The notation $k''_I$ is chosen to match the more systematic notation used
 below in \S \ref{par:perv-cous-com}
  For  $(W,S)$ of type $A_n$ the
   decomposition  $\Sc^\Im$ was
   considered in \cite{KS-shuffle}. 
   
 We refer to the $X_I^\Im$ as {\em imaginary strata}.   Note that
   \[
   X_I^\Im \,=\,\bigsqcup_J \bigsqcup_{\bm\in \Xi(I,J)} U_\bm,
   \]
so  $\Sc^{(2)}$ refines $\Sc^\Im$. For any two  cells $U_\bm, U_\bn\subset X_I^\Im$, 
the inclusion $U_\bn\subset\ol  U_\bm$ is equivalent to $\bn \leq '' \bm$.

\begin{ex}
Let  $W$ be of type $A_1$. The stratification
$\Sc^{(2)}$ is this case was discussed in Example \ref{ex:sl2-1}
 depicted in Fig.\ \ref{fig.sl2-1}. In terms of $\Sc^{(2)}$, 
the locally closed decomposition  $\Sc^\Im$ is described as follows. It has two strata, 
 depicted in Fig.\ \ref{fig.sl2-2},
with $K_\emptyset$ being the union of three cells of $\Sc^{(2)}$ and $K_{\{\alpha\}}$
being the union of two such cells. We see, in particular, that $K_{\{\alpha\}}$ is not
a submanifold and therefore $\Sc^\Im$ is not a stratification in the sense of
Definition  \ref{def:lcd}(c).  For further comparison, the stratification
$\Sc^{(0)}$ in this case consists of $\{0\}$ and $\CC\- \{0\}$. 

 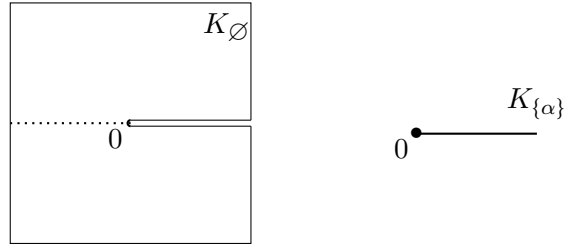
\begin{figure}[H]
 \centering
 
 \begin{tikzpicture}[scale=.4, baseline=(current bounding box.center)]

  \draw (0,0.1) -- (4, 0.1); 
    \draw (0,-0.1) -- (4, -0.1); 
    
    \draw [line width=0.8] (0,0.1) arc (90:270:0.1); 
  
  \draw (4,0.1)--  (4,4) --(-4, 4) -- (-4,-4) -- (4,-4) -- (4,-0.1); 
  
  \draw [dotted, line width=0.8] (0,0) -- (-4,0); 
  \node at (-.5, -0.5) {$0$}; 
 \node at (3.2, 3.2){$K_\emptyset$};  
  \end{tikzpicture}
  \quad\quad 
  \quad\quad 
   \begin{tikzpicture}[scale=.4, baseline=(current bounding box.center)]
 
 \node at (0,0){$\bullet$};
   \node at (-.5, -0.5) {$0$}; 
   \draw[line width=0.8] (0,0) --(4,0); 
 \node at (4,1){$K_{\{\alpha\}}$};

  \end{tikzpicture}
 \caption{Imaginary strata  for $W$ of type $A_1$.}\label{fig.sl2-2}
 \end{figure}
\end{ex}

\paragraph{ The Cousin complex.} Given a mixed Bruhat sheaf $E=(E(\bm), \del',\del'')$,
we can define a cellular sheaf $\wt\Ec_I$ on $X_I^\Im$ where stalk at 
any point of $U_\bm$ 
is $E(\bm)$ and the generalization map $E(\bn)\to E(\bm)$ for 
 $U_\bn\subset\ol  U_\bm$ is $\del''_{\bm,\bn}$. By the transitivity of the
 $\del''$-maps in   (MBS1), this gives a well-defined sheaf $\wt \Ec_I$. 
 We further put $\Ec_I = k''_{I*} \wt\Ec_I$, a sheaf on $\hW$. 
 
 \vskip .2cm
 
 For $I\subset S$ let $\k^I$ be the $\k$-vector space spanned by $I$, i.e., a vector
 space with basis $\{e_\alpha\}_{\alpha\in I}$. Let $\det(I)=\Lambda^{|I|} (\k^I)$ be
 the top exterior power of $\k^I$. For $I_1\subset I_2$ such that $|I_2|=|I_1|+1$,
 i.e., $I_2= I_1\sqcup \{\alpha\}$ for some $\alpha$, we have the map
 \[
 \eps_{I_1, I_2}: \det(I_1) \lra \det(I_2), \quad v\mapsto v\wedge e_\alpha. 
 \]
We now define the complex of sheaves
\be\label{eq:cousin}
\Ec^\bullet=\Ec^\bullet(E) \,=\,\biggl\{ \Ec_\emptyset \buildrel d\over\to  \bigoplus_{|I|=1} \Ec_I\otimes\det(I) \buildrel d\over\to  \bigoplus_{|I|=2} \Ec_I\otimes\det(I) \buildrel d\over\to \cdots
\buildrel d\over\to \Ec_{S}\otimes \det(S)
\biggr\}
\ee
graded so that $E_\emptyset$ in in degree $-{\dim_\CC(\hen)}$. 
The differential $d$ is induced by the maps $\del'$. More precisely, 
suppose that    $I_2= I_1\sqcup \{\alpha\}$ for some $\alpha$. Let $\bm\in\Xi(I_2, J)$ so 
 $U_\bm\subset X_{I_2}^\Im$. The definition $\Ec_{I_1} = (k''_{I_1})_*\wt\Ec_{I_1}$
 implies that  $\Ec_{I_1}$ is locally constant
 on $U_\bm$, hence globally constant because $U_\bm$ is simply connected. Its stalk there is  
 \[
 (\Ec_{I_1})_{U_\bm}\,:=\, \Gamma (U_\bm, \Ec_{I_1}) \,=\, \bigoplus_{ \bn\in\Xi(I_1,J)\atop
 \bn\geq ' \bm
 } E(\bn). 
 \]
 The stalk  $(\Ec_{I_2})_{U_\bm}$  is, by definition, $E(\bm)$. 
Now, the matrix element
\be\label{eq:d-I1-I2}
d_{I_1, I_2}: \Ec_{I_1}\otimes\det(I_1) \lra\Ec_{I_2}\otimes\det(I_2)
\ee
is defined over $U_\bm$  by the map
\[
d_{I_1, I_2,\bm}: 
\sum_{  \bn\in\Xi(I_1,J)\atop
 \bn\geq ' \bm
 } \del'_{\bn, \bm} \otimes \eps_{I_1, I_2} : \bigoplus_{  \bn\in\Xi(I_1,J)\atop
 \bn\geq ' \bm}
   E(\bn)\otimes\det(I_1) \lra E(\bm)\otimes\det(I_2). 
\]

\begin{prop} With the above definitions of $d_{I_1, I_2, \bm}$,
\begin{enumerate}
\item[(a)] The maps $d_{I_1, I_2 }$  is  a morphism of sheaves. 

\item[(b)] The morphisms  $d$ with  matrix elements $d_{I_1, I_2}$ define a
complex of sheaves $\Ec^\bullet$ as in \eqref{eq:cousin}, i.e., they  satisfy $d^2=0$. 
\end{enumerate}
\end{prop}

\noindent{\sl Proof:} (a) follows at once from (MBS2) (commutativity of $\del'$ with
$\del''$), while (b) follows from (MBS1) (transitivity of $\del'$). \qed

\vskip .2cm

We call $\Ec^\bullet$ the {\em Cousin complex} associated to $E$. 
Theorem \ref{thm:main} will be a consequence of the following more precise result.

\begin{thm}\label{thm:main-II}
\begin{enumerate}
\item[(a)] For $E\in \MBS_{(W,S)}$ the complex $\Ec^\bullet(E)$ is an object of $\Perv(\hW, \Sc^{(0)})$. 
In particular, it is (cohomologically) $\Sc^{(0)}$-constructible.

\item[(b)]    $\Ec^\bullet(E^\tau)$, defined by  \eqref{eq:MBS-dual},
 is naturally quasi-isomorphic to the twisted dual
$(\Ec^\bullet(E))^\tau$, defined by \eqref{eq:tw-verdier}.

\item[(c)] The functor $\GG: E\mapsto \Ec^\bullet(E)$ is an equivalence of categories
\[
\GG: \MBS_{(W,S)}\to\Perv(\hW,\Sc^{(0)}).
\]

\item[(d)]  For  $E\in \MBS_{(W,S)}$, the  restriction of  $\GG(E)$  to  the open stratum
  $W\backslash \hen^\reg$
  is isomorphic to the shifted  local system $\Lc_E\otimes \Lc_\sgn[r]$. 
  \end{enumerate}
\end{thm}

\paragraph{ The Fox-Neuwirth-Fuchs cells.}  
For any subset $S\subset \hen_\RR$ let $\Lin_\RR(S)$ be the $\RR$-linear
subspace spanned by $L$.

\vskip .2cm

Denote by $\Sc_\Hc^{(1)}$  the ``intermediate'', or  {\em Bj\"orner-Ziegler}
stratification of $\hen$ from 
  \cite{BZ} and  \cite[\S 2]{KS-hyp-arr}.   
  This is a quasi-regular cell decomposition of $\hen$ into cells $[C,D]$ labelled
by {\em face intervals}, i.e., pairs $(C,D)\in\Cc$ such that $C\leq D$. By definition,
$[C,D]$ consists of $x+iy\in \hen$ with $x,y\in\hen_\RR$ such that:
\begin{itemize}
\item[(a)] $y\in C$.

\item[(b)] $x$ is congruent to an element of $D$ modulo the subspace $\Lin_\RR(C)\subset\hen_\RR$. 
\end{itemize}
Thus, in the notation of the Appendix, 
\[
\Sc^{(2)}_\Hc \prec \Sc^{(1)}_\Hc \prec \Sc^{(0)}_\Hc. 
\]
The action of $W$ on $\hen$ preserves  $\Sc^{(1)}_\Hc$
and so $\Sc^{(1)}_\Hc$ descends to a stratification $\Sc^{(1)}$ of $\hW$ such that
\[
\Sc^{(2)} \prec \Sc^{(1)} \prec \Sc^{(0)}. 
\]
 
 \begin {prop}
 Every stratum of $\Sc^{(1)}$ is a topological cell, so $\Sc^{(1)}$ is a
 (not necessarily quasi-regular) cell decomposition of $\hW$. 
 \end{prop}
 
 \noindent{\sl Proof:} Since each $[C,D]$ is   a cell, it suffices to prove
 that $p: [C,D] \to p([C,D])$ is a homeomorphism, i.e.,  that if $z,z' \in [C,D]$ and $z'=wz$ for
 some $w\in W$, then $z'=z$. 
 
 Let $z=x+iy$ and  $z'=x'+iy'$ with $x,y,x',y'\in\hen_\RR$. Then $z'=wz$ implies $y'=wy$. But since
 $y,y'$ lie in the same face $C$ of $\Cc$, we have $y'=y$. This  also means that
  $w$ stabilizes $C$ pointwise and therefore stabiizes  
   the vector space $L:=\Lin_\RR(C)$, which
  is a flat of $\Hc$, pointwise. 
  
  By conjugating with an appropriate element of $W$, we can assume that $L$ is
  a flat whose generic part lies in the closure of the dominant Weyl chamber $C^+$,
  i.e., $L=\bigcap_{\alpha\in I} \alpha^\perp$ for some $I\subset S$. Then the
  subgroup in $W$ fixing $L$ pointwise is $W_I\subset W$. In this case the quotient
  arrangement $\Hc/L$ in $\hen_\RR/L$ (see  \cite[\S 2B] {KS-hyp-arr}) is the root
  arrangement associated to the semi-simplification of the Levi subalgebra associated to $I$.
  The group $W_I$ is the Weyl group of that semi-simplification. In particular,
  distinct on the same face of $\hen/L$ cannot be congruent under the action of $W_I$. 
  
  Now, the condition $z,z'\in[C,D]$ means $x=d+l$ and $x'=d'+l'$ with $d,d'\in D$ and $l,l'\in L$. 
  Moreover, $w(z)=z'$ implies $w(x)=x'$, so 
  \[
  w)+l = w(d+l) = w(x) = x'=d'+l', 
    \]
 i.e., $w(d)=d'+l'-l$,
which means that the images of $d$ and $d'$ in $\hen_\RR/L$, which lie in the
  same face of $\Hc/L$, are  congruent under the action of $w\in W_I$. 
  So these images are equal, i.e., $d'=d+l''$ for some $l''\in L$.
  Combining this with the above equality, we get
  \[
  wd = d+\lambda, \quad \text{where} \quad  \lambda := l''+l'-l\in L. 
  \]
 Note that $\lambda$ is fixed by $w$. For $\lambda\neq 0$ the last equality is impossible, since
 for any $n>0$ we have $w^n(d)= d+n\lambda$ which contradicts the fact that, $W$ being
 finite, we must have $w^n=\Id$ for some $n$. 
 So $\lambda=0$ and  $w(d)=d$, and therefore
 \[
 wx = w(d+l) = wd+l = d+l = x.
 \]
  Together with the equalty $w(y)=y$ proved earlier, this implies $z'=wz=z$. \qed
  
  \vskip .2cm
  
  We will call the cells $p([C,D])$ of the cell decomposition $\Sc^{(1)}$ the 
  {\em Fox-Neuwirth-Fuchs} cells of $\hW$. 
  
  \begin{exas}
  (a) For $W$ of type $A_1$ the Fox-Neuwirth-Fuchs cells  $\Sc^{(1)}$ of $\hW=\CC$ consist of 
  $\CC\- \RR_{\geq 0}$, $\RR_{>0}$
and $\{0\}$, see Fig.\ \ref{fig.sl2-2}.

(b) For $W = S_n$ acting on $\hen_\RR=\RR^n$ the Fox-Neuwirth-Fuchs cells of $\hW=\Sym^n(\CC)$
were discussed in \cite{KS-contingency}.
  \end{exas}
  
  Below, we refer to the Appendix for the meaning of the notations $\wedge, \vee$.
    
  \begin{prop}\label{prop:S-vee-tS}
  \begin{enumerate}
  \item[(a)] We have $\Sc^{(1)} = \Sc^\Im\wedge \Sc^{(0)}$. In particular, $\Sc^{(1)}$ refines both
  $\Sc^\Im$ and $\Sc^{(0)}$.

 \item[(b)] We also have $\Sc^{(0)}= \Sc^{(1)}\vee \tau(\Sc^{(1)})$, where $\tau: \hW\to\hW$
  is the involution \eqref{eq:tau-h}. In particular
  (Proposition \ref{prop:strat-append}(b)),  any sheaf which is $\Sc^{(1)}$-constructible
  and $\tau(\Sc^{(1)})$-constructible is $\Sc^{(0)}$-constructible. 
  \end{enumerate}
  \end{prop}
  
  \noindent{\sl Proof:} As $\Sc^{(2)}$ refines $\Sc^\Im$,  $\Sc^{(1)}$ and $\Sc^{(0)}$, we have
   equivalence relations $\equiv_\Im$, $\equiv_{\Sc^{(1)}}$ and $\equiv_{\Sc^{(0)}}$
  on the set
  $\Xi$ labelling cells of $\Sc^{(2)}$, which describe when two cells lie in the
  same stratum of the corresponding coarser decomposition. 
  
  Now, by definition, $\equiv_\Im$ is the equivalence  closure of  the relation $\geq''$.
  At the same time, $\Sc^{(1)}$ is the equivalence closure of the relation $R''$  where
   $\bm R'' \bn$ if, first,  $\bm\geq'' \bn$ and, second,  the inequality is anodyne, i.e.,
  $\bm\equiv_{\Sc^{(0)}} \bn$. This implies (a). 
  
  \vskip .2cm
  
  Let us prove (b). Note that  $\Sc^{(2)}$ also refines $\tau(\Sc^{(1)})$ and so we have the
  equivalence relation $\equiv_{\tau(\Sc^{(1)})}$ on $\Xi$ describing how   $\Sc^{(2)}$-cells
  are arranged into $\tau(\Sc^{(1)})$-cells. As before,  $\equiv_{\tau(\Sc^{(1)})}$
  is the equivalence closure of the relation $R'$ defined
  by $\bm R' \bn$ if, first,  $\bm\geq' \bn$ and, second,  the inequality is anodyne. 
  It follows that $R'\cup R''$, considered 
  as a subset of $\Xi\times\Xi$, is contained in $\equiv_{\Sc^{(0)}}$,
  and so $(R'\cup R'')^\sim$, the equivalence closure of $R'\cup R''$,
   is contained in $\equiv_{\Sc^{(0)}}$. 
 Now, $\Sc^{(1)}\vee\tau(\Sc^{(1)})$ is, by definition,  the partition of $\hW$ into unions of
 $\Sc^{(2)}$-cells corresponding to the classes of $(R'\cup R'')^\sim$. So each part
 of $\Sc^{(1)}\vee\tau(\Sc^{(1)})$ is contained in a single $\Sc^{(0)}$-stratum.

 Conversely, given an  $\Sc^{(0)}$-stratum $S$, we consider the corresponding
  $\equiv_{\Sc^{(0)}}$-class $\Xi_S\subset\Xi$. Since the $U_\bm$ with $\bm\in\Xi_S$
  form a quasi-regular cell decomposition of $S$, the  sub-poset $\Xi_S\subset\Xi$
  is closed under taking intermediate extensions. That is, if 
  $\bm,\bn, \bp\in\Xi$ are such that  $\bm \leq\bn \leq \bp$ and $\bm, \bp\in\Xi_S$,
  then $\bn\in\Xi_S$. Now, $S$ being connected,  any two $\bm,\bn\in \Xi_S$
  are connected by a chain of anodyne $\geq, \leq$- inequalities. But any anodyne
  inequality, say  $\bm\geq \bn$,  factors into two anodyne inequalities
  $\bm\geq' \bm'\geq'' \bn$. This means that $\Xi_S$ is a class of $(R'\cup R'')^\sim$,
  thus proving  (b). \qed
  
  \paragraph{ Perversity of the Cousin complex.}\label{par:perv-cous-com}
   Here we prove parts (a) and (b) of
  Theorem \ref{thm:main-II}. The argument is similar to that of 
  \cite [\S 5-6] {KS-hyp-arr-II}, so we give a more condensed presentation.
  
  \vskip .2cm
  
  For $I\subset S$ let $X_I^\Re=\tau(X_I^\Im)$. Thus $X_I^\Re = \Ren^{-1}(K_I)$,
  where $\Ren: \hW\to W\backslash \hen_\RR$ is induced by 
  the ``real part'' map $\Re: \CC\to\RR$. 
  
  For $\bm\in\Xi(I,J)$ we have a commutative diagram of embeddings
  \[
  \xymatrix{
  U_\bm \ar[d]_{j''_\bm}
  \ar[r]^{j'_\bm} & X_I^\Im
  \ar[d]^{k''_I}
  \\
  X_J^\Re \ar[r]_{k'_J}& \hW. 
  }
  \]
 
  \begin{lem}\label{lem:*!=!*}
  For any $V\in\Vect_\k$ we have canonical isomorphisms
  \[
  \begin{gathered}
  (k'_J)_! \, (j''_\bm)^*\, \ul V_{U_\bm} \buildrel \simeq \over\lra (k''_I)_*\,  (j'_\bm)_!\,  \ul V_{U_\bm}, 
  \\
  (k''_I)_!\,  (j'_\bm)_* \,  \ul V_{U_\bm} \buildrel \simeq \over\lra (k'_J)_*\,  (j''_\bm)_!\,   \ul V_{U_\bm}. 
  \end{gathered}
  \]
  \end{lem}
  
  \noindent{\sl Proof:} Let $\bm=W(\bc,\bd)$. By Proposition \ref{prop:p-iA+B-homeo}, 
  $p:\hen \to \hW$ induces a homeomorphism of closures $i\ol A_\bc+ \ol A_\bd \to \ol U_\bm$.
  So our statement reduces to the similar statement about the product diagram
  \[
  \xymatrix{
  iA_\bc + A_\bd \ar[r]\ar[d]& i\hen_\RR + A_\bd \ar[d]
  \\
  iA_\bc + \hen_\RR \ar[r]& \hen,
  }
  \]
  which is  \cite[Prop. 5.2]  {KS-hyp-arr-II}. \qed
  
  \vskip .2cm
  
For any $J\subset S$ let $\ol J= S\- J$
  be the complement. Let also $r=\dim_\CC(\hW)$.
  
  \begin{prop}\label{prop:E-and D}
    For any $E\in\MBS$  
 we have a natural  isomorphism
  $\DD(\Ec(E))\simeq \tau^* \Ec(E^\tau)$  in the derived category. 
  \end{prop}
  
  \noindent{\sl Proof:} 
  By the definition of $\wt\Ec_I = \wt\Ec_I(E)$
  as a cellular sheaf on $X_I^\Im$ with stalks $E(\bm)$
  and   generalization maps  $\del''_{\bm,\bn}$,   the sheaf $\Ec_I = j_{I*}\wt\Ec_I$
  is isomorphic to
   the total complex of:
  \[
  \bigoplus_{|\ol J|=0} \, \bigoplus _{\bm\in\Xi(I,J)} (k''_I)_* \, (j'_\bm)_! \, \ul{E(\bm)}_{U_\bm} 
  \lra 
   \bigoplus_{|\ol J|=1}\,  \bigoplus _{\bm\in\Xi(I,J)} (k''_I)_* \, (j'_\bm)_! \, \ul{E(\bm)}_{U_\bm}[1] 
   \to\cdots 
  \]
see \cite[(1.12)] {KS-hyp-arr}.  So the Cousin complex $\Ec^\bullet(E)$ is isomorphic to
the total object of the double complex in  the derived category whose $(p,q)$th term is
 \be\label{eq:DC-1}
 \bigoplus_{   |I| =r+q\atop
 | J|=r-p } \, \bigoplus _{\bm\in\Xi(I,J)} (k''_I)_* \, (j'_\bm)_! \, \ul{E(\bm)}_{U_\bm}[p],
\ee
the horizontal differentials (corresponding to the generalization maps of the $\wt \Ec_I$) 
being given by the $\del''$, and the vertical differentials  by the $\del'$. 

Let us apply the Verdier duality $\DD$ to this double complex. Recall that $\DD$
interchanges $*$ with $!$, and that for a constant sheaf on a cell, we have
\[
\DD\left(\ul{E(\bm)}_{U_\bm}\right) \,=\, \ul{E(\bm)^*}_{U_\bm} \, [\dim_\RR U_\bm]. 
\]
If $\bm\in\Xi(I,J)$, then $\dim_\RR U_\bm=2r-|I|-|J|$. 
Therefore $\DD(\Ec^\bullet (E))$ is quasi-isomorphic to the total object of the double complex
in the derived category whose $(p,q)$th term is
\be\label{eq:DC-2}
\bigoplus_{|I|=r-q\atop |J|=r+p} (k''_I)_!\, (j'_\bm)_* \, \ul{E(\bm)^*}_{U_\bm} [q],
\ee
the horizontal differentials given by the duals to the $\del''_E$ and the
vertical differentials given by the duals to the $\del'_E$. Taking into account Lemma 
\ref{lem:*!=!*}, we recognize in \eqref{eq:DC-2} a version of the double complex
 \eqref{eq:DC-1}  but with $E^\tau$ instead of $E$ and with the roles of the real and imaginary
 parts exchanged. In other words, we recognize $\tau^*\Ec^\bullet(E^\tau)$, whence the
 claim. \qed
 
 \begin{cor}\label{cor:E(E)-constr}
 $\Ec^\bullet(E)$ is (cohomologically) $\Sc^{(0)}$-constructible. 
 \end{cor}
 
 \noindent{\sl Proof:} Recall that, in the notation of
 the proof of Proposition \ref{prop:S-vee-tS}  the equivalence relation
 $\equiv_{\Sc^{(1)}}$ on the set $\Xi$ parametrizing the cells of $\Sc^{(2)}$, is the equivalence
 closure of the relation $R''$. Since $\del''_{\bm,\bn}$ is an isomorphism for anodyne 
 $\bm\geq'' \bn$,  each $\Ec_I(E)$ is $\Sc^{(1)}$-constructible and so
 $\Ec^\bullet(E)$ is   cohomologically $\Sc^{(1)}$-constructible. 
 
 On the other hand, Proposition \ref{prop:E-and D} implies, in the same way,
 that $\DD(\Ec^\bullet(E))$ is quasi-isomorphic to a complex of which each term
 is $\tau(\Sc^{(1)})$-constructible, and so it is  cohomologically $\tau(\Sc^{(1)})$-constructible. 
Since Verdier duality preserves cohomologically constructible complexes, $\Ec^\bullet(E)$
is also cohomologically $\tau(\Sc^{(1)})$-constructible. Our statement now follows
from Proposition \ref{prop:S-vee-tS}(b). \qed

\begin{prop}\label{prop:E(E)-perv}
$\Ec^\bullet(E)$ is perverse. 
\end{prop}

\noindent{\sl Proof:} Let us prove ($\Perv^-$). Let $q\in\ZZ$. 
Since $\Ec^\bullet(E)$  is $\Sc^{(0)}$-constructible, the set 
$Z=\Supp\, \ul H^q(\Ec^\bullet(E))$ is a complex algebraic subvariety in $\hW$.  Also, 
\[
Z\,\subset \, \Supp \, \Ec^q(E) \,=\, \bigcup_{|I|\geq r+q}  X_I^\Im. 
\]
Let $\wt Z=p^{-1}(Z)\subset\hen$. We conclude that for any $z=x+iy\in\wt Z$ with
$x,y\in\hen_\RR$, the point $y$ lies in the union of  the (real) flats of $\Hc$ of
codimension $\geq r+q$. But because $\wt Z$ is a complex algebraic subvariety of $\hen$,
this implies that $\wt Z$ lies in the union of  the (complex) flats of $\Hc_\CC$
of (complex) codimension $\geq r+q$. Since $p:\hen\to\hW$ is a finite map,
$\dim_\CC(\wt Z)=\dim_\CC(Z)$, and we obtain $(\Perv^-)$ for $\Ec^\bullet(E)$.
The condition $(\Perv^+)$ follows from this by Proposition \ref{prop:E-and D}. \qed

\vskip .2cm

Propositions  \ref{prop:E-and D}, \ref{prop:E(E)-perv} and Corollary \ref{cor:E(E)-constr} imply parts (a) and (b) of Theorem \ref{thm:main-II}. 


\section{The Cousin complex of a perverse sheaf}\label{sec:cous-perv}

Here, we prove part (c) of Theorem  \ref{thm:main-II} by constructing a quasi-inverse
to the functor
\[
\GG: \MBS_{(W,S)}\lra \Perv(\hW),\quad E\mapsto \Ec^\bullet(E).
\]
For this, similarly to \cite{KS-hyp-arr,  KS-shuffle,  KS-hyp-arr-II}, we start from a perverse sheaf
$\Fc$ and construct geometrically a Cousin-type resolution of $\Fc$.

\paragraph{Cousin complex II.} Recall the embedding $k''_I: X_I^\Im\hra \hW$ of the imaginary stratum.
Again, $r=\dim_\CC\hW$.

\begin{prop}
Let $\Fc\in\Perv(\hW)$. Then:
\begin{enumerate}

\item[(a)] The complex $(k''_I)^!\Fc$ is quasi-isomorphic to a single sheaf $\wt\Ec_I=\wt\Ec_I(\Fc)$
in degree $|I|-r$.

\item[(b)] The complex $k''_{I*}\wt\Ec_I$ is quasi-isomorphic to a single sheaf
$\Ec_I=\Ec_I(\Fc) = R^0k''_{I*} (k''_I)^!\Fc$.

\item[(c)] $\Fc$ has an explicit  Cousin resolution $\Ec^\bullet(\Fc)$ of the form
\[
 \Ec_\emptyset(\Fc) \buildrel \delta\over\to \bigoplus_{|I|=1}\Ec_I(\Fc)
\buildrel\delta\over\to \cdots \buildrel\delta\over\to  \Ec_{S}(\Fc),
\]
graded so that $\Ec_\emptyset(\Fc)$ is in degree $-r$. 
\end{enumerate}
\end{prop}

\noindent{\sl Proof:}  (a) This is analogous to \cite[Prop. 2.3.7(a)]{ KS-shuffle} which  is
particular case $W = S_n$. 
To start, note that $p^*\Fc\in\Perv(\hen,\Sc^{(0)}_\Hc)$ and $\Fc=(p_*p^*\Fc)^W$. 
In terms of the  Cartesian square
\be\label{eq:cart-sq}
\xymatrix{
\hen \ar[rr]^p&& \hW
\\
\bigsqcup\limits_{C\in\Cc\atop p(C)=K_I} \hen_\RR+iC \ar[rr]_{\hskip 1cm p_I} 
\ar[u]^{l_I} 
&&X_I^\Im
\ar[u]_{k''_I}
}
\ee
we have 
\be\label{eq:jI*F}
(k''_I)^!\Fc\,\simeq \, (p_{I*} \, l_I^!\,  p^*\Fc)^W. 
\ee
But $l_I$ is the disjoint union of the embeddings $l_C: \hen_\RR+iC\hra\hen$.
It remains to notice that each $l_C^!\, p^*\Fc$ is quasi-isomorphic to a single sheaf
in degree $\codim(C)-r = |I|-r$, by \cite[Cor. 4.11] {KS-hyp-arr}. We also need to take into account
the difference in normalizations  for perverse sheaves in {\em loc. cit}. 

\vskip .2cm

(b) Let $z\in\hW$. The stalk at $z$ of $R^q  k''_{I*}\, \wt\Ec_I$ is, by definition,
 $H^q(U\cap X_I^\Im, \wt\Ec_I)$
for a small ball $U$ around $z$. By \eqref{eq:jI*F} and the Cartesian square 
\eqref{eq:cart-sq}, this $H^q$ is a subspace in the direct sum
\[
\bigoplus_{p(C)=K_I} H^q\bigl( p^{-1}(U)\cap (\hen_\RR+iC), \, l_C^!\, p^*\Fc\bigr). 
\]
 Let us prove that for $q>0$ each summand  vanishes. Indeed, since
 $p^{-1}(U)$ is a disjoint union of balls, this vanishing follows from a similar statement about
 perverse sheaves on  real hyperplane arrangements, namely
 \[
 R^q l_{C*} \, ( l_C^!\, p^*\Fc[\codim(C)-r]) =0 
 \]
 for $q\neq 0$, which is \cite[Cor. 4.11(a)]{KS-hyp-arr}. 

\vskip .2cm

(c) This is a  purely formal consequence of
 (a) and (b) via  the Postnikov system associated
to $\Fc$ and the increasing filtration of $\hW$ by closed subspaces 
$X_{\leq m}^I = \bigcup_{|S\- I|\leq m} X_I^\Im$, see  \cite[\S 1B]{KS-hyp-arr}.  

\paragraph{From a perverse sheaf to a mixed Bruhat sheaf.} 
Let $\Fc\in\Perv(\hW)$. Because of Proposition \ref{prop:S-vee-tS}(a), each $\Ec_I(\Fc)$
is $\Sc^{(1)}$-constructible, hence $\Sc^{(2)}$-constructible. 
For $\bm\in\Xi(I,J)$ let $E(\bm)$ be the stalk of $\Ec_I(\Fc)\otimes\det(I)^{\otimes (-1)}$
at $U_\bm$. The generalization maps of the $\Ec_I(\Fc)$ and the differential $\delta$ in
the complex $\Ec^\bullet(\Fc)$ translate directly into linear maps $\del''_{\bm,\bn}$,\
$\del'_{\bm,\bn}$  that which satisfy (MBS1)-(MBS2)
as in Definition \ref{def:MBS}.
More precisely, the transitivity of the generalization maps gives 
the transitivity of the $\del''$,
the condition $\delta^2=0$ gives the transitivity of the $\del'$, and the fact that
$\delta$ is a morphism of cellular sheaves gives (MBS2). 

\begin{prop}\label{prop:E=MBS}
The diagram $E=(E(\bm),\del', \del'')$  also satisfies (MBS3), so it is a mixed Bruhat sheaf. 
\end{prop}

\noindent{\sl Proof:} $\Sc^{(1)}$-constructibility of each $\Ec_I(\Fc)$ gives one half
of (MBS3): $\del''_{\bm, \bn}$ is an isomorphism for anodyne $\bm\geq'' \bn$. Let us
prove the other half. Similarly
to the proof of Proposition \ref{prop:E-and D}, 
  we can represent $\Ec^\bullet(\Fc)$,
in the derived category
by the total object of the double complex
  consisting of shifted sheaves of the form
$(k_I'')_*\, (j_\bm')_! \, \ul{E(\bm)}_{U_m}$, then apply  Verdier duality. This will give
an explicit complex of sheaves $\Gc^\bullet$ which, on one hand,  is quasi-isomorphic
to $\DD(\Fc)$ and, on the other hand, has the form
\[
 k'_{\emptyset *}\,\wt\Gc_\emptyset \lra \bigoplus_{|I|=1}
k'_{I*}\wt\Gc_I \lra \cdots
\]
with the leftmost term in degree $-r$. Here $k'_I: X_I^\Re= \tau(X_I^\Im)\hra\hW$ is the embedding and
$\wt\Gc_I$ is an $\Sc^{(2)}$-cellular sheaf on $X_I^\Re$. Explicitly, the stalks of
$\wt\Gc_I$ are the $E(\bm)^*$ and the generalization maps are the $(\del'_{\bm,\bn})^*$.
Note that for such $\Gc^\bullet$ we necessarily have
\be\label{eq:tildeGI=}
\wt\Gc_I \,\simeq\, (k'_I)^! \, \Gc^\bullet [r-|I|],
\ee
because for $I_1\subsetneq I$ we have $(k'_{I_1})^! k'_{I_*}\wt\Gc_I=0$. 
This means that $\Gc^\bullet$, as an explicit complex of sheaves, is isomorphic
to the intrinsic Cousin complex of the perverse sheaf $\DD(\Fc)\simeq\Gc^\bullet$
but formed using the $X_I^\Re$ instead of $X_I^\Im$. In particular, each
$\wt\Gc_I$ is constructible with respect to $\Sc^\Re\wedge \Sc^{(0)}=\tau(\Sc^{(1)})$.
This means that its generalization maps  $(\del'_{\bm,\bn})^*$ associated to anodyne $\bm\geq'\bn$
are isomorphisms, and so the corresponding $\del'_{\bm,\bn}$   are
themselves
isomorphisms, which gives (MBS3) for $E$. \qed

\vskip .2cm

The proposition means that we have a  functor
\[
\EE: \Perv(\hW)\lra \MBS_{(W,S)}.
\]
It is clear that $\GG\circ \EE\simeq\Id$, as $\Fc$ is quasi-isomorphic to 
$\Ec^\bullet(\Fc)=\GG(\EE(\Fc))$. Conversely, given $E\in\MBS_{(W,S)}$ and 
writing 
$\Fc=\GG(E)=\Ec^\bullet(E)$, we see that $\Ec^\bullet(E)$, as an explicit complex,
is isomorphic to the intrinsic Cousin complex of the perverse sheaf $\Fc$.
This follows from the identification $\Ec_I(E) \,=\, (k''_I)^! \Ec^\bullet(E)[r-|I|]$
obtained  in the same way as \eqref{eq:tildeGI=}. This means  $\EE(\GG(E))\simeq E$
which finishes the proof of  parts (a)-(c)  of Theorem \ref{thm:main-II}.  

\paragraph{Origin of the sign twist.} Let us now prove part (d) of of Theorem \ref{thm:main-II}. 
Denote  by  $k: W\backslash \hen^\reg\hra \hW$  the embedding. 
Let $E\in\MBS_{(W,S)}$. 
As $\Ec(E)\in\Perv(\hW)$, the restriction  $k^* \Ec(E)$ has the
form $\Lc[r]$, where $\Lc$ is the  local system given  explicitly by 
\[
\Lc \,=\, \Ker  \biggl( k^* \Ec_\emptyset(E) \buildrel d\over\to \bigoplus_{|I|=1}k^* \Ec_I(E)\otimes\det(I) \biggr), 
\]
(see \eqref{eq:cousin}). 
Denote by  $ l: X_\emptyset^\Im  \hra W\backslash \hen^\reg$ 
the embedding whose composition
with $k$ is $k''_\emptyset: X_\emptyset^\Im\hra\hW$.  
Now, by definition $\Ec_\emptyset(E) = k''_{\emptyset *}\,  \, \wt\Ec_\emptyset(E)$, 
where $\wt\Ec_\emptyset(E)$ is the local system on the imaginary stratum $X_\emptyset^\Im$
with stalks $E(\bm)$ and generalization maps $\del''_{\bm,\bn}$ for anodyne $\bm\geq''\bn$
such that $U_\bn, U_\bm\subset X_\emptyset^\Im$. In other words,
$\wt\Ec_\emptyset(E)=l^*\Lc_E$, and 
$\Ec_\emptyset(E)= j_{\emptyset *} l^*\Lc_E$. So 
over $X_\emptyset^\Im\subset W\backslash \hen^\reg$,
 we  have 
\[
(k''_\emptyset)^* \Ker(d) \simeq (k''_\emptyset)^* \Ec_\emptyset(E) \simeq 
 \Lc_E
\]
(no sign yet!).   

\vskip .2cm

For two faces $C,D\in\Cc$ we have the orbit $W(C,D)\in\Xi$ and we denote
by $U_{C,D}=U_{W(C,D)}\in\Sc^{(2)}$ the corresponding cell. 
Let $\Sc$ be the cell decomposition of $\hW^\reg$ induced by $\Sc^{(2)}$.
Denote by $U\subset\hW^\reg$ the union of all cells of $\Sc$ of codimension
$\leq 1$, i.e., of the $U_{C,D}\subset\hW^\reg$ such that $\codim(C)+\codim(D)\leq 1$.
Obviously, $U$ is connected.

\begin{lem}
The homomorphism of the fundamental groups $\pi_1(U)\to\pi(\hW^\reg)$
induced by the embedding, is surjective. 
\end{lem}

\noindent{\sl Proof:} Note that $\hW^\reg$ is a $C^\infty$-manifold. Let $\Sc^*$
be the {\em Poincar\'e dual cell decomposition} of $\hW^\reg$, so that $p$-dimensional
cells of $\Sc^*$ are in bijection with codimension $p$ cells of $\Sc$. 
Explicitly, the cells of $\Sc^*$ can be given as some unions of
simplices of the barycentric subdivision of $\Sc$, cf. \cite{basak}. 
The set $U$ is homotopy equivalent to the $1$-skeleton of $\Sc^*$, while
$\hW^\reg$ is identified with the union of all the cells of $\Sc^*$. Therefore
the fundamental groupoid of $\hW^\reg$ with base points at the $0$-cells of $\Sc^*$
has the well known presentation: the generators correspond to $1$ cells
of $\Sc^*$, while the relations correspond to $2$-cells of $\Sc^*$. Passing from groupoids
to groups with a single base point, this implies that the fundamental group
of $U$ (where no relations are imposed) maps surjectively to that
of $\hW^\reg$ (where the relations are  imposed). \qed

\vskip .2cm
By the lemma, in order to
 identify  two
local systems on the entire 
 $W\backslash \hen^\reg$  it suffices to do so over $U$, i.e., 
   over the union of the $U_{C,D}$ with
$\codim(C)+\codim(D)\leq 1$. 

\vskip .2cm 

If $\codim(C)=0$ and $\codim(D)\leq 1$, then $U_{C,D}\in X_\emptyset^\Im$,
so by the above we have an identification $ \Ker(d) = \Lc_E$ over
$U_{C,D}$. 

\vskip .2cm 

Suppose that $\codim(C)=1$ and $\codim(D)=0$. Then $C$ lies in the closure of
exactly two codimension $0$
faces, say $C_1$ and $C_2$. This means that the open set $X_\emptyset^\Im$ approaches
the cell $U_{C,D}$ from two sides, similarly to what is depicted in the
left part of Fig. \ref{fig.sl2-2}. So the local system structure on $\Lc_E$
gives an identification of $\Ec_\emptyset(E) = k''_{\emptyset *} l^*\Lc_E$
{\em with the direct sum  of two copies of $\Lc_E$} over $U_{C,D}$.
The relevant part of the local system structure on $\Lc_E$ is given by the inverses of 
the anodyne $\del'$ corresponding to the inequalities $W(C_i,D) \geq' W(C,D)$, $i=1,2$. 
At the same time, the differential $d$ of which
we take the kernel, is given by these same $\del'$ (not inverses). This means that
sections of $\Ker(d)$ near $U_{C,D}$
 will be pairs of sections of $\Lc_E$  on the two sides of the ``cut'' $U_{C,D}$
whose values on $U_{C,D}$ (with respect to the local system structure on $\Lc_E$)
sum up to $0$. Such pairs can be seen as sections of $\Lc_E\otimes\Lc_\sgn$. 
This finishes the proof of  part (d) of of Theorem \ref{thm:main-II}, so the theorem is
proved.


\section{Geometry of Bruhat orbits}\label{sec:geom-bruhat}\label{sec:GBO}

\paragraph{Parabolic Bruhat decomposition.}\label{par:Par-BD}

We  now apply the results of \S 1-4 to the following situation:

\vskip .2cm

$\gen_\CC\supset\ben_\CC \supset\hen=\hen_\CC$ is  a split reductive Lie algebra over $\CC$, its
chosen Borel and Cartan subalgebras. It is standard that these data are in fact defined over $\ZZ$.
In particular, we have  the real vector space $\hen_\RR$, the real part of $\hen$.

\vskip .2cm

$\hen_\RR^* \supset\Delta\supset \Delta^+\supset \Delta_\simp$ is the space of real weights,
with the subsets of all roots (weights of $\gen$), positive roots (weights of $\ben$)
and simple roots. The set of positive roots defines the dominant cone
$C^+\subset\hen_\RR$. 

\vskip .2cm

$W$ is  the Weyl group acting on $\hen$. For $\alpha\in\Delta$ we denote by 
$s_\alpha\in W$ the corresponding reflection. Thus the set $S$ of simple reflections
of $W$ is identified with $\Delta_\simp$. 

 \vskip .2cm
 
$
\Hc \,=\, \bigl\{ \hen_\RR^\alpha = (\alpha^\perp)_\RR, \, \alpha\in\Delta^+\bigr\}
$ is 
the arrangement of root hyperplanes in $\hen_\RR$ i.e., of the reflection hyperplanes
for $W$. 

 \vskip .2cm

 Let $\KK$ be a field. Since $\gen_\CC$ is in fact defined over $\ZZ$,
 we have the corresponding split reductive Lie $\KK$-algebra $\gen_\KK$,
 and similarly for $\hen_\KK\subset\ben_\KK\subset\gen_\KK$. The root system
 $\Delta\supset\Delta^+\supset\Delta_\simp$ is then embedded into the $\KK$-vector space $\hen_\KK^*$. 
For $\alpha\in\Delta$ we denote by $e_\alpha\in\gen$ the Chevalley root generator
corresponding to $\alpha$. 

 \vskip .2cm

 Let $G$  be a split reductive algebraic group over $\KK$ with Lie algebra $\gen_\KK$ and 
  $T\subset B\subset G$ the maximal torus and Borel subgroup with Lie algebras
 $\hen$ and $\ben$ respectively. 
   A parabolic subgroup $P\subset G$  (resp. parabolic subalgebra $\pen\subset\gen_\KK$) is called
 {\em standard}, if $P\supset B$ (resp. $\pen\supset\ben_\KK$). 
  As is well known, standard parabolics correspond to
 subsets $I\subset\Dsim$.  We denote 
 \[
 P_I= G_I U_I,    \quad \pen_I = \gen_I\oplus \uen_I, \quad \pen_I=\Lie(P_I), \, \gen_I=\Lie(G_I), 
 \, \uen_I=\Lie(U_I)
 \]
 the standard parabolic subgroup corresponding to $I$ with its standard  Levi subgroup $G_I$ and
 unipotent radical $U_I$, as well as the corresponding standard parabolic subalgebra  $\pen_I$
 with its standars Levi $\gen_I$ and nilpotent radical $\uen_I$. 
 Thus $\pen_I$ is generated by $\ben_\KK$ and the root generators
 $e_{-\alpha}$, $\alpha\in I$. 
 
 \vskip .2cm
 
 A parabolic subgroup $P\subset G$, resp. subalgebra $\pen\subset\gen_\KK$ will be
  called {\em semi-standard}, if $P\supset T$ resp. $\pen\supset\hen$. 
 Again,  the following is well known.
 
 \begin{prop}\label{prop:SSP}
 \begin{enumerate}
  \item[(a)] Semi-standard parabolics are in bijection with faces $C\in\Cc$
 of the Coxeter complex. Given $C\in\Cc$, the corresponding semi-standard parabolic subgroup and subalgebra with its Levi
and uni/nilpotent radical
 \[
 P_C = G_CU_C, \quad \pen_C=\gen_C\oplus \uen_C, \quad \pen_C=\Lie(P_C),\,
 \gen_C=\Lie(G_C),\, \uen_C=\Lie(U_C)
 \]
are characterized by the following conditions: 
\begin{itemize}
 \item[(SSP)]  The roots of $\pen_C$ are those $\alpha\in\Delta$ for which $\alpha|_C\geq 0$. 
 Among these, 
 the roots of $\gen_C$ are the $\alpha$ satisfying $\alpha|_C=0$ and the roots of $\uen_C$
 are the $\alpha$ satisfying $\alpha|_C>0$. 
 \end{itemize}
 
\item[(b)]  Two semi-standard parabolics  are conjugate with respect to $G$, 
 if and only if they are conjugate
with respect to the normalizer
 $N(T)\subset G$, and such conjugation
  corresponds to the action of $W=N(T)/T$ on $\Cc$. \qed
  \end{enumerate} 
 \end{prop}

We denote by $F_I=G/P_I$ the flag space associated to $I\subset\Dsim$.
 We consider it as an algebraic variety over $\KK$. 
 As is well known, $G/P_I$  parametrizes parabolic subgroups $P\subset G$ conjugate to $P_I$ as well as parabolic subalgebras $\pen\subset\gen$ conjugate to $\pen_I$. 
 We refer to such parabolics as {\em parabolics of type} $I$.  If $I_1\subset I_2$, then
 $P_{I_1}\subset P_{I_2}$ so we have the projection
 \be\label{eq:q-I1-I2}
 q_{I_1, I_2}: F_{I_1} \lra F_{I_2}. 
 \ee
 
 \vskip .2cm

 By a {\em Bruhat orbit}  of type $(I,J)$ we will mean a $G$-orbit $O$ on $F_I\times F_J$.
 Such an $O$ is a quasi-projective variety over $\KK$ which we think of as consisting
 of pairs of parabolics $(P,P')$. 
 The parabolic Bruhat decomposition   can be formulated as follows. 

\begin{prop}\label{prop:par-bruhat}
 Let $I,J\subset\Delta_\simp$. We have a bijection
 \[
 G\backslash (F_I\times F_J) \,\simeq \, W\backslash ((W/W_I)\times (W/W_j))
 \,=\,\Xi(I,J)
  \,\subset\, 
 \Xi = W\backslash (\Cc\times \Cc). 
 \]
 More precisely, each $G$-orbit on $F_I\times F_J$ contains a pair of semi-standard parabolics
 $(P_C, P_D)$ for some pair of faces $(C,D)\in\Cc\times \Cc$ defined uniquely up to a simultaneous
 $W$-action. 
 
\end{prop}

\noindent{\sl Proof:} 
The standard formulation in e.g.,  \cite[\S 14.16]{borel}  or 
\cite[Ch. 4, \S 2.5, Rem. 2]{bourbaki}, 
   is in terms of  an identification of the sets of double cosets
\[
P_I\backslash G/P_J \,\simeq W_I\backslash W /W_J. 
\]
To get the statement in our form,  recall that for any group $H$ and subgroups $K,L$
we have a bijection
\[
H\backslash ((H/K)\times (H/L)) \buildrel\simeq\over\lra K\backslash H/L, 
\quad H(h_1K, h_2L) \,\mapsto\,  K (h_1^{-1}h_2) L. 
\]
The remaining details are left to the reader. \qed

\vskip .3cm
 
 Thus the 2-sided Coxeter complex $\Xi$ parametrizes Bruhat orbits in all the $F_I\times F_J$.
 For $\bm\in \Xi(I,J)$ we denote by 
 \[
 F_I\buildrel p'_{\bm}\over  \lla O_\bm\buildrel p''_{\bm}\over  \lra F_J
 \]
the corresponding Bruhat orbit with its projections to the factors. According to the above proposition this diagram may be identified with
\[
 G/P_C\buildrel p'_{\bm}\over  \lla G/(P_C\cap P_D) \buildrel p''_{\bm}\over  \lra G/P_D.
 \]

 \paragraph{Bruhat order on $\Xi(I,J)$.} 
 The identification $\Xi(I,J) \simeq G\backslash (F_I\times F_J)$ makes manifest the
  {\em Bruhat order}  on $\Xi(I,J)$, which we denote
$\teq$. It reflects the relation of inclusion of orbit closures. That is,
$\bm\teq\bn$ iff $O_\bm\subset \ol{O_\bn}$. With respect to the identification
$\Xi(I,J) \simeq W_I\backslash W / W_J$, this order is induced by the
two-sided Bruhat order on $W$.  This latter identification
implies the following. 
 
 \begin{prop}
 The contraction maps \eqref{eq:contr-maps} are monotone with respect to
  the Bruhat orders $\teq$
 in their source and target. 
 \end{prop}
 
 \noindent{\sl Proof:} 
 For example, 
$\phi'_{(I_1, I_2|J)}: \Xi(I_1, J) \to  \Xi(I_2, J)$,  $I_1\subset I_2$, is  the map
\[
W_{I_1}\backslash W/ W_J \lra W_{I_2}\backslash W/W_J
\]
induced by the inclusion $W_{I_1}\subset W_{I_2}$. Since the  orders  $\teq$ on the source
and target of this map are induced by the same Bruhat order on $W$, the map is monotone. 
\qed

\paragraph{Structure of the orbits. }   As for any real hyperplane arrangement,
the set $\Cc$ of faces of $\Hc$ carries the {\em composition}, or {\em Tits product}
operation $\circ$, see \cite[\S 2.30]{tits}, \cite{BZ} or \cite[\S 2B]  {KS-hyp-arr}. 
For two faces $C,D$
the new face $C\circ D$ can be described, geometrically, 
 as follows. Choose any $c\in C, d\in D$ and
draw a straight line interval $[c,d]\subset\hen_\RR$. Then $C\circ D$ is the face containing 
the points $c'\in [c,d]$ which are very close to $c$ but not equal to $c$. 
By construction, $C\leq C\circ D$. The operation $\circ$ is not commutative. 

\vskip .2cm

Let now $\bm\in\Xi(I,J)$. Representing  $\bm$ as an orbit $\bm=W(C,D)$,
we have two $W$-orbits $W(C\circ D)$ and $W(D\circ C)$
associated to $\bm$. Define two subsets $\Hor(\bm), \Ver(\bm)\subset\Dsim$
called the {\em horizontal} and {\em vertical readings} of $\bm$ by the conditions
\[
W(C\circ D) \,\ni \, K_{\Hor(\bm)}, \quad W(D\circ C) \,\ni\, K_{\Ver(\bm)} 
\]
(see 1.A for the notation $K_I$).
Note that $C\leq C\circ D$ implies $\Hor(\bm)\subset I$ and
$D\leq D\circ C$ implies $\Ver(\bm)\subset J$. 

\begin{prop}\label{prop:shadows}
Let $\bm\in\Xi(I,J)$ and $O_\bm\subset F_I\times F_J$ be the corresponding Bruhat orbit. 
\begin{enumerate}

\item[(a)] 
For any pair of parabolic subgroups $(P,P')\in O_\bm$ with unipotent radicals $U,U'$,
the subgroup 
$$
P\circ P' := (P\cap P')U\subset P
$$ 
is a parabolic subgroup in $G$ of type $\Hor(\bm)$,
and $P'\circ P = (P\cap P')U'\subset P'$ is a parabolic subgroup in $G$ of type $\Ver(\bm)$.

\item[(b)] 
 Associating to $(P,P')$  the subgroups $ P\circ P'$ 
 and $P'\circ P$ defines projections $r'_\bm$, $r''_\bm$ in the  commutative
diagram
\[
\xymatrix{
F_{\Hor(\bm)}
\ar[d]_{q_{\Hor(\bm), I}}&& \ar[ll]_{r'_{\bm}} O_\bm
\ar[lld]^{p'_\bm} 
\ar[rrd]_{p''_\bm}
 \ar[rr]^{r''_{\bm}}&&F_{\Ver(\bm)}\ar[d]^{q_{\Ver(\bm), J}}
\\
F_I &&&&F_J. 
}
\]

\item[(c)] The fibers of $r'_\bm, r''_\bm$ are affine spaces. 
\end{enumerate}
\end{prop}

\noindent{\sl Proof:} (a) By Proposition \ref{prop:par-bruhat}, we can assume $P$ and $P'$ 
semi-standard: $P=P_C$, $P'=P_D$ for some faces $C,D\in\Cc$. We claim that
$$
P_C\circ P_D = P_{C\circ D}.
$$ 
It suffices to prove 
the equality of the Lie algebras
\[
(\pen_C\cap \pen_D)\oplus \uen_C \,=\, \pen _{C\circ D}. \quad 
\]

For this, we recall the algebraic definition of $C\circ D$ in  \cite{BZ} or
 \cite [\S 2C] {KS-hyp-arr}.
That is, consider the set $\{0,+,-\}$ with the partial order $0 <  +$ and $0<-$ while  $+$ and $-$ are
non-comparable. For any $\alpha\in\Delta$ and any  $C\in \Cc$ we have the sign
$\sgn(\alpha|_C)\in \{0,+,-\}$. Then 
\[
\sgn(\alpha|_{C\circ D}) \,=\,\begin{cases}
\sgn(\alpha|_D),& \text{ if } \sgn(\alpha|_C)<\sgn(\alpha|_D), 
\\
\sgn(\alpha|_C), & \text{otherwise}. 
\end{cases}
\]
This means :
\[
\begin{gathered}
\alpha|_{C\circ D}=0\quad \Leftrightarrow \quad \alpha|_C=\alpha|_D=0, 
\\
\alpha|_{C\circ D}>0\quad \Leftrightarrow \quad \bigl(  ( \alpha|_D>0, \alpha|_C=0) \,\,\text{or}\,\,
\alpha|_C>0 \bigr), 
\end{gathered} 
\] 
 which, in view of Condition (SSP) of Proposition \ref{prop:SSP}(a), gives precisely the roots of 
 $(\pen_C\cap \pen_D)\oplus \uen_C$.
 
This proves part (a). Part (b) is now clear. To see (c), we agan look at  the  semi-standard
 representatives
above. In this case $O_\bm=G/(P_C\cap P_D)$, as $P_C\cap P_D$ is the
stabilizer in $G$  of the point $(P_C, P_D)\in F_I\times F_J$. Now, the subgroup $P_C\cap P_D\subset G$
may not be parabolic but has the same Levi quotient as  $(P_C\cap P_D)U_C$. Therefore
the fibers of $r'_\bm$ are isomorphic to 
\[
(P_C\cap P_D)U_C\bigl/ (P_C\cap P_D) \,=\,
U_C\bigl/ (U_C\cap P_C\cap P_D) \,=\, U_C\bigl/ (U_C\cap P_D)
\]
which is the factor of a unipotent group by a unipotent subgroup so it is isomorphic to an
affine space. Similarly for $r''_\bm$. 
 \qed
 
 The {\it fundamental diagram} 5.5 (b) may be rewritten as
 
 \[
\xymatrix{
G/P_{C\circ D}
\ar[d]_{q_{\Hor(\bm), I}}&& \ar[ll]_{r'_{\bm}} G/(P_C\cap P_D)
\ar[lld]^{p'_\bm} 
\ar[rrd]_{p''_\bm}
 \ar[rr]^{r''_{\bm}}&& G/P_{D\circ C} \ar[d]^{q_{\Ver(\bm), J}}
\\
G/P_C &&&& G/P_D. 
}
\]
 
 \paragraph{The diagram of Bruhat orbits.}  Let $\bm, \bn\in\Xi$ and $\bm\geq \bn$.
 In particular, if $\bm\in\Xi(I_1, J_1)$ and $\bn\in\Xi(I_2, J_2)$, then $I_1\subset I_2$ and
 $J_1\subset J_2$, so we have the projection
 \[
 q_{I_1, I_2}\times q_{J_1, J_2}: F_{I_1}\times F_{J_1} \lra F_{I_2}\times F_{J_2}. 
 \]
 
 \begin{prop}\label{prop:diag-pmn}
 \begin{enumerate}
 \item[(a)] If $\bm\geq\bn$, then $q_{I_1, I_2}\times q_{J_1, J_2}$ takes $O_\bm$ to $O_\bn$, so we have
 a   projection $p_{\bm, \bn}: O_\bm\to O_\bn$.

\item[(b)] The projections $p_{\bm,\bn}$, $\bm \geq\bn$, are transitive, so they form a contravariant functor
 from $(\Xi,\leq)$ to the category of algebraic varieties over $\KK$. 
 \end{enumerate}
 \end{prop}
 
 \noindent{\sl Proof:} (a) If $\bm\geq\bn$, then we can represent $\bm=W(A,B)$, $\bn=W(C,D)$,
 where $A,B,C,D\in\Cc$ are such that $A\geq C$ and $B\geq D$. 
 This means that $P_A\subset P_C$ and $P_B\subset P_D$. But  the first inclusion means 
 that $P_C$, considered as a point of
 $F_{I_2}$, is the image of $P_A$, considered as a point of $F_{I_1}$, under $q_{I_1, I_2}$.
 Similarly for $P_B$ and $P_D$. This shows that one point of $O_\bm$, namely
 $(P_A, P_B)\in F_{I_1}\times F_{J_1}$ is mapped into a point of $O_\bn$, namely
 $(P_C,P_D)\in F_{I_2}\times F_{J_2}$. Since both $O_\bm$ and $O_\bn$ are 
 $G$-orbits and the projection in question is $G$-equivariant, we conclude that
 $O_\bm$ is mapped onto $O_\bn$, in a surjective way, thus proving (a).
 Now part (b) is obvious because of the transitivity of the projections $q$ in \eqref{eq:q-I1-I2} for any  three
 subsets $I_1\subset I_2\subset I_3$. \qed
  
 \paragraph{Maps of orbits and maps of flag varieties.} 
      
  \begin{prop}\label{prop:Om-On-square}
  \begin{enumerate}
    \item [(a1)] If $\bm\geq' \bn$, then $\Ver(\bm)\subset \Ver(\bn)$, and we have a commutative
     diagram
     \[
 \xymatrix{
 O_\bm\ar[rr]^{p_{\bm,\bn}}
 \ar[d]_{r''_\bm}&& O_\bn\ar[d]^{r''_\bn}
 \\
 F_{\Ver(\bm)} \ar[rr]_{q_{\Ver(\bm), \Ver(\bn)}} && F_{\Ver(\bn)}
 }
  \]
 
 \item[(a2)] If, moreover, $\bm\geq' \bn$ is anodyne, then $\Ver(\bm)=\Ver(\bn)$.

  \item[(b1)] If $\bm\geq '' \bn$, then $\Hor(\bm)\subset\Hor(\bn)$, and we have
 a commutative diagram 
     \[
 \xymatrix{
 O_\bm\ar[rr]^{p_{\bm,\bn}}
 \ar[d]_{r'_\bm}&& O_\bn\ar[d]^{r'_\bn}
 \\
 F_{\Hor(\bm)} \ar[rr]_{q_{\Hor(\bm), \Hor(\bn)}} && F_{\Hor(\bn)}
 }
  \]
  
  \item[(b2)] If, moreover, $\bm\geq''\, \bn$ is anodyne, then $\Hor(\bm)=\Hor(\bn)$. 
  \end{enumerate}
     \end{prop}

 \noindent{\sl Proof:}  It suffices to prove (a1-2), the other two statements being
 similar. Suppose $\bm\geq'\bn$. Then
 we can represent $\bm=W(A,C)$, $\bn=W(B,C)$ with $A,B,C\in \Cc$ such that
 $A\geq B$.  Now, the Tits product $\circ$ is monotone in the second argument, see  \cite[Prop. 2.7(a)] {KS-hyp-arr}.
 Therefore $A\geq B$ implies $C\circ A \geq C\circ B$, which means
 $\Ver(\bm)\subset\Ver(\bn)$. The commutative diagram  follows  directly
 from the definitions of the maps, thus
 proving (a1). 
 
 \vskip .2cm
 
 To prove (a2), introduce the following notation. 
 For any two flats $M,N\subset\hen_\RR$ of $\Hc$ let $M\vee N$ be the minimal flat  of $\Hc$ containing
 them both. 
 Suppose now that $\bm\geq'\bn$ is anodyne.  This means that we can represent 
 $\bm=W(A,C)$, $\bn=W(B,C)$ with $A,B,C\in \Cc$ such that
 $A\geq B$ and the product cells $iA+C, iB+C\subset\hen_\CC$ lie in the same
 stratum of $\Sc ^{(0)}_\Hc$, i.e., in the generic part of the complexification $L_\CC$
 of the same real flat of $\Hc$. This last condition means that
 \[
  \Lin_\RR(A)\vee \Lin_\RR(C)\,= \, L \,=\, \Lin_\RR(B)\vee \Lin_\RR(C). 
 \]
 As before, we have $C\circ A \geq C\circ B$. Recall now that 
  $\Lin_\RR(C\circ A) = \Lin_\RR(C)\vee \Lin_\RR(A)$,
  by  \cite {KS-hyp-arr} Prop. 2.7(b), and similarly for $C\circ B$. 
   This means that $C\circ A$ and
  $C\circ B$ have the same linear envelope and thus $C\circ A=C\circ B$. This implies that
  $\Ver(\bm) = \Ver(\bn)$, proving (a2). \qed

   \paragraph{Maps of  orbits and  correspondences between flag varieties.}    
     For future use, we record a companion result to Proposition \ref{prop:Om-On-square},
     dealing with the other type of projections. Namely, if $\bm\geq'\bn$, then,
     in general, $\Hor(\bm)\not\subset\Hor(\bn)$, so $F_{\Hor(\bm)}$ and $F_{\Hor(\bn)}$
     are not connected by a map. However, they are connected by a correspondence,
     as the following proposition shows.
     
     \begin{prop}\label{prop:Corr-Z}
     \begin{enumerate}
     \item[(a)] If  $\bm\geq' \bn$, then the image under $p_{\bm,\bn}$ of any fiber of
     $r'_\bm$ is a union of fibers of $r'_\bn$. Therefore we have a commutative diagram
     with a Cartesian square
     \[
     \xymatrix{
     O_\bm\ar[rr]^{p_{\bm,\bn}} 
     \ar[dr]^{s}
     \ar[d]_{r'_\bm}
    && O_\bn
     \ar[d]^{r'_\bn}
     \\
     F_{\Hor(\bm)} & \ar[l]^{\hskip 0.5cm \rho_1} Z \ar[r]_{\hskip -0.5cm \rho_2}& F_{\Hor(\bn)},
     }
     \]
     where
      \[
     Z\,=\,   \bigl\{ (x_1, x_2)\in F_{\Hor(\bm)}\times F_{\Hor(\bn)}\bigl| (r'_\bm)^{-1}(x_1) \supset 
     (r'_\bn)^{-1}(x_2) \bigr\}, 
     \]
     $\rho_1$ and $\rho_2$ are the projections to the first and second factor, and
     \[
     s(o) \,=\,  (r'_\bm(o), r'_\bn(p_{\bm,\bn}(o)),\quad o\in O_\bm. 
     \]
     Further, $Z$ is a single $G$-orbit. More precisely, if $\bm= W(A,C)$,
     $\bn=W(B,C)$ with $A,B,C\in\Cc$ and $A\geq B$, then $Z=O_\bx$,
     where $\bx=W(A\circ C, B\circ C)$.

    \item [(b)] Likewise, if $\bm\geq''\bn$, then the image under $p_{\bm,\bn}$ of any fiber of
     $r''_\bm$ is a union of fibers of $r'_\bn$, and we have a commutative diagram
     similar to (a). 
     \end{enumerate}
  
     \end{prop}

     \noindent{\sl Proof:} We prove (a), since (b) is similar. Consider first the following general situation.
     Let $K_1\supset H_1\subset H_2\subset K_2$ be subgroups of $G$, so we have a diagram
     of projections of homogeneous spaces
     \be\label{eq:diagram-groups}
     \xymatrix{
     G/H_1\ar[d]_{r_1}
     \ar[r]^p& G/H_2\ar[d]^{r_2}
     \\
     G/K_1 & G/K_2. 
     }
     \ee
The condition that $p(r_1^{-1}(x_1))$, for any coset  $x_1\in G/K_1$, is a union of fibers $r_2^{-1}(x_2)$,
is, by homogeneity,  equivalent to the condition that it is such a union for  single $x_1$, e.g. for 
 the coset $K_1$. In this case $p(r_1^{-1}(x_1))\subset G/H_2$ is    the set of left cosets by $H_2$
 contained in the right $H_2$-invariant subset $K_1H_2\subset G$. The condition that it is a union
 of some $r_2^{-1}(x_2)$ is then that  the set $K_1H_2$ is a union of left cosets of $K_2$,
 i.e., it is invariant under right multiplication with $K_2$. This can be expressed as $K_1H_2=K_1K_2$. 
 
 \vskip .2cm
 
 We now apply this to our situation as follows. As $\bm\geq'\bn$, we can write 
 $\bm = W(A,C)$ and $\bn=W(B,C)$ for $A,B,C\in\Cc$ with $A\geq B$. Let $P_A, P_B, P_C$
 be the corresponding semi-standard parabolics, with unipotent radicals $U_A,U_B, U_C$.
 Since $A\geq B$, we have $P_A\subset P_B$, and the unipotent radicals are included in
 the opposite direction:  $U_B\subset U_A$. Our original situation is then a particular case of
 \eqref{eq:diagram-groups} corresponding to
 \[
 \begin{gathered}
 K_1=(P_A\cap P_C)U_A \,=\, U_A(P_A\cap P_C), 
 \\
 H_1= P_A\cap P_C, \,\,\, H_2 = P_B\cap P_C, 
 \\
 K_2 = (P_B\cap P_C)U_B = U_B(P_B\cap P_C). 
 \end{gathered}
 \]
 So $K_1H_2=U_A(P_B\cap P_C)$, while
 \[
 K_1K_2 \,=\, \bigl( U_A (P_A\cap P_C)\bigr) \bigl((P_B\cap P_C)U_B
 \bigr) \,=\, U_A(P_B\cap P_C)U_B \,=\,
 U_AU_B(P_B\cap P_C)\,=\, U_A(P_B\cap P_C),
 \]
 which is the same.  This shows the existence of the diagram, in particular,
 of the correspondence $Z$. 
 
 Further, the morphism $s: O_\bm\to Z$ is
 surjective, since $p_{\bm,\bn}$ and $r'_\bn$ are surjective. Therefore $Z$ is
 a single $G$-orbit. To show that it is exactly the orbit $O_\bx$ as claimed,
 it suffices to find the image of the point $(P_A, P_C)$. Now, $p_{\bm,\bn}$
 takes $P_A$ to $P_B$, so $s(P_A, P_C)= (P_{A\circ C}, P_{B\circ C})$,
 whence the statement. 
 \qed

 \paragraph
   {Example: Associated parabolics and intertwiner correspondences.}\label{par:associ}
   For future reference we recall some elementary instances of the above constructions.
   
   \vskip .2cm
   
   Two faces $C,D\in\Cc$ will be called {\em associated}, if  $\Lin_\RR(C)=\Lin_\RR(D)$.
   Denote this latter space $L$ and put $m=\dim(L)=\dim(C)=\dim(D)$. 
 We further call  two associated faces
   $C$ and $D$ {\em adjacent}, if they are separated by an $(m-1)$-dimensional face $\Pi$,
   that is $C>\Pi<D$. For adjacent $C$ and $D$ we put
   \be\label{eq:delta-C-D}
   \Delta(C,D) \,=\,\{\alpha\in\Delta: \, \alpha|_ C>0, \, \alpha|_D<0\}, 
   \quad \delta(C,D)=|\Delta(C,D)|. 
   \ee
   
 \vskip .2cm

   A {\em gallery} joining two associated faces $C,D$  is a sequence of
   $m$-dimensional faces $(C_0=C,C_1,\cdots, C_l=D)$ all lying in $L$ such that
   for each $i=1,\cdots, l$, the faces $C_{i-1}$ and $C_i$ are adjacent: $C_{i-1}>\Pi_i < C_i$,
   $\dim(\Pi_i)=m-1$.  The number $l$ is
   called the {\em length} of the gallery.
   A {\em minimal gallery} is a gallery of minimal possible length. 
   The length of a minimal gallery is called the {\em face distance}
   between associated faces.

   \vskip .2cm

     Two semi-standard parabolics $P_C, P_D$, $C,D\in\Cc$, are called
    associated, resp. adjacent, if $C$ and $D$ are associated, resp. adjacent.
    
    \vskip .2cm
    
 Assume that $C,D$ are associated.    Note that in this case $C\circ D=C$
     and $D\circ C=D$. Let $I(C), I(D)\subset\Dsim$ be the types of $P_C, P_D$.
     Putting $\bm=W(C,D)\in \Xi(I(C), I(D))$, we find that $\Hor(\bm)=I(C)$ and
     $\Ver(\bm)=I(D)$. 
     
     Let us denote for simplicity $F_C=F_{I(C)}$ (the space of parabolics conjugate
     to $P_C$) and similarly $F_D=F_{I(D)}$.
     Note that $\dim(F_C)=\dim(F_D)$, since $P_C$ and $P_D$ have the same Levi.
      Denote $O_{C,D} = O_\bm\subset F_C\times F_D$ the orbit corresponding to $\bm$.      
     Proposition \ref{prop:shadows} shows that in the diagram 
     \be\label{eq:intertwiner}
 F_C \buildrel p'_{C,D}\over  \lla O_{C,D} \buildrel p''_{C,D}\over  \lra F_D
 \ee
the fibers of both projections are affine spaces (of the same dimension).  
 This diagram is the 
 classical {\em intertwiner correspondence} used to define principal
 series intertwiners (and, more generally, intertwiners between parabolically
 induced representations).  
 
 \begin{prop}\label{prop:associated}
 Let $C,D$ be associated faces. Then: 
 \begin{enumerate}
 
  \item [(a)]  The dimension of the fibers of $p'_{C,D}$ and
 $p''_{C,D}$ is equal to $\delta(C,D)$.

\item[(b)] If   $(C_0=C,C_1,\cdots, C_l=D)$ is a minimal gallery
 joining $C$ and $D$, then the correspondence $O_{C,D}$ is the 
 fiber product of the correspondences
 \[
 F_{C_{i-1}} \buildrel p'_{C_{i-1}, C_i}\over\lla O_{C_{i-1}, C_i} 
  \buildrel p''_{C_{i-1}, C_i}\over\lra F_{C_i}, \quad i=1,\cdots, l.
 \]
 In particular, $\delta(C,D)=\sum_{i=1}^l \delta(C_{i-1}, C_i)$.

 \item[(c)]  Consider chambers (faces  of maximal dimension)
 $\wt C, \wt D$ of $\Hc$ such that $\wt C\geq C$ and $\wt D\geq D$. Then, the minimal face distance possible between
 such $\wt C$ and $\wt D$ is $\delta(C,D)$. If $\wt C, \wt D$ are chambers with this
 face distance, then in the following diagram of projections the squares
 are Cartesian:
 \[
 \xymatrix{
 F_{\wt C} \ar[d] &  \ar[l]_{p'_{\wt C, \wt D}} O_{\wt C, \wt D}\ar[d]
  \ar[r]^{p''_{\wt C, \wt D}} & F_{\wt D} 
  \ar[d]
 \\
  F_{C} &  \ar[l]_{p'_{C,  D}} O_{ C,  D} \ar[r]^{p''_{C,  D}} & F_{ D}.
 }
 \]
 \end{enumerate}
 \end{prop}
 
 \noindent {\sl Proof:} This is classical material. Let us only comment on part (c). 
 Note that $F_{\wt C}=F_{\wt D}$ is the full flag space $G/B$, so the orbit $O_{\wt C, \wt D}$
 corresponds to some  element  $w$ of the Weyl group $W$.
  The condition that, say, the left square in the diagram is Cartesian
 is equivalent to the property  that the fibers of $p'_{\wt C, \wt D}$ project isomorphically
 to the fibers of $p'_{C,  D}$. Given that these fibers are affine spaces and given  the surjectivity
 of the projections, such a property
 reduces to the equality of the dimensions of these affine spaces, i.e., to the equality 
 $l(w)=\delta(C,D)$. Here $l(w)$ is the length of $w$, i.e., the face distance between
 $\wt C$ and $\wt D$. \qed
 
 
 \section{Bruhat orbits as a motivic Bruhat cosheaf}\label{sec:BO=MBC}
 
   The diagram $(O_\bm, p_{\bm,\bn})$ can be seen as a cellular cosheaf on
   $(W\backslash\hen_\CC, \Sc^{(2)})$
  with values in the category  of algebraic varieties over $\KK$. 
  This diagram will be the source of several examples of mixed Bruhat sheaves,
  obtained by
  applying various natural constructions (such as, e.g., passing to the spaces of functions). In this section we highlight the geometric properties
  of  $(O_\bm, p_{\bm,\bn})$ which will imply the axioms of a  mixed Bruhat sheaf
  for such constructions. 
 
\paragraph{ An analog of (MBS3).}
 \begin{prop}\label{prop:anodyne-affine}
 If $\bm\geq\bn$ is an anodyne inequality, then the fibers of $p_{\bm,\bn}$ are affine spaces. 
 \end{prop}
 
 \noindent{\sl Proof:} Any anodyne inequality $\geq$ factors into a composition of an anodyne
 $\geq'$ and an anodyne $\geq''$. So it suffices to prove the statement under additional
 assumption that $\bm \geq' \bn$ or $\bm\geq''\bn$.   
  Suppose $\bm\geq'\bn$ is anodyne. Then in the square of Proposition 
  \ref{prop:Om-On-square}(a) the lower horizontal arrow is the identity, and the
  fibers of the vertical arrows are affine spaces
   by Proposition
 \ref{prop:shadows}(c).  This means that taking a point $o\in O_\bm$,
 the Levi quotients of the $G$-stabilizers of $o$ and $p_{\bm,\bn}(o)$ will be the same, so
 the fibers of $p_{\bm,\bn}$ are affine spaces as well. 
  The case when $\bm\geq''\bn$ is anodyne is similar. \qed

  \begin{rem}
    Proposition \ref{prop:anodyne-affine}
 implies that after passing to the category  $\Dc_\Mc$ of Voevodsky motives
  (where $\AAA^1$-homotopy equivalences become isomorphisms, cf. \cite{beil-volog}), we get
  an $\Dc_\Mc$-valued  cosheaf on $\hW$ that is $\Sc^{(0)}$-constructible. 
 \end{rem}

     \paragraph{Fiber products of orbits. An analog of (MBS2).} 
     Let $\bm'\geq' \, \bn' \leq''\, \bn$ be elements of $\Xi$, so we have the projections
     \[
     \xymatrix{
     & O_{\bn'}\ar[d]^{p_{\bn', \bn}}
     \\
     O_{\bm'} \ar[r]_{p_{\bm', \bn}} & O_\bn. 
     }
     \]
 The following property is the geometric analog of the condition (MBS2) for mixed Bruhat sheaves. It is also analogous to Proposition  \ref{prop:sup=FP}.
 
 \begin{prop}\label{prop:MBS2-geom}
 For any  $\bm'\geq' \, \bn' \leq''\, \bn$
 we have the following decomposition into the
 union of orbits:
 \[
 O_{\bm'}\times_{O_\bn} O_{\bn'} \,   \,\simeq \bigsqcup_{ \bm\in\sup(\bm', \bn)} O_\bm. 
 \]
 \end{prop}
 
 \noindent{\sl Proof:}   The assumption  $\bm'\geq'\,  \bn' \leq' \bn$ implies that there
 are $I_1\subset I_2$ and $J_1\subset J_2$ such that 
 \[
 \bm' \in\Xi(I_1, J_2), \quad \bn'\in \Xi(I_2, J_1), \quad \bn\in\Xi(I_2, J_2). 
 \]
 Any $\bm\in\sup(\bm',\bn)$   must then lie in $\Xi(I_1, J_1)$. 
  Note that 
 \[
 O_{\bm'}\times_{O_{\bn}} O_{\bn'} = O_{\bm'}\times_{F_{I_2}\times F_{J_2}} O_{\bn'}.
 \] 
Note further that the square
\[
\xymatrix{
F_{I_1}\times F_{J_1}\ar[d] \ar[r]& F_{I_2}\times F_{J_1}\ar[d]
\\
F_{I_1}\times F_{J_2} \ar[r]& F_{I_2}\times F_{J_2}
}
\] 
 is Cartesian, being the external Cartesian product of two arrows 
 \[
 \bigl\{ F_{I_1}\buildrel q_{I_1, I_2}\over\lra F_{I_2} \bigr\} \times 
\bigl\{ F_{J_1}\buildrel q_{J_1, J_2}\over\lra F_{J_2}\bigr\}. 
 \]
 Therefore $ O_{\bm'}\times_{O_{\bn}} O_{\bn'}$ is contained in $F_{I_1}\times F_{J_1}$ and
 is the union of those orbits $O_\bm$ that project to $O_{\bm'}$ and $O_{\bn'}$, i.e.,
 of the $O_\bm$ with $\bm\in\sup(\bm',\bn)$, 
  as claimed. \qed

 
 \section{Functions on $\FF_q$-points.}\label{sec:Fq}
 
 \paragraph{Appearance as a bicube.} \label{par:F_q-bic}
 In this section we present the simplest example of
 a mixed Bruhat sheaf encoding the algebra behind parabolic induction and restriction. 
 Like other examples, it appears most immediately (and is well known)  in the form of a bicube,
 cf. \S \ref{sec:perv=mbs}\ref{par:bicube}
 
 \vskip .2cm
 
  We specialize the situation of \S \ref{sec:geom-bruhat} to the case when $\KK=\FF_q$ is a finite field. 
  This field has to be distinguished from
  the ``coefficient'' field $\k$ as in \S  \ref{sec:perv=mbs}\ref {par:MBS}
  In this section we assume that $\k$ is algebraically closed of characteristic
  $0$.  
  For a variety $X/\FF_q$ we denote $\Fun(X)$ the $\k$-vector space of functions $X(\FF_q)\to\k$. 
 If $f: X\to Y$ is a morphism of varieties over $\FF_q$, we denote
 $f^*: \Fun(Y)\to\Fun(X)$, $f_*: \Fun(X)\to\Fun(Y)$ the inverse image (pullback)
 and the direct image (sum over the fibers) of functions on $\FF_q$-points.
 
 \vskip .2cm
 
 For $I\subset\Dsim$ we have the flag space $F_I=G/P_I$.  For $I\subset J\subset\Dsim$ we have
 the projection $q_{IJ}: F_I\to F_J$, cf. \eqref{eq:q-I1-I2}. We define a $\Dsim$-bicube $Q$ by
 \[
 Q_I = \Fun(F_I), \,\, v_{IJ} = (q_{IJ})_*: \Fun(F_I)\to\Fun(F_J), \,\, u_{IJ} = q_{IJ}^*: \Fun(F_J)\to\Fun(F_I). 
 \]
 Thus 
 $
 Q_I \,=\, \Ind_{P_I(\FF_q)}^{G(\FF_q)}\,\k 
 $
 is the simplest parabolically  induced representation.
 We now proceed to upgrade this bicube to a mixed Bruhat sheaf. 
 
 \paragraph{Definition of the diagram $E_q$.}\label{par:E_q}

 \vskip .2cm
 
 For $\bm\in\Xi$ we have the orbit $O_\bm$ and define the $\k$-vector space $\wt E(\bm)=\Fun(O_\bm)$.
 If $\bm\geq'\bn$, we define
 \[
 \del'_{\bm,\bn} =(p_{\bm,\bn})_*: \wt E(\bm) = \Fun(O_\bm) \lra \Fun(O_\bn) = \wt E(\bn).
 \]
 If $\bm\geq''\bn$, we define 
  \[
 \del''_{\bm,\bn} =p^*_{\bm,\bn}: \wt E(\bn) = \Fun(O_\bn) \lra \Fun(O_\bm) = \wt E(\bm). 
 \]
 
 \begin{prop}
 The diagram $\wt E=(\wt E(\bm), \del'_{\bm,\bn}, \del''_{\bm,\bn})$ satisfies the conditions
 (MBS1-2) of Definition \ref {def:MBS}. 
 \end{prop}
 
 \noindent{\sl Proof:} (MBS1), i.e., the transitivity of the $\del'$ and the $\del''$,  follows from
 the transitivity of the $p_{\bm,\bn}$ and from  the compatibility with the direct and inverse images
 with  composition. The condition (MBS2) follows directly from Proposition 
 \ref{prop:MBS2-geom}. More precisely, that proposition implies that 
  for any  $\bm'\geq' \, \bn' \leq''\, \bn$
 we have a Cartesian square
 of finite sets
 \[
 \xymatrix{
  \bigsqcup\limits_{ \bm\, \in\, \sup(\bm',\bn)} O_\bm(\FF_q)
  \ar[d]
   \ar[r]& 
   O_{\bn'}(\FF_q) \ar[d]^{p_{\bn', \bn}}
     \\
     O_{\bm'}(\FF_q) \ar[r]_{p_{\bm', \bn}} & O_\bn(\FF_q). 
 }
 \]
 So $\del''_{\bn', \bn} \del'_{\bm', \bn'} = p_{\bn', \bn}^* (p_{\bm',\bn})_*$ is equal, by the
 base change formula, to the result of first pulling back to the disjoint union of the $O_\bm(\FF_q)$
 and then pushing forward to $O_{\bn'}(\FF_q)$, which is precisely the right hand side of (MBS2). 
 \qed
 
 \vskip .2cm
 
 However, the diagram $\wt E$ does not satify (MBS3). So for each $\bm$  we define the subspace
 \[
 E_q(\bm) \,=\, (r'_\bm)^*\,  \Fun(F_{\Hor(\bm)}) \,\subset \, \wt E(\bm), \quad
 r'_\bm: O_\bm\lra F_{\Hor(\bm)},
 \]
 to consist of functions pulled back from $F_{\Hor(\bm)}$. 
 
 \begin{thm}\label{thm:Fq}
 \begin{enumerate}
 
 \item[(a)] The maps $\del'$ and $\del''$ of $\wt E$ preserve the subspaces $E_q(\bm)$.

\item[(b)] The diagram $E_q=(E_q(\bm), \del'_{\bm,\bn}, \del''_{\bm, \bn})$ is
 a mixed Bruhat sheaf, i.e., it satisfies all three conditions (MBS1-3). 
 
 \item[(c)] The   bicube $Q$ above is obtained from $E_q$ by procedure
 (functor $\Qc$) described
 in \S \ref{sec:perv=mbs}\ref{par:bicube}
 \end{enumerate}
 \end{thm}
 
 \paragraph{ Proof of Theorem \ref{thm:Fq}.} 
 We start with part (a). 
 Look first at the map $ \del''_{\bm,\bn} =p^*_{\bm,\bn}$, $\bm\geq'' \bn$. 
 In this case we have a commutative square of Proposition \ref{prop:Om-On-square}(b1)
 so  $p_{\bm,\bn}^*$ takes functions pulled back by $r'_\bn$ to
 functions pulled back by $r'_\bm$. 
 
 \vskip .2cm
 
 Look now at the map $\del'_{\bm,\bn} = (p_{\bm,\bn})_*$, $\bm\geq' \bn$. 
 We then have the diagram of Proposition \ref{prop:Corr-Z}(a). 
Using the base change for the Cartesian square in that diagram  we get, for any 
$f\in  \Fun(F_{\Hor(\bm)})$:
\[
(p_{\bm,\bn})_* \, (r'_\bm)^*\, f \,=\, (p_{\bm,\bn})_* \, s^* \, \rho_1^* \, f \,=\,(r'_\bn)^* \,
(\rho_2)_* \, (\rho_1^*) \, f, 
\]
so $\del'_{\bm,\bn} = (p_{\bm,\bn})_* $ takes pulled back functions to pulled back functions,
i.e., $E_q(\bm)$ to $E_q(\bn)$. This finishes the proof of part (a). 

\vskip .2cm

We now prove (b). The conditions (MBS1-2) for $E$ follow from the validity of these conditions for $\wt E$. 
Let us prove   (MBS3). 

\vskip .2cm

Let $\bm\geq''\bn$ be anodyne. Then, by Proposition  \ref{prop:Om-On-square}(b2),
$\Hor(\bm)=\Hor(\bn)$, so the square in part (b1) of that proposition becomes a triangle. 
This shows that $\del''_{\bm,\bn}= (p_{\bm, \bn})^*$ is an isomorphism  
from $E_q(\bn)= (r'_\bn)^* \Fun(F_{\Hor(\bn)})$  
 onto $E_q(\bm)= (r'_\bm)^* \Fun(F_{\Hor(\bm)})$. 
 
 \vskip .2cm
 
 Let now $\bm\geq ' \bn$ be anodyne. Proposition \ref{prop:Corr-Z}(a) shows that
 the map $\del'_{\bm,\bn}: E_q(\bm)\to E_q(\bn)$ is isomorphic to the map
 $\rho_{2*}\, \rho_1^*: \Fun(F_{\Hor(\bm)}) \to \Fun(F_{\Hor(\bn)})$. 
  Now, if $\bm=W(A,C)$ and $\bn=W(B,C)$ with $A\geq B$, then the
  condition that  $\bm\geq ' \bn$ is  anodyne means that, similarly to the proof of
  Proposition \ref{prop:Om-On-square}, we have
  \[
 \Lin_\RR(A\circ C) \,=\,  \Lin_\RR(A)\vee \Lin_\RR(C) \,=\,\Lin_\RR(B)\vee\Lin_\RR(C)
 \,=\,\Lin_\RR(B\circ C). 
  \]
 In other words, $A\circ C$ and $B\circ C$ are associated faces, cf. \S
 \ref{sec:GBO}\ref{par:associ} The last claim in Proposition \ref{prop:Corr-Z}(a)
 implies then that $Z$ is a particular case of the intertwiner correspondence
 \eqref{eq:intertwiner} for two associated parabolics. So the isomorphicity of
  $\del'_{\bm,\bn}$ in this case is a particular case of the following classical fact.
  
  \begin{prop}\label{prop:inter-iso}
  Let $P_C, P_D$ be two associated semi-standard parabolics. Then the
  intertwiner
  \[
  (p''_{C,D})_* \, (p'_{C,D})^*: \Fun(F_C) \lra \Fun(F_D)
  \]
  is an isomorphism. 
  \end{prop}
  
  \noindent{\sl Proof:} For convenience of the reader we recall the argument
  by reduction to the simplest case. 
   First, the   Cartesian
  squares in the diagram
  in part (c) of Proposition \ref{prop:associated}  show that the intertwiner for $\wt C, \wt D$ takes functions pulled back
   from $F_C$ to functions pulled back from $F_D$, so it is enough to prove
    that such a Borel intertwiner is an isomorphism. 
    Next, in this case $O_{\wt C, \wt D}$
   corresponds to some element $w\in W$, and by  \ref{prop:associated}(b)
   it is enough to consider the case when $w=s_\alpha$ is a simple reflection.
   In this case we have the $\PP^1$-fibration
   $
   q_\alpha: G/B \to G/P_\alpha
   $
   and the corresponding orbit consists of $(x,y)\in G/B\times G/B$
   such that $q_\alpha(x)=q_\alpha(y)$ but $x\neq y$. That the
   intertwiner is an isomorphism in this case reduces to the case of a single fiber of $q_\alpha$,
   i.e., to the case of the correspondence
   \[
   \PP^1\lla (\PP^1\times \PP^1)\- \text{diag.} \lra \PP^1, 
   \]
   where the isomorphism is obvious. \qed
   
   \vskip .2cm
   
   This finishes the proof of part (b) of Theorem \ref{thm:Fq}. To prove (c), we use
   the notation of  \S \ref{sec:perv=mbs}\ref{par:bicube} In that notation, we observe that
   \[
   \Hor(\bm'_I) =\Hor(\bm_I) = \Hor(\bm''_I) = I
   \]
   and, moreover,
   \[
   O(\bm'_I) = O(\bm_I) = O(\bm''_I) = F_I,
   \]
   the corresponding projections $O\to F_\Hor$ being the identities. 
   More precisely $O(\bm'_I)$ is the space $\pt\times F_I$ (a single orbit),
    and $O(\bm''_I)$ is $F_I\times\pt$ (likewise a single orbit) while $O(\bm_I)$
    is the diagonal orbit in $F_I\times F_I$. 
    We colnclude that the map $\phi_I$ is also the identity. 
    From this, we identiy  $\Qc(E)_I= \Fun(F_I) = Q_I$. With respect to these
    identifications, 
    the map $v_{IJ}$ in $\Qc(I)$, being $\del_{\bm'_I, \bm'_J}$, is 
    identified with $(q_{IJ})_*$  and $u_{IJ}$, being $\phi_I^{-1}\circ \del''_{\bm''_I, \bm''_J}
    \circ \phi_J$, is identified with $q_{IJ}^*$. 
    
     Theorem \ref{thm:Fq}. is proved.\qed

   \paragraph{The case of  ``$\FF_1$-points''.}
    Let us indicate an even simpler example corresponding
   to the formal limit $q=1$ where, instead of groups of $\FF_q$-points, we consider Weyl
   groups.  
   
   For $\bm\in\Xi(I,J)$ put $E_1(\bm)=\Fun(\bm)$, where we consider $\bm$ as a subset
   (orbit) in $(W/W_J)\times (W/W_J)$. If $\bm\geq\bn$, we have a $W$-equivariant surjection
   $\pi_{\bm,\bn}:\bm\to\bn$, as in \eqref{eq:pi-m-n}. 
    If $\bm\geq'\bn$, we define
 \[
 \del'_{\bm,\bn} =(\pi_{\bm,\bn})_*:  E_1(\bm) = \Fun(\bm) \lra \Fun(\bn) = E_1(\bn).
 \]
 If $\bm\geq''\bn$, we define 
  \[
 \del''_{\bm,\bn} =\pi^*_{\bm,\bn}: E_1(\bn) = \Fun(\bn) \lra \Fun(\bm) = E_1(\bm). 
 \]
 
    \begin{prop}
    The diagram $E_1=(E_1(\bm), \del', \del'')$ is a mixed Bruhat sheaf of type $\gen$. 
    \end{prop}
    
    \noindent{\sl Proof:} (MBS1) follows from the transitivity of the $\pi_{\bm,\bn}$.
    The condition (MBS2) follows from base change and Proposition \ref{prop:sup=FP}. 
    Finally, (MBS3) follows from Proposition \ref{prop:ano=bij}. \qed
    
    \vskip .2cm
    
    The mixed Bruhat sheaf $E_1$ consists of $W$-modules, so it gives a perverse sheaf
    $\Fc_1\in\Perv(\hW)$ with $W$-action. In particular, for any irreducible  $W$-module $V$
    we have the mixed Bruhat sheaf and perverse sheaf formed by the vector spaces of
    multiplicities of $V$:
    \be\label{eq:E_1^V}
    E_1^V = (E_1\otimes_\k V)^W \,=\, \bigl(E_1^V(\bm)= (E_1(\bm)\otimes_\k V)^W
    \bigr), \quad \Fc_1^V = (\Fc_1\otimes_\k V)^W. 
    \ee
    If we choose a representative $(C,D)\in\bm\subset \Cc\times\Cc$, then we have the
    ``parabolic'' subgroup $W^{C,D}\subset W$, the stabilizer of the pair $(C,D)$.
    It is conjugate to the ``standard'' parabolic subgroups $W_{\Hor(\bm)}$ and $W_{\Ver(\bm)}$. 
    The choice of $(C,D)$ allows us to identify
   \be\label{eq:E_1^V=inv}
        E_1^V(\bm) \,\simeq \, V^{W^{C,D}}.
    \ee
     Thus we can say that  the mixed Bruhat sheaf $E_1^V$ is formed by the spaces of invariants in $V$
     with respect to the various parabolic subgroups in $W$. 
     
     To identify $\Fc_1$ and $\Fc_1^V$, consider the diagram of projections and open embeddings
     \be\label{eq:p-and-p-reg}
     \xymatrix{
     \hen^\reg \ar[d]_{\wt j}
     \ar[r]^{\hskip -0.3cm p^\reg}& \hW^\reg
     \ar[d]^j
     \\
     \hen\ar[r]_{\hskip -0.2cm p}& \hW.
     }
     \ee
     For any $W$-module $V$ we denote by $\Lc_V$ the corresponding local system on $\hW^\reg$. 
     As before, let $r=\dim_\CC\hen$. 
     
     \begin{prop}
     \begin{enumerate}
    \item [(a)] The perverse sheaf $\Fc_1$ is identified with $p_*\ul \k_\hen[r]$, with the $W$-action  being
     the natural $W$-action on the direct image twisted by
     the sign character. In particular, $\Fc_1$ reduces to a single sheaf in degree $-r$.

   \item  [(b)] For any irreducible $W$-module $V$ we have an identification
     \[
     \Fc_1^V \,\simeq \, (R^0 j_*  \Lc_{V\otimes\sgn})[r]
      \,=\, 
     j_{!*} (\Lc_{V\otimes\sgn}[r]). 
     \]
     In particular, $\Fc_1^V$ reduces to a single sheaf in degree $-r$. 
     \end{enumerate}
      \end{prop}
      
      Here, for a local system $\Nc$ on $\hW^\reg$, the notation $j_{!*}(\Nc[r])$ means
      the image, in the  abelian category $\Perv(\hW)$, of the natural map
      $j_!\Nc [r]\to Rj_*\Nc[r]$ (we note that since $j$ is an affine embedding,
      both the source and the target of that map are perverse sheaves). 
      In our case  $\Nc=\Lc_{V\otimes\sgn}$  this image coincides
      with   $(R^0j_* \Nc)[r]$.

      \vskip .2cm
      
      \noindent{\sl Proof:} (a) We start with general remarks. By a $d$-cell we mean a topological
      space  $\sigma$ homeomorphic to $\RR^d$. For such $\sigma$ its $1$-dimensional
      orientation $\k$-vector space is defined as $\orr(\sigma)=H^d_c(\sigma,\k)$.
      If $L$ is an $\RR$-vector subspace in a finite-dimensional $\RR$-vector space $V$, then
      we have a canonical isomorphism $\orr(V/L)\otimes \orr(L) \simeq \orr(V)$. 
      For $C\in\Cc$ let, as before,  $\Lin_\RR(C)\subset \hen_\RR$ be the $\RR$-linear span of $C$.
      We have a canonical isomorphism $\orr(C)\simeq \orr(\Lin_\RR(C))$. We also define
      the  coorientation space of $C$ as $\coor(C)=\orr(\hen_\RR/\Lin_\RR(C))$, the
      orientation space of the normal space. Thus $\coor(C)\otimes\orr(C)\simeq \orr(\hen_\RR)$.
      As a $W$-module, $\orr(\hen_\RR)\simeq\sgn$ is  the $1$-dimensional sign representation.

      Now, since $\ul\k_\hen[r]$ is perverse, $p_* \ul\k_\hen[r]$ is perverse and
      lies in $\Perv(\hW)$. Let $E'$ be the mixed Bruhat sheaf associated to it.
      Let us construct an isomorphism $E'\simeq E_1\otimes\sgn$ of mixed Bruhat sheaves with
      $W$-action.

       By definition,  the $E'(\bm)$ are the stalks
      of the terms of the Cousin complex of $p_* \ul\k_\hen[r]$, i.e., of the sheaves
      \be\label{eq:cousin-p_*k}
      j_{I*} j_I^! p_*( \ul\k_\hen) \otimes \det(I) [|I|-r], 
      \ee
      where $j_I: X_I^\Im\hra\hW$ is the embedding
      of the imaginary stratum. 
      These sheaves can be found ''upstairs'' in $\hen$, as we are dealing with a direct image  from 
      $\hen$.       The preimage  $p^{-1}(X_I^\Im)$ is the union of the tube cells 
      $\hen_\RR+iC\buildrel
      j_C\over \hra \hen$ over all faces $C\in\Cc$ in the $W$-orbit of $K_I$.  
   For any such tube cell, $j_C^!\ul\k_\hen$  is, by the local Poincar\'e duality,  
   the constant sheaf in degree $\codim_\RR(C)$
   with stalk being $\coor(C)$.  
   
   Note further that $\det(I)$ is identified with $\orr(K_I)$. For each $C\in \Cc$ in the $W$-orbit of
   $K_I$, the projection $p$ identifies $C$ with $K_I$ and so identifies $\orr(C)$ with
   $\orr(K_I)=\det(I)$.

  To descend back to $\hW$, let $\bm\in\Xi(I,J)$, which we think of as a subset (orbit) in
    $\Cc\times\Cc$.
   Then by the above, $E'(\bm)$, which is  the stalk of  \eqref{eq:cousin-p_*k} at the cell $U_\bm$,
   is   identified, $W$-equivariantly as
   \[
 E'(\bm) \,\simeq \,   \bigoplus_{(C,D)\in\bm} \coor(C)\otimes\det(I)
 \,\simeq \,  \bigoplus_{(C,D)\in\bm} \coor(C)\otimes\orr(C) \,\simeq \, 
  \bigoplus_{(C,D)\in\bm} \orr(\hen_\RR). 
 \]
 The latter space is nothing but $\Fun(\bm)\otimes\sgn \,=\, E_1(\bm)\otimes\sgn$. 
  It remains to show that the identifications $E'(\bm)\to E(\bm)\otimes\sgn$
  thus constructed, are compatible with the maps $\del', \del''$. This is straightforward. 
   
   \vskip .2cm 
   
   (b)
   Obviously,
   \[
   \ul\k_\hen \,=\, R^0 \,\wt j_* \,\,   \wt j^* \,\, \ul\k_\hen \ \,=\, \wt j_{!*}\, \,  \wt j^*\,\, \ul\k_\hen. 
    \]
    This implies that
    \[
    p_* \ul\k_\hen \,=\, R^0 j_* j^* p_*\ul\k_\hen \,=\, j_{!*}j^* p_*\ul\k_\hen, 
    \]
    and so the same relation will hold after we take the space of multiplicities of
    any irreducible $W$-module $V$. Now, 
     from (a) it follows that $j^*\Fc_1$ is the local system corresponding to the regular representation
     of $W$ but with the ``external'' $W$-action twisted by sign. Therefore
     $j^* \Fc_1^V \simeq \Lc_{V\otimes\sgn}$, whence the statement.
     \qed
     
     \paragraph{The perverse sheaf $\Fc_q$. The Hecke algebra picture.} 
   Let us use the notation $F=G/B$ for  the full flag space of $G$ and 
    $G_q \supset B_q$ for the finite groups $G(\FF_q) \supset B(\FF_q)$. 
    
    \vskip .2cm
    
    Return  to the mixed Bruhat sheaf $E_q$ from \S \ref{par:E_q}
    Let $\Fc_q\in\Perv(\hW)$ be the corresponding perverse sheaf.
    As $E_q$ and $\Fc_q$ consist of 
      $G_q$-modules,  for any $G_q$-module $V$ we have the mixed Bruhat sheaf
      $E_q^V$ and the perverse sheaf $\Fc_q^V$ formed by the multiplicities of $V$, as in 
      \eqref{eq:E_1^V}. As in \eqref{eq:E_1^V=inv}, we can say that $E_q^V$ ``consists of''
     spaces of  invariants in $V$ under various parabolic subgroups in $G_q$. These subgroups are,
      however, not necessarily the standard ones $P_I(\FF_q)$ so certain
      conjugations are involved. 
      
      \vskip .2cm
      
      Call a $G_q$-module $V$ {\em special}, if it appears in $\Fun(F)=\Ind_{B_q}^{G_q} \k$. 
      It is clear that $E_q^V$ and $\Fc_q^V$ are nonzero only if $V$ is special. 
      Let  $H_q = H(G_q, B_q)\subset \k[G_q]$  be the  {\em Hecke algebra} 
      formed by $B_q$-bi-invariant functions on $G_q$. 
       Let $\Sen$ be the set of irreducible special representations of $G_q$. 
   As in the case of  any finite group and subgroup, it is classical
      that  $\Sen$ in bijection with the set of  irreducible $H_q$-modules.
      More precisely,  we have
      the decomposition
      \be\label{eq:decom-Fun-F}
      \Fun(F) \,=\, \bigoplus_{V\in\Sen} V\otimes R(V), 
      \ee
      where  $R(V)$ is the irreducible
      $H_q$-module corresponding to $V$. It is also classical \cite{iwahori} that $H_q$
      can be given by generators $\sigma_\alpha$, $\alpha\in\Dsim$ subject to the
      braid relations and the quadratic relations
      \[
      (\sigma_\alpha+1)(\sigma_\alpha-q)\,=\,0.
      \] 
      In particular, we have a   morphism of algebras
      \[
      \k[\Br_\gen] \lra H_q, 
      \]
   and so any $H_q$-module $R$ gives a local system $\Lc_R$ on $\hW^\reg$. 
   
   \begin{prop}\label{prop:7-10}
   Let $V\in\Sen$ be an irreducible special representation of $G_q$. 
   The perverse sheaf $\Fc_q^V$ is isomorphic to $j_{!*}(\Lc_{R(V)}\otimes\Lc_\sgn[r])$,
   with $j$ as in   \eqref{eq:p-and-p-reg}.
   \end{prop}
   
   \noindent {\sl Proof:} 
   Consider first the local system of $G_q$-modules $\Lc_{E_q}$ on $\hW^\reg$  
   associated to
   the mixed Bruhat sheaf $E_q$ by Proposition \ref{prop:E-LS}. Its stalk
   at the cell $U_{W(C^+, 0)}$ is $\Fun(F)$, because the Bruhat  orbit associated to
   $W(C^+,0)$ is $F$. Further, the action of the braid group generators
   $\sigma_\alpha$ on that stalk of $\Lc_{E_q}$  is given by the  standard  intertwiners 
   $q_\alpha^*{q_{\alpha*}}-1$, cf. the proof of Proposition \ref{prop:inter-iso}. 
   These operators give precisely the action of $H_q$ on $\Fun(F)$ giving
   the decomposition \eqref{eq:decom-Fun-F}. This means that
   \[
   \Lc_{E_q} \,\simeq \, \bigoplus_{V\in\Sen} V\otimes \Lc_{R(V)}
   \]
   as a local system of $G_q$-modules. By Theorem \ref{thm:main}(b)
   this means that 
 \[
   j^*\Fc_q \, \simeq \, \bigoplus_{V\in\Sen} V\otimes \Lc_{R(V)}\otimes \Lc_\sgn [r]. 
   \]
 as a local system of
   $G_q$-modules. Therefore for $V\in\Sen$, 
   \[
   j^*\Fc_q^V \,\simeq \, \Lc_{R(V)}\otimes \Lc_\sgn [r], 
   \]
   a (shifted)  irreducible local system. Therefore $j_{!*} (\Lc_{R(V)}\otimes \Lc_\sgn [r])$,
   an irreducible perverse sheaf, is contained in $\Fc_q^V$. It remains to show
   that $\Fc_q^V$ is an irreducible perverse sheaf. For this it is enough to show that
   $E_q^V$ is an irreducible mixed Bruhat sheaf, which we do now. 
   We  start with a general lemma referring to the situation
    of \S \ref{sec:perv=mbs}\ref{par:bicube}
   
   \begin{lem}
   Let $E$ be a mixed Bruhat sheaf of type  $(W,S)$ and $Q=\Qc(E)$ is the $S$-bicube
   associated to $E$. If $Q$ is irreducible as an object of the category $\Bic_S$,
   then $E$ is irreducible as an object of the  category $\MBS_{(W,S)}$. 
   \end{lem}
   
   \noindent{\sl Proof of the lemma:} Mixed Bruhat sheaves of type $(W,S)$ are
   the same  as  covariant functor $\Mc_{(W,S)}\to\Vect_\k$, where $\Mc_{(W,S)}$ is
   a natural universal $\k$-linear category. More precisely, $\Mc_{(W,S)}$
    has   objects $[\bm]$ corresponding to  $\bm\in\Xi$, 
   generating morphisms
   $\den'_{\bm, \bn}: [\bm]\to [\bn]$ for $\bm\geq'\bn$ and $\den''_{\bm,\bn}: [\bn]\to[\bm]$
   for $\bm\geq''\bn$, which are subject to the relations expressing (MBS1)-(MBS3). 
   In particular, these relations include the  invertibilty of the morphisms
    $\den'$ and $\den''$ associated
   to anodyne inclusions. 
   
   Similarly,  $S$-bicubes are the same as covariant functors $\Bc_S\to\Vect_\k$,
   where $\Bc_S$ is the $\k$-linear category with fornal objects $[I]$, $I\subset S$
   and generating morphisms $\ven_{IJ}: [I]\to [J]$ and $\uen_{IJ}: [J]\to [I]$
   for $I\subset J$  which are subject to the relations expressing the axioms of a bicube. 
   
   From this point of view, the  functor $\Qc$ associating to a mixed Bruhat sheaf $E$ the bicube $\Qc(E)$
   is given by pullback via a natural functor $\beta: \Bc_S \to \Mc_{(W,S)}$
   whose construction was essentially described in \S \ref{sec:perv=mbs}\ref{par:bicube}
   
   Suppose now that $E$ is reducible in $\MBS_{(W,S)}$, i.e., it has a 
      mixed Bruhat subsheaf $E'\subset E$ with  $0\neq E'(\bm)\subsetneq E(\bm)$
      for some $\bm$. Note that  the object $[\bm]\in\Mc_{(W,S)}$
 can be connected, by a chain of
   isomorphisms (anodyne $(\den')^{\pm 1}$, $(\den'')^{\pm 1}$)   to an object $[\bn]$
    where $\bn$
   is of the form $W(K_I,0)$. Denoting $\fen: [\bm]\to[\bn]$ the composite
   isomorphism of that chain, we have  that $E'(\bn)$ is the result of applying $\fen$ 
   to $E'(\bm)$. 
   This implies that $0\neq E'(\bn)\subsetneq E(\bn)$ as well. 
   On the other hand, such $[\bn]$ is the image $\beta ([I])$. This implies
   that $0\neq \Qc(E')_I \subsetneq \Qc(E)_I$, i.e., $\Qc(E)$ is reducible
   in $\Bic_S$. Lemma is proved. \qed
   
   \vskip .2cm
   
   Applying the lemma to our situation, we invoke Theorem  \ref{thm:Fq} (c)
   which implies that $\Qc(E_q^V)= Q^V$, the 
   bicube of $V$-multiplicities in the bicube $Q$ of $G_q$-modules
   from  \S \ref{par:F_q-bic} So it is enough to show that 
    $Q^V$  
   is irreducible in the category of bicubes. 
   
   Now, $Q_\emptyset=\Fun(F)$, and for each $I\subset\Dsim$ the structure map
   $v_{\emptyset, I}: Q_\emptyset\to Q_I$ (direct image of functions) is surjective,
   while $u_{\emptyset, I}: Q_I\to Q_\emptyset$ (inverse image of functions)
    is injective. These properties will  still be true for the cube $Q^V$,
    since taking invariants under a finite group is an exact functor.
     So if we have a nonzero element $f$ in some $Q_I^V$, then 
     $u_{\emptyset, I}(f)$ is a nonzero vector in $Q_\emptyset^V= R(V)$,
     an irreducible $H_q$-module.  As   mentioned earlier in the proof,
     the action of the generators of $H_q$ is expressed
     as $q_\alpha^* q_{\alpha*}-1$, i.e.,  in terms of the bicube structure. 
     Therefore the minimal sub-bicube in $Q^V$ containing $f$, contains
     the entire $Q_\emptyset^V$, and since  
     $v_{\emptyset, I}: Q^V_\emptyset\to Q^V_I$
     is surjective, it contains each $Q_I^V$, so it coincides with $Q^V$. 
     This finishes the proof of Proposition \ref{prop:7-10}. 
     \qed
     
     \begin{rems}
 (a) The components of the   bicube $Q^V$ can be defined directly in
 terms of the $H_q$-module $R(V)$ as ``invariants''  with respect to the Hecke
 algebra of the standard Levi $G_I$. This suggests that   the full $E_q^V$
 can also be defined purely in Hecke-algebra terms.  Since Hecke algebras make
 sense for  more general Coxeter groups, 
 our construction seems to  generalize to such cases as well. 
 
 \vskip .2cm
 
 (b) Let us extend the correspondences $V\mapsto R(V), E_q^V, \Fc_q^V$ to
 arbitrary $G_q$-modules $V$ by additivity with respect to direct sums.
 Then (it is a general property of Hecke algebras) 
 $$
 R(V)\simeq V\otimes_{\k[G_q]}\k[G_q/B_q]\simeq V^{B_q}
 $$
 as $H_q$-modules, where $H_q$ acts on the second factor of the tensor product. 
 
 In particular, $R_{\Fun(F)}=H_q$ as
 a module over itself. Therefore $E_q^{H_q}(\bm)$ consists of functions
 on $O_\bm(\FF_q)$ pulled from $B_q$-invariant functions on $F_{\Hor(\bm)}(\FF_q)$,
  cf. \cite[\S 9]{harish-chandra}. 
 This space   has the same dimension as $E_1(\bm)=\Fun(\bm)$.
 This  makes it natural to  view  
  $\Fc_q^{H_q}=(\Fc_q)^{B_q}$  as a  ``$q$-deformation'' 
 of $\Fc_1$.  
 Note also that  $E_q^{H_q}(\bm)$ can be seen as a decategorified version of the parabolic category $\Oc$, cf. Remark \ref{rems:cat-O}(a) below.
 
\end{rems}
 
\begin{rem} We can go ``downstairs'', i.e., decategorify one more time
  and consider the dimensions of the spaces $\Fun(O_\bm)$. 
 For an algebraic subgroup $H\subset G$ defined over $\ZZ$ 
  consider the  number $n_{G/H}(q) = |(G/H)(\FF_q)|$. 
   Considered as a function of $q$, $n_{G/H}(q)$ is a polynomial. 


 \vskip .2cm
 
 Let $O_\bm \isom G/H$ be a Bruhat orbit.  The dimension  of  $\Fun(O_\bm)$
 is equal to  $n_{O_\bm}(q) = n_{G/H}(q)$. There are two possibilities: 
 \begin{itemize}
 \item[(a)] $H$ is a parabolic, in this case we call $O_\bm$ {\it compact} 
 (the space $O_\bm(\CC)$ is compact).
 
  \item [(b)] $H$ is proper intersection of two  parabolics, i.e., $H=P\cap P'$
  and $H\neq P, H\neq P'$.  In this case we call $O_\bm$ {\it noncompact}.
 \end{itemize}

 \noindent If $H$ is parabolic, say $H = P_I$, then the polynomial $n_{G/H}(q)\in\ZZ[q]$ is prime to  $(q - 1)$ and 
 $q$,
 since
 \[
 n_{G/H}(1) = \text{Card}(W/W_I), \quad n_{G/H}(0) = 1,
 \]
  the second equality following from the Bruhat decomposition.
 
 \vskip .2cm

 Notice that we have a $q$-analogue of Proposition 1.10 (iii): namely,
  $\bm \geq \bn$ is anodyne iff there exists 
 $i\in \ZZ_{\geq 0}$ such that 
 \[
 n_{O_\bm}(q) = q^in_{O_\bn}(q). 
 \]
 This is  true since the fibers of the projection $O_\bm \to O_\bn$  are affine spaces.
 
 \vskip .2cm
 Define a polynomial $\wt n_{O_\bm}(q)$ by 
 \[
 n_{O_\bm}(q) = q^i \, \wt n_{O_\bm}(q), \quad (\wt n_{O_\bm}(q), q)=1. 
 \]
Applying Proposition \ref{prop:shadows},  we see that:
 \begin{itemize}
 
 \item For any $\bm\in \Xi$ 
 the polynomial $n_{O_\bm}(q)$ is not divisible by $q-1$.
 
 \item $O_\bm$ is compact iff $n_{O_\bm}(q)$ is not divisible by $q$.
 
 \item $
 \wt n_{O_\bm}(q) = n_{F_{\Hor(\bm)}}(q) = n_{F_{\Ver(\bm)}}(q)
 $.

 \end{itemize}

     \end{rem}
     

     \section{Further directions and applications}\label{sec:further}
     
     Our approach can be pursued further in several directions. In this final section
     we sketch several such possibilities, leaving the details  for future work. 
     For simplicity we 
     assume that $\k$ is algebraically closed of characteristic $0$.

     \paragraph{Example: braided Hopf algebras.}\label{par:br-hopf}
     Let $\gen=\gen\len_n$. In this case:
     
     \vskip .2cm
     
     $\Delta$ is a root system of type $A_{n-1}$, so  elements of
     $\Dsim=\{\alpha_1,\cdots, \alpha_{n-1}\}$ correspond to the ``intervals'' between
     consecutive integers $\{1,\cdots, n\}$. 
    Subsets $I\subset\Dsim$ correspond to {\em ordered partitions} of $n$, i.e.,
    vectors $\alpha=(\alpha_1, \cdots, \alpha_p)$ of positive integers summing up to $n$.
    The space $F_I = F_\alpha$ consists of flags (filtrations)  of $L=\KK^n$
    \[
    V_\bullet = (V_1\subset \cdots \subset V_p=L), \quad \dim \gr_i^VL = \alpha_i. 
    \] 
    
    Let $\Xi=\Xi_n$ be the 2-sided Coxeter complex for $\gen\len_n$. 
     If $I ,J$ correspond to ordered partitions $\alpha,\beta$ as above, then
    $\Xi_n(I,J)$ is identified with the set of {\em contingency matrices with margins}
    $\alpha$ and $\beta$, i.e.,
    integer matrices $M=\|m_{ij}\|_{i=1,\cdots, p} ^{j=1, \cdots, q}$, $m_{ij}\geq 0$
     with row sums  being $\beta_j$ and column sums being $\alpha_i$. See
     \cite{petersen} \S 6 and \cite{KS-contingency}.

     \vskip .2cm
     
     The relation $M\geq' N$, resp.   $M\geq'' N$, means that $N$ is obtained from $M$ by summing some groups
     of adjacent columns, resp. rows. 
          
     \vskip .2cm
     
The orbit $O_M\subset F_\alpha\times F_\beta$ consists of pairs of filtrations 
$(V_\bullet, V'_\bullet)$ such that
 $\dim \gr_i^V \gr_j^{V'}L  = m_{ij}$, cf. \cite{BLM}.       Note that $ \gr_i^V \gr_j^{V'}L\simeq \gr_j^{V'}\gr_i^V L$
 (Zassenhaus lemma). 
 
 \vskip .2cm
 
 Our methods, specialized to the case   $\gen=\gen\len_n$, 
  lead to a very simple proof and a clear understanding of the main result of \cite{KS-shuffle} 
  (developing a part of \cite{gaitsgory}) on
 braided Hopf algebras.  More precisely:
 
 \vskip .2cm 
 
 Note that Theorem \ref{thm:main} 
can be formulated and proved for perverse sheaves with values in any abelian category
$\Vc$. The concept of a ``sheaf''  can be understood as a   sub-analytic sheaf,
 similarly for complexes as in \cite{KS-shuffle}. The Verdier dual of a constructible
 sub-analytic complex is understood as taking  values in the opposite category $\Vc^\op$. 
 
 \vskip .2cm
 
 Let $(\Vc, \otimes, R, \1)$ be a braided monoidal abelian category with bi-exact $\otimes$,
 and let $A=\bigoplus_{i=0}^\infty A_i$, $A_0=\1$, be a graded bialgebra in $\Vc$ as in \cite{KS-shuffle},
 \S 2.4. To each $n$ we associate a mixed Bruhat sheaf $E=E_n$ on $\Xi_n$ by
 \[
 E_n(M) \,=\,\bigotimes_{i,j} A_{m_{ij}}.  
 \]
 Here the tensor product is  understood in the ``$2$-dimensional'' sense,
 using the interpretation of braided monoidal
 structures as having $N$-fold tensor operations labelled by arrangements of $N$ distinct
 points in the Euclidean plane $\RR^2$. 
We read the matrix structure of $M$  to position each factor $A_{m_{ij}}$ at the point $(-i,j)\in\ZZ^2$
 of a rectangular grid in $\RR^2$. After this, each  map $\del'_{M,N}$ is given by
 the multiplication in $A$, while $\del''_{M,N}$ is given by the comultiplication. 
 
 \vskip .2cm
 
 The space $\hW$ for $\gen\len_n$ is $\Sym^n(\CC)$, the $n$-th symmetric power of $\CC$. 
 Denoting $\Fc_n$ the perverse sheaf on $\Sym^n(\CC)$ corresponding to $E_n$
 by Theorem \ref{thm:main}, we get a system $(\Fc_n)_{n\geq 0}$ of perverse sheaves that is manifestly
  factorizable
 \cite[Def. 3.2.5]{KS-shuffle}  and the main result of \cite{KS-shuffle} (Theorem 3.3.1)
 follows easily. 
 
 \vskip .2cm
 
 Note that the bicube associated to $E_n$ consists of ``$1$-dimensional'' 
 (linearly ordered) tensor products
 \[
 A_{\alpha_1}\otimes \cdots\otimes A_{\alpha_p}, \quad \alpha_1+\cdots + \alpha_p=n, 
 \]
with the $u$-maps given by multiplication  (bar-construction) and the $v$-maps given by the comultiplication (cobar-construction). 

\vskip .2cm

 Passing from a bicube to a mixed Bruhat sheaf in this and other examples can be
 seen as ``unfolding''  of a naive  $1$-dimensional structure to a more fundamental $2$-dimensional one.

\paragraph{Eisenstein series and constant terms.} An example of a braided Hopf algebra
is given by $\Hc(\Ac)$, the Hall algebra of a hereditary abelian category $\Ac$
with appropriate finiteness conditions \cite{green}.
One can apply this approach to  $\Ac=\Coh(X)$, the category of coherent sheaves on
a smooth projective
curve $X/\FF_q$ as in \cite{k-eis-series, k-schiffman-vass}.
 Considering  functions supported on vector bundles,  we get
a graded, braided Hopf algebra $\Hc^\Bun = \bigoplus_{n\geq 0} \Hc^\Bun_n$ where $\Hc_n^\Bun$
  consists of  {\em   unramified automorphic forms} for the group
 $GL_n$ over the function field $\FF_q(X)$. The multiplication is given by forming
(pseudo) Eisenstein series and comultiplication  by taking the constant
term of an automorphic form. 
Because $\Coh(X)$ does not fully satisfy the finiteness conditions  (an object may have
infinitely many subobjects, but only finitely many subobjects of any given degree),
the comultiplication in $\Hc^\Bun$ must be understood using generating functions
or rational functions of a spectral parameter as in \cite {k-schiffman-vass}. 
With this taken into account
(i.e., after extending the field of scalars to allow  the dependence on the extra parameter),
 $\Hc^\Bun$ gives, for each $n$, a mixed Bruhat sheaf
on $\Xi_n$ and so a perverse sheaf on $\Sym^n(\CC)$, as explained in  \S \ref{par:br-hopf}

\vskip .2cm

If now $G$ is a general split reductive group over $\ZZ$ with Lie algebra $\gen$,
 we still have the   classical theory of unramified automorphic forms and Eisenstein series
 for $G$ over  $\FF_q(X)$ in \cite{moeglin}. It usually  appears in the form of a bicube
 $Q$, where
 \[
 Q_I \,=\,\Fun(\Bun_{G_I}(X)), \quad I\subset\Dsim
 \]
 is the space of automorphic forms for  the standard Levi $G_I$, i.e., of functions on
the set of isomorphic classes of principal $G_I$-bundles on $X$. For $I\subset J$
the map $v_{IJ}: Q_I\to Q_J$ is given by taking the (pseudo) Eisenstein series and
$u_{IJ}: Q_J\to Q_J$ is given by taking the constant term of an automorphic form. 
 
\vskip .2cm
 
 For general $G$, this theory does not have a  Hopf algebra interpretation. 
 However, one can extend the above bicube $Q$ to a mixed Bruhat sheaf $E$ and so
 obtain a perverse sheaf on $\hW$. For this, given $\bm\in\Xi(I,J)$, 
  one should consider the moduli space
 $\Bun_{G,\bm}(X)$ formed by principal $G$-bundles together with a 
 $P_I$-structure and a $P_J$-structure (i.e., sections of the associated bundles
 with fibers $G/P_I$ and $G/P_J$), everywhere in relative position $\bm$. 
 The corresponding $E(\bm)$ is then found inside the space of functions on 
   $\Bun_{G,\bm}(X)$.

    \paragraph{Categorical upgrade: mixed Bruhat schobers.}
    The concept of a mixed Bruhat sheaf is very convenient for a categorical upgrade, 
          i.e., 
    replacing   vector spaces   with $\k$-linear  dg-enhanced  triangulated categories
    (simply  ``triangulated categories'' below).
    The possibility of such upgrade of the theory of perverse sheaves was
    raised in  \cite{KS-schobers}, where such hypothetical objects were called
    perverse schobers.  See also \cite{bondal-KS, donovan}. 
    
    \vskip .2cm
 
        Recall \cite {AL} that a diagram of   triangulated categories
      \[
     \xymatrix{
\Cen \,    \ar@<.5ex>[r]^{f}&\,  \Den
\ar@<.5ex>[l]^{g}
    }
\]
consisting of a (dg-) functor $f$ 
and its right adjoint $g=f^*$ is called a {\em spherical adjunction}
(and $f$ is called a {\em spherical functor}), if the cones of the unit and counit of the adjunction
\[
T_\Cen = \Cone\bigl\{e:  \Id_\Cen \lra  gf\bigr\} [-1], \quad T_\Den = \Cone\bigl\{ \eta: fg\lra \Id_\Den\bigr\}
\]
     are equivalences  (i.e., quasi-equivalences of dg-categories). As noticed in \cite{KS-schobers}, such a diagram
     can be seen as a categorical upgrade of a perverse sheaf $\Fc\in\Perv(\CC,0)$
     given by $(\Phi, \Psi)$-description, see Example \ref{ex:sl2-1}.
     
     \vskip .2cm 
     
      Now, $\Perv(\CC,0)=\Perv(\hW)$
     for $\gen=\sen\len_2$, so Theorem  \ref{thm:main} suggests a generalization
     of the concept of a spherical functor to arbitrary $\gen$. Let us call such structures
     {\em mixed Bruhat schobers} and sketch the main features of the definition.
     
     \vskip .2cm
     
     So a mixed Bruhat schober $\Een$ should consist of 
  triangulated categories $\Een(\bm)$, $\bm\in\Xi$ and dg-functors 
     \[
\begin{gathered}
\den'_{\bm,\bn} = \den'_{\bm,\bn,\Een}: \Een(\bm)\lra \Een(\bn), \quad \bm\geq' \bn,
\\
\den''_{\bm,\bn} = \den''_{\bm,\bn, \Een} : \Een(\bn)\lra \Een(\bm), \quad \bm\geq'' \bn,
\end{gathered}
\]
satisfying the following analogs of (MBS1-3). First, 
 so that the $\den'$, as well as the $\den''$ must be transitive (up to coherent homotopies).
 Second,   (MBS2) is upgraded into the data  of a ``filtration'' on the functor 
 $\den''_{\bn, \bn'} \den'_{\bm', \bn'}$, $\bm'\geq' \bn' \leq''\bn$
 with ``quotients'' being  the functors
 $ \den'_{\bm,\bn,} \den''_{\bm, \bm'}$ for $\bm$ running in the poset $(\sup(\bm',\bn), \teq)$. 
 Such a filtration can be understood either as a Postnikov system (see
 \cite  {KS-hyp-arr} \S 1A) or as a Waldhausen diagram
 (see \cite[\S 5]{DK-TSTC}  or \cite[\S 7.3]{DK-HSS}), adapted for the case of
 a partially-ordered indexing set. 
 The analog of the condition (MBS3) is that $\den'_{\bm,\bn}$ for any anodyne
 $\bm\geq'\bn$ and  $\den''_{\bm,\bn}$ for any anodyne $\bm\geq''\bn$ must be an equivalence
 (i.e., a quasi-equivalence of dg-categories). 
 Further, we should impose natural adjointness conditions meaning that
 $\den'_{\bm,\bn}$ is identified with the right adjoint of $\den''_{\bm^\tau,\bn^\tau}$
 after composing with appropriate ``homotopies" connecting $\bm$ with $\bm^\tau$
 and $\bn$ with $\bn^\tau$ (note that the cells $U_\bm$ and $U_{\bm^\tau}$ always
 lie in the same stratum of $\Sc^{(0)}$, and so $\bm$ and $\bm^\tau$ can be connected
 by a chain of anodyne   $\leq', \leq''$ or their inverses). 
 
 Precise details will be given in a subsequent paper. Let us list two natural sources of
 such structures.

 \paragraph{Constructible sheaves on Bruhat orbits.}  
We can upgrade the constructions of  \S \ref{sec:Fq}  by replacing the space of
 functions on $\FF_q$-points of a variety with the category of constructible
 complexes.  
 
 We consider the simplest setting when $\KK=\CC$. For an algebraic variety $X/\CC$
 let $D(X)=D^b_\constr(X)$ be the derived category of bounded complexes with
 cohomology sheaves constructible with respect to some $\CC$-algebraic stratification as in 
 Appendix \ref{app:strat}.   
 \vskip .2cm
 
 We have a $\Dsim$-bicube  $\Qen$ of triangulated categories similar to that \S \ref{sec:Fq}\ref{par:F_q-bic}
 It consists of the categories $D(F_I)$ and functors 
 \[
(q_{IJ})_* = (q_{IJ})_!: D(F_I) \lra D(F_J), \,\,\, (q_{IJ})^*: D(F_J) \to D(F_I), \quad I\subset J. 
 \]
 Note that for $\gen=\sen\len_2$ the bicube reduces to the diagram
 \[
 \xymatrix{
D(\CC\PP^1) \,    \ar@<.7ex>[r]^{\pi_*}&\,  D(\pt)  
\ar@<.7ex>[l]^{\pi^*}
    }, \quad \pi: \CC\PP^1\to\pt,
 \]
 which is a proto-typical example of a spherical adjunction, $\CC\PP^1$ being the sphere $S^2$,
 cf. \cite[Ex. 1.10]{KS-schobers}. 
  To extend the bicube $\Qen$,   we proceed
  similarly to  \S \ref{sec:Fq}. 
 \vskip .2cm

 Given $\bm\in\Xi$, we consider first the category $  D(O_\bm)$. For $\bm\geq'\bn$
 we define $\den'_{\bm,\bn}: D(O_\bm) \to D(O_\bm)$ to be the
 functor $(p_{\bm,\bn})_!$, the (derived) direct image with proper supports. 
 For  $\bm\geq''\bn$ we define
 $\den''_{\bm,\bn}: D(O_\bn)\to D(O_\bm)$ to be the functor $(p_{\bm,\bn})^*$. 
 We define the category $\Een(\bm)$ 
 to be the  essential image of the
 pullback functor $(r'_\bm)^*: D(F_{\Hor(\bm)}) \to D(O_\bm) $. 
 As in \S  \ref{sec:Fq}, we see that $\den', \den''$ preserve the $\Een(\bm)$,
 so we have a diagram of triangulated categories
 \[
 \Een \,=\, (\Een(\bm), \den',\den'')
 \]
 upgrading the mixed Bruhat sheaf $E_q$ of Theorem \ref{thm:Fq}.

 \begin{rems}\label{rems:cat-O}
 
 (a)  We can also  consider the diagram $\Een^B$ formed by $B$-equivariant objects in the 
 $\Een(\bm)$. Then $\Een^B(\bm)$ is identified with the category of $B$-equivariant
 constructible complexes of $F_{\Hor(\bm)}$, so  the diagram  consists of
   various (graded  derived versions of) parabolic 
 categories $\Oc$.
 
 \vskip .2cm

(b)  Instead of $D(O_\bm)$, we can use other types of  ``categories of sheaves''  on  $O_\bm$
 which possess an appropriate formalism of pullbacks and pushforwards. For example, we
 can use
  the category of mixed motives over $O_\bm$ in \cite{cisinski}.
  
  \vskip .2cm

  (c) We can also  take the ``quasi-classical'' approach, i.e., consider, instead of
  constructible complexes (i.e., complexes of holonomic regular  D-modules) on the $O_\bm$, complexes of 
  coherent sheaves on $T^*(O_\bm)$, thus establishing a connection with
  the braid group actions on the coherent derived categories of such cotangent bundles
  via flops \cite{bezr-riche, cautis, khovanov-thomas}. 
 
 \end{rems}

 \paragraph{Parabolic induction and restriction.} 
  We consider the simplest case of finite Chevalley groups.
  That is, 
   take $\KK=\FF_q$.   For any $I\subset\Dsim$ let 
$\Cen_I $ be the derived category of finite-dimensional $\k$-linear
representations of the finite group
$G_I(\FF_q)$. 
 If $I\subset J$, then $G_I\subset G_J$ and 
 we have the classical
{\em parabolic induction} and {\em restriction} functors
\[
\begin{gathered}
\Ind_{I,J}: \Cen_I  \lra  \Cen_J, \quad M\mapsto \Ind_{(P_I\cap G_J)(\FF_q)}^{G_{J}(\FF_q)} \, M, 
\\
\Res_{I,J}: \Cen_J\lra \Cen_I, \quad N\mapsto N^{(U_I\cap G_J)(\FF_q)},
\end{gathered}
\]
  These two functors are both left and right adjoint to each other. 
 Further, the $\Ind$ and $\Res$-functors are transitive, so we have a bicube $\Cen$ of 
 triangulated categories.
 The passing to the derived categories, seemingly unnecessary for this simple case
 (all the functors are exact at the level of abelian categories) makes the following example
 more salient.
 
 \begin{ex}
 Let $\gen=\sen\len_2$. Then the bicube has the form
 \[
  \xymatrix{
\Cen_\emptyset = D^b\Rep(\FF_q^*)  \,    \ar@<.7ex>[r]^{\hskip -0.5cm \Ind} &\,  D^b\Rep(SL_2(\FF_q)) = \Cen_{\{1\}}
\ar@<.7ex>[l]^{\hskip -0.5cm \Res},
    }
 \]
 where $\Ind$ is the functor of forming the principal series representation and $\Res$
 is the functor of invariants with respect to the standard unipotent subgroup. 
 It follows from the elementary theory of representations of $SL_2(\FF_q)$, cf.  
 \cite{fulton-harris}, 
 that this is in fact a spherical adjunction.  
   \end{ex}
  
  To extend the bicube $\Cen$ 
  to more general  Bruhat  orbits, we  associate to each $\bm\in\Xi$
  the derived category of $G(\FF_q)$-equivariant vector bundles $V$ on
  the discrete set  $O_\bm(\FF_q)$  such   that for any element  $x\in O_\bm(\FF_q)$,  
  the unipotent radical of the stabilizer of $x$ acts in the fiber $V_x$ trivially. 
  Such category is equivalent to the derived category of representations of the standard Levi
  $G_{\Hor(\bm)}(\FF_q)$.      
 

\appendix

\section{Stratifications and constructible sheaves}\label{app:strat}

\renewcommand{\theparagraph}{\Alph{section} \arabic{paragraph}.} 

\paragraph{Stratifications.} 
 We fix some terminology to be used in the rest of the paper. 
By a {\em space} we mean a 
real analytic space.  For a space $X$ we can speak about subanalytic subsets in $X$,
see, e.g., \cite[\S 8.2]{KaSha} and references therein. Subanalytic subsets
form a Boolean algebra. 

\begin{Defi}\label{def:lcd}
Let $X$ be a space.

 \begin{enumerate}

\item[(a)]   A {\em  partition}  of $X$ is a  finite family $\Sc=(X_a)_{a\in A}$ of 
subanalytic subsets in $X$ such that we have a disjoint
decomposition $X=\bigsqcup_{a\in A} X_a$.
of $X$ as a disjoint union of subanalytic sets. The sets $X_a$ are called
the {\em strata} of the  partition $\Sc$.

\item[(b)] A {\em locally closed decomposition (l.c.d.)} of $X$ is a partition
  $\Sc=(X_a)_{a\in A}$ such that each  $X_a$ is locally closed and the closure of each
  $X_a$ is a union
of strata.  In this case  the set $A$ becomes
partially ordered by $a\leq b$ if $X_a\subset \ol X_b$.

\item[(c)] A {\em stratification} of $X$ is an l.c.d. such that each $X_a$ is an analytic submanifold and
the Whitney conditions are satisfied. 
A {\em stratified space} is a real analytic space with a stratification. 
 \end{enumerate}
\end{Defi}

  Further a {\em cell decomposition} of $X$ is a stratification $\Sc$ such that each stratum
is homeomorphic to an open $d$-ball $B^d$ for some $d$. A cell
decomposition is called {\em regular}, if for each cell (stratum) $X_a$
there exists a homeomorphism $B^d\to X_a$ which extends to an embedding
of the closed ball $\ol B^d\to X$ whose image is a union of cells. We will say that
$(X,\Sc)$ is a {\em regular cellular space}.

\vskip .2cm

A cell decomposition of $X$ if called {\em quasi-regular}, if $X$, as a stratified space,
can be represented as $Y\- Z$ where $Y$ is a regular cellular space and $Z\subset  Y$
is a closed cellular subspace.

\vskip .2cm

 Given two partitions $\Sc$ and $\Tc$ of $X$, we say that $\Sc$ {\em refines}
 $\Tc$ and write $\Sc \prec \Tc$, if each stratum of $\Tc$ is a union of strata of $\Sc$.

\vskip .2cm

Given two partitions $\Sc =(X_a)$ and $\Tc= (Y_b)$ of $X$, their 
{\em maximal common refinement} $\Sc\wedge \Tc$ is the  partition
into  subsets defined as connected components of the
 $X_a\cap Y_b$. If $\Sc, \Tc$ are l.c.d.'s or stratifications, then so is
  $\Sc\wedge\Tc\prec \Sc, \Tc$.

\vskip .2cm

We also have the  {\em maximal common coarsening}  
$\Sc\vee \Tc$. This is a partition 
consisting of equivalence classes of the equivalence
relation $\equiv$ on $X$ defined as follows.
We first form the relation $R$
defined by:  $xRy$ if $x$ and $y$ lie in the same stratum $X_a$ of $\Sc$ or in the same
stratum $Y_b$ of $\Tc$. 
Then we define $\equiv$ as 
 the equivalence closure of $R$. Thus the strata of $\Sc\vee\Tc$
 are certain unions of the $X_a\cap Y_b$.
 We do not know whether the assumption that 
 $\Sc, \Tc$ are  both  
  l.c.d.'s (resp. stratifications), implies that  $\Sc\vee\Tc$ 
  is an l.c.d (resp.  a stratification).

 \vskip .2cm
 
  We will be particularly interested in the cases when $\Sc, \Tc$ and $\Sc\vee\Tc$ are
  all stratifications.

\paragraph{Constructible sheaves.} Let $\k$ be a field and $(X,\Sc)$ be a stratified
space. As usual, a sheaf $\Gc$ of $\k$-vector spaces on $X$ is called
$\Sc$-{\em constructible} if the restriction of $\Gc$ on each stratum
is  locally constant of finite rank.   For $V\in \Vect_\k$ we denote
$\ul V_X$  the constant sheaf on $X$ with stalk $V$. 

\vskip .2cm

We denote by $\Sh(X,\Sc)$ the
category of $\Sc$-constructible sheaves. 
A complex $\Fc$ of sheaves is called (cohomologically) $\Sc$-constructible, if each
cohomology sheaf $\ul H^q(\Fc)$ is $\Sc$-constructible. We denote
$D^b_\Sc\Sh(X)$ the derived category of  $\Sc$-constructible complexes
with only finitely many nonzero cohomology sheaves. 
It carries the Verdier duality $\DD$, see \cite{KaSha}. 
The following is clear.

\begin{prop}\label{prop:strat-append}
\begin{enumerate}

\item[(a)] Let $\Sc,\Tc$ be two stratifications of $X$ such that 
 $\Sc\prec\Tc$. Then each $\Tc$-constructible
sheaf is $\Sc$-constructible.

\item[(b)] Let $\Sc, \Tc, \Uc$ be three stratifications of $X$ such that $\Uc=\Sc\vee\Tc$.
Suppose a sheaf $\Gc$ is both $\Sc$-constructible and $\Tc$-constructible.
Then $\Gc$ is $\Uc$-constructible. \qed
\end{enumerate}
\end{prop}

 Let $\Sc=(X_a)_{a\in A}$ be a quasi-regular cell decomposition of $X$,
so $(A,\leq )$ is naturally a poset under the closure relation on cells. 
 Recall,  
  cf. \cite {KS-hyp-arr} \S 1D,  that an $\Sc$-constructible
(cellular) sheaf $\Gc$ on $X$ is uniquely determined by the data of its
{\em stalks} $G_a = \Gamma(X_a, \Gc)$ at the cells  and {\em generalization maps}
$\gamma_{a,b}: G_a\to G_b$, $a\leq b$, which satisfy the transitivity conditions
\be\label{eq:G-transit}
\gamma_{a,a}=\Id, \quad \gamma_{a,c} = \gamma_{b,c}\circ\gamma_{a,b}, \,\, a\leq b\leq c.
\ee
A datum $R=(G_a, \gamma_{ab})$  formed by  finite-dimensional vector spaces
$G_a$ and linear maps $\gamma_{ab}$ satisfying  \eqref {eq:G-transit}, will be called
a {\em representation} of $A$. It is simply  a covariant functor from $(A,\leq)$ 
(considered as a category)
to $\Vect_\k$. Representations of $A$ form an abelian category $\Rep(A)$; we denote by 
$D^b\Rep(A)$ the corresponding bounded derived category. 
 Given $R = (G_a, \gamma_{ab})\in\Rep(A)$,
one defines directly the cellular sheaf $\sh(R)$ with stalks $G_a$ and generalization maps
$\gamma_{ab}$, thus giving a functor $\sh: \Rep(A)\to \Sh(X,\Sc)$. 
The above discussion can be formulated more precisely as follows,
see, e.g., \cite[Prop. 1.8]{KS-hyp-arr} :

\begin{prop}\label{prop:cell-sheaves}
Let $\Sc=(X_a)_{a\in A}$ be a quasi-regular cell decomposition of $X$. Then:
\begin{enumerate}

\item[(a)] The functor $\sh: \Rep(A)\to \Sh(X,\Sc)$ is an equivalence of abelian categories. 

\item[(b)] The  termwise extension of $\sh$ to complexes 
 defines  an equivalence
of triangulated categories
$D\sh: D^b(\Rep(A))\to D^b_\Sc\Sh(X)$.  \qed
 \end{enumerate}
\end{prop}


\vskip 1cm

\small{

M.K.: Kavli IPMU, 5-1-5 Kashiwanoha, Kashiwa, Chiba, 277-8583 Japan, 
\hfil\break
{\tt mikhail.kapranov@protonmail.com}

\smallskip

 V.S.: Institut de Math\'ematiques de Toulouse, Universit\'e Paul Sabatier, 118 route de Narbonne, 
31062 Toulouse, France, 
 {\tt schechtman@math.ups-tlse.fr }
 }


\addcontentsline{toc}{section}{References}
\begin{thebibliography}{100}

\bibitem{AL} R. Anno, T. Logvinenko. Spherical DG-functors.
 {\em  J. Eur. Math. Soc.} {\bf 19}  (2017) 2577-2656.  
 
 \bibitem{basak} T. Basak. Combinatorial cell complexes and Poincar\'e duality. 
 {\em Geom. Dedicata} {\bf 147} (2010) 357-387. 


\bibitem{BBD} A. Beilinson, J. Bernstein, P. Deligne, O. Gabber.
Faisceaux pervers. {\em Ast\'erisque} {\bf 100} 1980. 

\bibitem{BLM}  A. Beilinson, G. Lusztig, R. MacPherson. A geometric setting for
the quantum deformation of $GL_n$. {\em Duke Math. J.} {\bf 61} (1990) 655-677. 

\bibitem{bezr-tilting} A. Beilinson, R. Bezrukavnikov, I. Mirkovic. Tilting exercises. 
{\em Mosc. Math. J.} {\bf  4}  (2004)  547--557. 

\bibitem{beil-volog} A. Beilinson, V.  Vologodsky. A DG guide to Voevosky motives. arXiv:math/0604004. 

\bibitem{bernstein-bezr}   J. Bernstein, R.  Bezrukavnikov, D.  Kazhdan. 
Deligne-Lusztig duality and wonderful compactification. arXiv:1701.07329. 

\bibitem{bernstein-zelevinsky} I. N. Bernstein, A. V. Zelevinsky. Induced representations of reductive
p-adic groups I. {\em Ann, Sci. \'Ecole Norm. Sup.} {\bf 10} (1977) 441-472. 

\bibitem {BFS} R. Bezrukavnikov, M. Finkelberg, V. Schechtman.  Factorizable Sheaves and Quantum Groups, 
\newblock {\em  Lecture Notes in Math.} {\bf 1691}, 
Springer-Verlag, 1998. 

\bibitem{bezr-riche} R. Bezrukavnikov, S.  Riche. Affine braid group actions on derived categories of Springer resolutions. {\em Ann. Sci. 
\'Ec. Norm. Sup.} 
{\bf  45}   (2012) 535-599. 



\bibitem{BZ} A.~Bj\"orner, G.~Ziegler.
\newblock Combinatorial stratification of complex arrangements.
\newblock {\em Jour. AMS}, {\bf 5} (1992), 105-149. 


\bibitem{bondal-KS} A. Bondal, M. Kapranov, V. Schechtman. Perverse schobers
and birational geometry. {\em Selecta Math.} {\bf 24} (2018) 85-143. 

\bibitem{borel} A. Borel. Linear Algebraic Groups. Springer-Verlag, 1991.

\bibitem{bourbaki} N. Bourbaki. Lie Groups and Lie Algebras, Chapters 4-6. 
Springer-Verlag, 2008. 

\bibitem{cautis} S. Cautis. Flops and about: a guide. arXiv:1111.0688. 

\bibitem{cisinski} D.-C. Cisinski, F. D\'eglise. Triangulated categories of mixed motives.
arXiv:0912.2110. 

\bibitem{deconcini}
C. De Concini, C. Procesi.  Complete symmetric varieties. In: ``Invariant theory'' (Montecatini,
1982) p. 1-44, {\em Lecture Notes in Math.} {\bf  996} Springer, Berlin, 1983.

\bibitem{diaconis} P. Diaconis. A. Gangolli. Rectangular arrays with fixed margins, in: 
``Discrete Probability and Algorithms" (Minneapolis MN 1993) p. 15-41. IMA Vol. Math. Appl.
{\bf 72}, Springer-Verlag, 1995. 

\bibitem{donovan} W. Donovan. Perverse schobers on Riemann surfaces:
constructions and examples. {\em Eur. J. Math.} {\bf 5} (2019) 771-797. 

\bibitem{DK-TSTC} T. Dyckerhoff, M. Kapranov. Triangulated surfaces in
triangulated categories. {\em J. Eur. Math. Soc.} {\bf 20} (2018) 1473-1524.

\bibitem{DK-HSS} T. Dyckerhoff, M. Kapranov. Higher Segal Spaces.  
{\em Lecture Notes in Math.}  {\bf 2244}.  Springer-Verlag, 2019. 

\bibitem{westerland} J. S. Ellenberg, T. T. Tran, C. Westerland. Fox-Neuwirth-Fuks cells, quantum shuffle algebras
and Malle's conjecture for functional fields. arXiv:1701.04541. 


\bibitem{fox-neuwirth} R. Fox, L. Neuwirth. Braid groups. {\em Math. Scand.} {\bf 10} (1962) 119-126. 


\bibitem{fuchs} D.B. Fuks. Cohomology of the braid group mod 2. 
\newblock {\em Funkc. Anal. i Pril.} {\bf 4} (1970), N. 2,  62-73.

\bibitem{fulton-harris} W. Fulton, J. Harris. Representation theory. A First Course.
Springer-Verlag,  1991. 

 

 \bibitem{gaitsgory} D. Gaitsgory. Notes on factorizable sheaves.  Preprint (2008) available at
  <http://www.math.harvard.edu/~gaitsgde/GL/FS.pdf>  

 

\bibitem {GM-cusp} M.~Granger, Ph.~Maisonobe. 
\newblock Faisceaux pervers relativement \`a un point de rebroussement.
{\em C.R. Acad. Sci. Paris. S\'er. I} {\bf 299} (1984), 567-570. 

\bibitem{green} J. A. Green. Hall algebras, hereditary algebras and quantum groups.
{\em Invent. Math. } {\bf 120} (1995) 361-377. 

\bibitem{harish-chandra} Harish-Chandra. Eisenstein series over finite fields. 
In: ``Functional Analysis and Related Fields''
(E. F. Browder  ed.) p. 76-88,  Springer-Verlag, 1970. 

 

\bibitem{iwahori} N. Iwahori. On the structure of a Hecke ring of a Chevalley group
over a finite field. {\em J. Fac. Sci. Univ. Tokyo}, {\bf 10} (1964) 215-236. 

\bibitem{joyal-street} A. Joyal, R. Street. The category of representations of the general
linear groups over a finite field. {\em J. Algebra}, {\bf 176} (1995) 908-946. 

\bibitem{k-eis-series} M. Kapranov. Eisenstein series and quantum affine algebras.
{\em J. Math. Sci. (N.Y.)} {\bf 84} (1997) 1311-1360. 


 \bibitem{KS-schobers}  M. Kapranov, V. Schechtman. Perverse  schobers. 
 arXiv:1411.2772. 
 
   \bibitem  {KS-hyp-arr} M. Kapranov, V. Schechtman. Perverse sheaves on real hyperplane arrangements. {\em Ann. Math.} {\bf 183} (2016), 619 - 679. 
 
 \bibitem {KS-shuffle} M. Kapranov, V. Schechtman. Shuffle algebras and perverse sheaves.  {\em Pure and Appl. Math. Quarterly} {\bf 16} (2020) 573-657. 
 
 
 \bibitem{KS-contingency} M. Kapranov, V. Schechtman (with an appendix by P. Etingof). 
 Contingency tables with
 variable margins.  {\em  SIGMA} {\bf 16} (2020), \#062. 
 
    \bibitem  {KS-hyp-arr-II} M. Kapranov, V. Schechtman. Perverse sheaves on real hyperplane arrangements II. arXiv:1910.01677.
    
    \bibitem  {KS-prob} M. Kapranov, V. Schechtman. PROBs and perverse sheaves I. Symmetric products. arXiv:2102.13321.
    
    \bibitem{k-schiffman-vass} M. Kapranov, O. Schiffmann, E. Vasserot.
    The Hall algebra of a curve. {\em Selecta Math.} {\bf 23} (2017) 117-177. 
 
 \bibitem{KaSha} M. Kashiwara, P. Schapira. Sheaves on Manifolds. Springer-Verlag, 1991. 
 
 \bibitem{khovanov-thomas} M.  Khovanov,  R.  Thomas.
 Braid cobordisms, triangulated categories, and flag varieties.
 {\em Homology Homotopy Appl.} 
    {\bf  9} (2007) 19-94.




 
 \bibitem{moeglin} C. Moeglin, J.-L. Waldspurger. Spectral Decomposition and Eisenstein
 Series. Cambridge Univ. Press, 1995. 
 


 \bibitem{petersen} T. K. Petersen. A two-sided analog of the Coxeter complex.  
 {\em Electronic Journal of
Combinatorics} {\bf 25} (2018) \#P4.64. 
 
 \bibitem{soibelman} Y.  Soibelman. Meromorphic tensor categories. arXiv:q-alg/9709030.

 \bibitem{tits} J. Tits, Buildings of spherical type and finite BN pairs.
 {\em Lecture Notes in Math.} 
    {\bf  386}, Springer-Verlag, 1974. 
  
\end{thebibliography}
\end{document}